\documentclass[12pt,leqno]{article}
\usepackage{amsmath,amssymb,amscd}
\newcommand\qed{{\unskip\nobreak\hfil\penalty50\hskip2em\vadjust{}
    \nobreak\hfil$\Box$\parfillskip=0pt\finalhyphendemerits=0\par}}

\numberwithin{equation}{section}

\newtheorem{lemma}{Lemma}[section]
\newtheorem{theorem}{Theorem}[section]
\newtheorem{corollary}{Corollary}[section]
\newtheorem{proposition}{Proposition}[section]
\newcommand{\bl}{\begin{lemma}}
\newcommand{\el}{\end{lemma}}
\newcommand{\bp}{\begin{proposition}}
\newcommand{\ep}{\end{proposition}}
\newcommand{\bth}{\begin{theorem}}
\newcommand{\et}{\end{theorem}}
\newcommand{\bco}{\begin{corollary}}
\newcommand{\ec}{\end{corollary}}
\newcommand{\be}{\begin{equation}}
\newcommand{\ee}{\end{equation}}
\newcommand{\bal}{\begin{align}}
\newcommand{\eal}{\end{align}}
\newcommand{\bi}{\begin{itemize}}
\newcommand{\ei}{\end{itemize}}
\newcommand{\la}{\label}

\newcommand{\bs}{{\bigskip}}
\newcommand{\ms}{{\medskip}}
\newcommand{\noi}{\noindent}
\newcommand{\1}{{1\!\!1}}
\newcommand{\bz}{{\bold z}}

\newcommand{\bq}{{\bold q}}

\newcommand{\bb}{{\bold b}}

\newcommand{\bw}{{\bold w}}

\renewcommand{\b}{\beta}
\renewcommand{\d}{\delta}
\newcommand{\D}{\Delta}
\newcommand{\e}{\varepsilon}
\newcommand{\g}{\gamma}
\newcommand{\G}{\Gamma}
\renewcommand{\l}{\lambda}
\renewcommand{\L}{\Lambda}

\newcommand{\var}{\varphi}
\newcommand{\s}{\sigma}

\renewcommand{\th}{\theta}
\renewcommand{\O}{\Omega}
\renewcommand{\o}{\omega}
\newcommand{\z}{\zeta}

\newcommand{\x}{\times}
  
\renewcommand{\i}{\infty}
\newcommand{\p}{\partial}

\newcommand{\bE}{{\mathbb E}}
\newcommand{\bN}{{\mathbb N}}
\newcommand{\bP}{{\mathbb P}}
\newcommand{\bR}{{\mathbb R}}

\newcommand{\bW}{{\mathbb W}}

\newcommand{\bn}{{\bf{n}}}
\newcommand{\by}{{\bf{y}}}

\newcommand{\cK}{{\mathcal K}}
\newcommand{\cD}{{\mathcal D}}

\newcommand{\cL}{{\mathcal L}}

\newcommand{{\cA}}{{\mathcal A}}
\newcommand{\cM}{{\mathcal M}}

\newcommand{\cB}{{\mathcal B}}

\newcommand{\cQ}{{\mathcal Q}}

\begin{document}
\title{\Large\bf{Kinetic description of scalar conservation laws with Markovian data}}

\author{FRAYDOUN REZAKHANLOU
\\
Department of Mathematics\\UC Berkeley}

\maketitle

\begin{abstract}
  We derive a kinetic equation to describe the 
  statistical structure of solutions $\rho$ to scalar conservation
  laws $\rho_t=H(x,t,\rho )_x$, 
  with certain Markov initial conditions. When the Hamiltonian function 
is convex and increasing  in $\rho$, we show that 
  the solution $\rho(x,t)$ is a Markov process in $x$ (respectively $t$)
with $t$ (respectively $x$) fixed. Two classes of Markov conditions are considered in this article. In the first class, the initial data is characterize by a drift $b$ which satisfies a linear PDE, 
and a jump density $f$ which
 satisfies a kinetic equation as time varies. 
In the second class, the initial data is a concatenation of fundamental
solutions that are characterized by a parameter $y$, which is a Markov jump process with a jump density $g$ satisfying a kinetic equation. 
When $H$ is not 
increasing in $\rho$, 
the restriction of $\rho$ to a line in $(x,t)$ plane 
is a Markov process of the same type, provided that 
the slope of the line satisfies an inequality. 
 \end{abstract}

\section{Introduction}
\label{sec1}

Hamilton–Jacobi equation (HJE) is one of the most popular and
studied PDE which enjoys vast applications in numerous areas
of science. Originally HJEs were formulated in connection with
the completely integrable Hamiltonian ODEs of celestial mechanics.
They have also been used to study the evolution of the value functions in control and differential game theory.
 Several growth models in physics and biology are described by
HJEs. In these models, a random interface
separates regions associated with different phases and the interface can be locally
approximated by the graph of a solution to a HJE. To make up for the lack of exact information or/and the presence of impurity, it is common to assume that the Hamiltonian function which appears in our HJE is random. 
Naturally we would like to understand how the randomness affects
the solutions and how the statistics of solutions are propagated with time.

In dimension one, the differentiated version of a Hamilton–Jacobi equation becomes
a scalar conservation law for the inclination of the one-dimensional interface, and may be used to model an one-dimensional fluid. In the context of fluids, we wish to obtain some qualitative information about
 the structure of shocks and their fluctuations.

The primary purpose of this article is to derive an evolution equation for the statistics of solutions to a HJE in dimension one. We achieve this by utilizing a kinetic description for the shock densities of piecewise smooth solutions.

Given a $C^2$ Hamiltonian function $H:\bR\x[0,\i)\x\bR\to\bR$, we consider the HJE 
\be\la{eq1.1}
u_t=H(x,t,u_x),\ \ \ \ t\ge t_0,
\ee
or the corresponding 
scalar conservation law
\begin{equation}
  \label{eq1.2}
      \rho_t = H(x,t,\rho )_x ,\ \ \ \ t\ge t_0.
               \end{equation}
We assume that the Hamiltonian function $H(x,t,\rho)$ is convex in 
the {\em momentum} variable $\rho$. 
As our main goal, we show that the statistics of $\rho(x,t)$
admits an exact kinetic description when the initial data
$\rho^0(x) = \rho(x,t_0)$ is an  inhomogeneous  Markov process.

\subsection{Main result I}

For our first result, we assume that  the initial data
$\rho^0 = \rho^0(x)$ is a piecewise-deterministic  inhomogeneous  Markov process (PDMP) Markov process 
determined by a generator $\mathcal{A}_{x}^0=\cA_{x,t_0}$ acting on test
functions $\psi(\rho)$ according to
\begin{equation}
  \label{eq1.3}
  (\mathcal{A}_x^0 \psi)(\rho) = b^0(x,\rho) \psi'(\rho) 
  + \int_{\rho}^\infty \big(\psi(\rho_*) - \psi(\rho)\big) \, f^0(x,\rho, \rho_*) \ d\rho_*.
\end{equation}
The random path $\rho^0(x)$ may be constructed by solving
(deterministically) the ODE $d\rho^0/dx = b^0(x, \rho^0)$, interrupted by jumps which occur stochastically: the rate density at which $\rho^0$ makes a
jump at $x$ is $f^0( x,\rho^0(x), \rho_*)$. As our main result, we show that the process $x \mapsto \rho(x,t)$ (for fixed $t > t_0$)
is again a PDMP, with generator 
\begin{equation}
  \label{eq1.4}
  \big(\mathcal{A}_{x,t} \psi\big)(\rho) = b(x,t,\rho ) \psi'(\rho) + \int_\rho^\infty \big(\psi(\rho_*)
  - \psi(\rho)\big) \, f(x,t, \rho, \rho_*) \, d\rho_*.
\end{equation}
Here $b(x,t,\rho)$ and $f(x,t,\rho_-,\rho_+)$ are obtained from their initial
($t = t_0$) conditions 
\be\la{eq1.5}
b(x,t_0,\rho)=b^0(x,\rho),\ \ \  \ \ f(x,t_0,\rho_-,\rho_+)=f^0(x,\rho_-,\rho_+),
\ee
by solving a semi-linear PDE, 
\begin{equation}
  \label{eq1.6}
  b_t+H_xb_\rho-H_{\rho}b_x
  =H_{\rho\rho}b^2+2H_{\rho x}b+H_{xx},
  \end{equation}
and a kinetic (integro-)PDE 
\begin{equation}
  \label{eq1.7}
    f_t-(vf)_x - C(f) =Q(f) ,
\end{equation}
where 
\be\la{eq1.8}
v(x,t,\rho_-,\rho_+):=
\frac {H(x,t,\rho_-)-H(x,t,\rho_+)}{\rho_--\rho_+},
\ee
$Q(f)=Q^+(f)-Q^-(f)$ is a coagulation-like collision operator, and $C(f)=C^+(f)+C^-(f)$ is a 
linear first order differential operator.  More precisely, 

\ms\noi
{\bf(i)} $Q^+$ is a quadratic operator and
 $Q^+(f)(x,t,\rho_-,\rho_+)$ is defined as 
 \begin{align}
&\int_{\rho_-}^{\rho_+}
\big(v(x,t,\rho_*,\rho_+)
-v(x,t,\rho_-,\rho_*)\big)f(x,t,\rho_-,\rho_*)f(x,t,\rho_*,\rho_+)
\ d\rho_*.\la{eq1.9}
\end{align}

\ms\noi
{\bf(ii)} The quadratic operator $Q^-$ is of the form 
$Q^-(f)=f Jf$, for a linear operator $J$.  
Given $f$, the function $(Jf)(x,t,\rho_-,\rho_+)$ is defined as
\be\la{eq1.10}
A(vf)(x,t,\rho_+)-A(vf)(x,t,\rho_-)-
v(x,t,\rho_-,\rho_+)
\big((A f)(x,t,\rho_+)-(A f)(x,t,\rho_-)\big),
\ee
for linear operators $A$ defined by
\begin{align*}
A h(x,t,\rho_-)&=\int_{\rho_-}^\i h(x,t,\rho_-,\rho_+)\ d\rho_+.
\end{align*}

\ms\noi
{\bf(iii)} 
Given a $C^1$ kernel $f$, 
\be\la{eq1.11}
 \big(C^+f\big)(x,t,\rho_-,\rho_+)=\left[K(x,t,\rho_+,\rho_-)
f(x,t,\rho_-,\rho_+)\right]_{\rho_+},
\ee
where
\begin{align*}
K(x,t,\rho_+,\rho_-)&=b(x,t,\rho_+)v(x,t,\rho_-,\rho_+)-\b(x,t,\rho_+),\ \ \ \ {\text{ with }}\\
\b(x,t,\rho )&=\big(H_x+bH_\rho\big)(x,t,\rho ).
\end{align*}
Here and below, by the expression $X_a$ we mean the partial derivative of $X$ with respect to the variable 
$a$. For example the
right-hand side of \eqref{eq1.11} represents the partial derivative of the expression inside the brackets with respect to $\rho_+$.

\ms\noi
{\bf(iv)}
Given a $C^1$ kernel $f$, 
\begin{align}\nonumber
\big(C^-f\big)(x,t,\rho_-,\rho_+)=&b(x,t,\rho_-)
(vf)_{\rho_-}(x,t,\rho_-,\rho_+)-\b(x,t,\rho_-)f_{\rho_-}(x,t,\rho_-,\rho_+) \\
=&b(x,t,\rho_-)
\big(v_{\rho_-}f\big)(x,t,\rho_-,\rho_+)+
K(x,t,\rho_-,\rho_+)f_{\rho_-}(x,t,\rho_-,\rho_+) .\la{eq1.12}
\end{align}

\bs
\noi
{\bf Remark 1.1}
For a more compact reformulation of our equations
\eqref{eq1.6} and \eqref{eq1.7}, let us write
\be\la{eq1.13}
x_1=x,\ \ \ x_2=t,\ \ \  f^1=f,
\ \ \ f^2=v f^1,\ \ \ b^1=b,\ \ \ b^2=\b .
\ee
Recall $Ag(\rho)=A(g)(\rho)=\int g(\rho,\rho_*)\ d\rho_*,$ and
define
\be\la{eq1.14}
(g\otimes k)(\rho_-,\rho_+)=g(\rho_-,\rho_+) k(\rho_+),\ \ \ \
(k\otimes g)(\rho_-,\rho_+)=k(\rho_-)g(\rho_-,\rho_+).
\ee
A more symmetric rewriting of the equations
\eqref{eq1.6} and \eqref{eq1.7} read as
\be\la{eq1.15}
b^1_{x_2}-b^2_{x_1}=b^1b^2_\rho-b^2b^1_{\rho},\ \ \ \ \
f^1_{x_2}-f^2_{x_1}= \cQ(f^1,f^2)-\cQ(f^2,f^1),
\ee 
where
\be\la{eq1.16}
\cQ(f^j,f^i)=f^j*f^i-A(f^j)\otimes f^i-f^j
\otimes A(f^i)+b^j\otimes f^i_{\rho_-}-(f^j\otimes b^i)_{\rho_+},
\ee
where 
\[
(f^j*f^i)(\rho_-,\rho_+)=\int f^j(\rho_-,\rho_*)f^i(\rho_*,\rho_+)
\ d\rho_*.
\]
\qed

\bs
We now formulate our assumptions on the initial drift $b^0$,
the  initial jump
 rate kernel $f^0$, and  the Hamiltonian function $H(x,t,\rho)$.

\bs\noi
{\bf Hypothesis 1.1(i)}
The Hamiltonian function 
$H :\bR\x[t_0,T]\x [P_-,P_+] \to \bR$
 is  a $C^2$ function. Additionally, $H$ is increasing and  convex  in $\rho$.
    
\ms\noi
 {\bf(ii)} The PDE \eqref{eq1.6}
has a bounded $C^1$ solution $b\le 0$ for $t\in[t_0,T]$. We set  
$b^0(x,\rho):=b(x,t_0,\rho)$.

  \ms\noi
 {\bf(iii)} The PDE \eqref{eq1.7} has a   solution 
$ f:\hat \L\to[0,\i)$, where
$\hat \L:=\bR\x [t_0,T]\x \L(P_-,P_+)$, with
 \[
   \L(P_-,P_+):=\L\cap [-P,P]^2:= \big\{(\rho_-,\rho_+) :\
 P_-\le\rho_- \leq \rho_+ \leq P_+\big\}.
  \]
 We assume that $f$ is $C^1$ in the interior of $\hat\L$, and that
$f$ is continuous in $\hat \L$. Moreover,
$f(x,t, \rho_-,\rho_+)>0$, when $P_-<\rho_-<\rho_+<P_+$,
and $f(x,t, \rho_-,\rho_+)=0$, whenever $\rho_-$ or $\rho_+\notin(P_-,P_+)$.
To ease our notation, we extend the domain of the definition
of $f$ to $\bR\x[t_0,T]\x \bR^2$, by setting
 $f(x,t,\rho_-,\rho_+)=0$, whenever $\rho_-$ or $\rho_+\notin\L(P_-,P_+)$.
We also write $f^0(x,\rho_-,\rho_+)$ for 
$f(x,t_0,\rho_-,\rho_+)$.

\ms\noi
 {\bf(iv)} We assume that $\rho(x,t)$ is an entropy solution of 
\eqref{eq1.2}, and that its 
 initial condition $\rho^0(x):=\rho(x,t_0)$ is $0$ for $x<a_-$, 
and is a Markov process for $x\ge a_-$ that
  starts at $\rho^0(a_-) = m_0$. This Markov process  has an infinitesimal   generator in the form \eqref{eq1.3} for a drift
  $b^0$ and a jump rate density $f^0$.  
\qed

\bs
Our statistical description consists of a 
one-dimensional marginal, a drift,  and
a rate kernel generating the rest of the path.  The evolution of the
drift and the rate kernel are given by \eqref{eq1.6} and
the kinetic  equation \eqref{eq1.7}. Evolution of the marginal will be described in terms of the solutions
to  these equations. We continue with some definitions.

\ms\noi
{\bf Definition 1.1(i)} We define the linear operator $\cA^i$ by
\begin{equation}
  \label{eq1.17}
  \big(\cA^i_{x,t} \psi\big)(\rho) =  \big(\cA^{i} \psi\big)(\rho) = 
b^i(x,t,\rho ) \psi'(\rho) + \int_\rho^\infty f^i(x,t, \rho, \rho_+) \big(\psi(\rho_+)  - \psi(\rho)\big)  \, d\rho_+,
\end{equation}
for $i=1,2$. Note that $\cA^1=\cA$ of \eqref{eq1.4}, and $f^i$ was defined in \eqref{eq1.13}.
We write $\cA^{i*}$
 for the adjoint of the operator $\cA^i$
which acts on measures. When
the measure $\nu$ is absolutely continuous with respect to the Lebesgue measure with a $C^1$ Radon-Nykodym derivative, then $\cA^{i*}\nu$ is also absolutely continuous with respect to the Lebesgue measure. The action of the operator $\cA^{i*}$ on $\nu$ can be described in terms of its action on the corresponding Radon-Nykodym derivative. By a slight abuse of notation, we write 
$\cA^{i*}$ for the corresponding operator that now acts on $C^1$ functions.
More precisely, for a probability density $\nu$, we have
\begin{align*}
\big(\cA_{x,t}^{i*}\nu\big)(\rho)=&\left[\int_{-\i}^\rho 
f^i(x,t,\rho_*,\rho ) 
     \ \nu(\rho_*) \ d\rho_*\right]-A(f^i)(x,t,\rho)\nu(\rho)
 -\big(b^i(x,t,\rho )\nu(\rho)\big)_\rho.
 \end{align*}

    \ms\noi
    {\bf(ii)} We write $\cM$ for the set of measures and 
    $\cM_1$ for the set of probability measures.
    \qed

\begin{theorem}
  \label{th1.1}
  Given a $C^1$ rate $f$, and $m_0\in\bR$, 
  assume $\ell:[t_0,\i)\to\cM_1$ satisfies
  $\ell (t_0,d\rho_0)=\d_{m_0}(d\rho_0)$, and
  \be\la{eq1.18}
  \frac {d\ell}{dt}=\cA_{a_-,t}^{2*}\ell, 
  \ \ \ \ t>t_0.
  \ee
  When Hypothesis~1.1 holds, the
  entropy solution $\rho$ to \eqref{eq1.1}
   for each fixed $t > t_0$ has $x = a_-$ marginal given by
  $\ell(t,d\rho_0)$ and for $a_- < x < \infty$ evolves according to 
  a Markov process with the generator $\cA^1_{x,t}$.
Moreover, the process $t\mapsto \rho(a,t)$ is a Markov 
process with generator $\cA^2_{a,t}$, for every $a\ge a_-$.
\end{theorem}

\bs\noi
{\bf Remark 1.2(i)} According to Hypothesis~1.1, the function $H$ is increasing.  This condition is needed to guarantee that $f^2\ge 0$, which in turn guarantee that $\cA^2$ is a generator of a Markov process.
 This restriction on $H$ can be relaxed almost completely.
The main role of the condition $H_\rho>0$ is that all shock discontinuities of $\rho$ travel with negative velocities so that they 
cross any fixed location, say $x=a$ eventually. This allows us to assert that 
if $\rho(a,t)$ is known, then the law of $\rho(x,t)$ can be determined uniquely for all
$x>a$. In general, 
we may try to determine
$\rho(x,t)$ for $x>a(t)$, provided that $\rho(a(t),t)$ is specified. 
The condition $H_\rho>0$, allows us to choose $a(t)$ constant. If instead
we can find a negative constant $c$ such that 
$H_\rho> c$, then $\hat \rho(x,t):=\rho(x-ct,t)$ satisfies 
\[
\hat \rho_t=\hat H(x,t,\hat \rho)_x,
\]
for $\hat H(x,t,\rho)=H(x-ct,t,\rho)-c\rho$, which is increasing. Hence, the process
$t\mapsto\hat\rho(x,t)=\rho( x-ct,t)$ is now Markovian with a generator  $\hat\cA^2$ which is obtained from $\cA^2$
by replacing $H$ with $\hat H$. Even an upper bound on $H'$ can lead to a result similar to Theorem~1.2. For example if $H_\rho<0$,
then $x\mapsto \rho(x,t)$ is a Markov process but now as we decrease $x$.    

\ms\noi
{\bf(ii)} To guarantee the existence of a solution to \eqref{eq1.6}
in an interval $[t_0,T]$,
let us assume that  $H_{\rho x}$ and
$H_{xx}$ are uniformly bounded, and that $H_{xx}\le 0$
in this interval. Under such assumptions, we claim that if
initially at time $t=t_0$ the drift is nonpositive and bounded, then
the no blow-up condition of Hypothesis {\bf{1.1(ii)}} is met
 because $b$ remains bounded and nonnegative.
To see this,  assume that the function $b$ solves the equation \eqref{eq1.6}, and  write $\Theta_s^t(a,m )$ for the flow of the Hamiltonian ODE
\be\la{eq1.19}
\dot x=-H_{\rho}(x,t,\rho ),\ \ \ \ \ \dot\rho =H_x(x,t,\rho ) .
\ee
In other words $(\rho(t),x(t))=\Theta_s^t(a,m))$ solves \eqref{eq1.19}, subject to the initial conditions $x(s)=a,$ and 
$\rho(s)=m$. To ease the notation,
we write $b(x,\rho,t)$ and $H(x,\rho,t)$ for
$b(x,t,\rho),$ and $H(x,t,\rho)$ respectively. Evidently, 
$\hat b(x,\rho,t )=
b\big(\Theta_s^t(x,\rho),t\big)$ 
satisfies
\be\la{eq1.20}
\hat b_t=A\hat b^2+2B\hat b+C,
\ee
where
\[
A(x,\rho,t ):=H_{\rho\rho}\big(\Theta_s^t(x,\rho),t\big),\ \ \ \  
B(x,\rho,t ):=H_{\rho x}\big(\Theta_s^t(x,\rho),t\big),\ \ \ \ 
C(x,\rho,t ):=H_{xx}\big(\Theta_s^t(x,\rho),t\big).
\]
Since the right-hand side
of \eqref{eq1.20} is nonpositive when $\hat b=0$, we deduce that 
$\hat b(t)=\hat b(x,\rho,t)$ remains nonpositive for $t\in[t_0,T]$, if this is true initially at $t=t_0$. 
Note that since $\hat b_t\ge 2B\hat b+C$, with $B$ and $C$ bounded, $b$ is also bounded from below in $[t_0,T]$, if this is so initially.

\ms\noi
{\bf(iii)} The existence of a classical solutions to \eqref{eq1.7}
and \eqref{eq1.18} can be found in [KR2] and
[OR] when $H$ is independent 
of $(x,t)$, and $b$ is either constant, or $f$ is independent of $(x,t)$.
The same type of arguments can be worked out in our setting.

\ms\noi
{\bf(iv)} As a consequence of Hypothesis 1.1{\bf(iii)}, the density 
$\rho(x,t)\in[P_-,P_+]$ almost surely. This restriction will be needed 
in Sections 2 and 3 when we derive a forward equation for the law 
of $\rho(\cdot,t)$. The boundedness of $\rho(x,t)$ is needed only when we restrict $\rho$ to a bounded set of the form 
$\L:=[a_-,a_+]\x [t_0,T]$ (see Theorem 2.1 below). Note however that for a $C^1$
drift $b$, the density is always bounded below in $\L$, because the random 
jump only increases the density. So we only need to require an upper bound on the density, and the requirement  $P_-\le \rho_-$ is redundant. In other words, 
we can find $P_-$ that depends on $b$, and a lower bound of $\rho^0$,
such that the condition $P_-\le \rho_-$ holds in $\L$.

In Theorem 1.3 below, we will learn how to relax the boundedness requirement
on the density.
\qed

\subsection{Main result II}

Our Hypothesis 1.1{\bf(ii)} is rather stringent requirement because the right-hand side of  the PDE \eqref{eq1.6} is quadratic in $b$. Our main Theorem 1.1 applies only when no new shock discontinuity is created in the time interval $[t_0,T]$. Indeed a blow-up of the drift occurs exactly when a new jump discontinuity is formed for a local continuous solution that is represented by the ODE $\rho_x=b(x,t,\rho)$. In Remark 1.1{\bf{(ii)}} we stated conditions that would prevent a blow-up, but these conditions exclude many important stochastic growth models that are governed by HJE associated with random Hamiltonian (see Example 1.1{\bf(ii)} below). 

We emphasis
that Theorem 1.1 offers a kinetic description for the interaction between
the existing shock discontinuities, not those which are created after the initial time. To go beyond what is offered by Theorem 1.1, we need to enlarge the class of Markovian profiles that has been used so far. 
We offer a way to achieve this by considering profiles that are 
Markovian concatenations of {\em fundamental solutions} of \eqref{eq1.2}. 

\bs\noi
{\bf Definition 1.2} Given $z=(y,s)\in\bR^2$, by a {\em fundamental
solution} $W(\cdot;z):\bR\x(s,\i)\to\bR$ associated with $z$
we mean
\be\la{eq1.21}
W(x,t;z)=\sup\left\{\int_s^tL\big(\xi(\th),\th,\dot \xi(\th) \big)\ d\th:\ 
\xi\in C^1\big([s,t];\bR\big),\ \xi(s)=y,\ \xi(t)=x\right\},
\ee
where $L$ is the Legendre transform of $H$ in the
$p$-variable:
\[
L(x,t,v)=\inf_p\big(p\cdot v+
H(x,t,p)\big),\ \ \ \ H(x,t,p)=\sup_v\big(L(x,t,v)-p\cdot v\big).
\]
We also set $M(x,t;z)=W_x(x,t;z)$ for the $x$-derivative of $W$.
\qed

\bs
Under our conditions on $H$, the function $W$ is a Lipschitz function of
 $(x,t)$ for $t>s$,
and $M(x,t)=M(x,t;z)$ is well-defined a.e.. A representation of $M$ is given as follows. For each $(x,t)$, we may find a maximizing piecewise $C^1$ path 
$\xi(\th)=\xi(\th;x,t;z)$ that is differentiable at time $\th=t$. 
The function $M$ is continuous at $(x,t)$ if and only if
the maximizing path is unique. When this is the case,  we simply have
\be\la{eq1.22}
M(x,t)=L_v\big(\xi(t),t,\dot \xi(t)\big)=L_v\big(x,t,\dot \xi(t)\big).
\ee
In general $M(x,t)$ could be multi-valued; for each maximizing path, the right-hand side of \eqref{eq1.22} offers a possible value for $M(x,t)$.

The Cauchy problem associated with \eqref{eq1.1} has a representation
of the form
\be\la{eq1.23}
u(x,t)=\sup_y\big(u^0(y)+W(x,t;y,t_0)\big).
\ee
In other words, $u$ given by \eqref{eq1.23}, satisfies \eqref{eq1.1}
in viscosity sense for $t>t_0$, and $u(x,t_0)=u^0(x)$.
The type of stochastic solutions we will be able to describe kinetically 
would look like
\[
u(x,t)=\sup_{i\in I}\big(g_i+W(x,t;z_i)\big),
\]
where $\big\{(z_i,g_i): \ i\in I\big\}$ is a discrete set. Since our Markovian process is $\rho=u_x$, we consider
profiles of the form
\[
\rho(x,t)=\sum_{i\in I}M(x,t;z_i)\ \1\big(x\in [x_{i},x_{i+1})\big),
\]
for a discrete set $\big\{q_i=(x_i,z_i):\ i\in I\big\}$. (Note that because of the type of results we have in mind, we switched from $(g_i,z_i)$ to
$(x_i,z_i)$.)

We now give the definition of the Markov processes 
we will work with in this subsection.

\bs\noi
{\bf Definition 1.3(i)} Given $s,T$, with $s<T$, let $g(x,t;y_-,y_+)$ be a $C^1$
nonnegative (kernel) function that is defined 
for $x\in \bR,\ t\in [s,T],\ y_+\in (y_-,\i)$. We also write
\[
x_1=x,\ \ \ x_2=t,\ \ \ g^1=g,\ \ \ g^{2}=\hat v g,
\]
 where 
\be\la{eq1.24}
\hat v(x,t,y_-,y_+)=
\frac{H\big(x,t,M(x,t;y_+,s)\big)-H\big(x,t,M(x,t;y_-,s)\big)}
{M(x,t;y_+,s)-M(x,t;y_-,s)}.
\ee
We write $\cB^i_{x,t}$ for the operator 
\[
\cB^i_{x_1,x_2} F(y)=\int_{y}^\i \big(F(y_*)-F(y)\big)\ g^i(x_1,x_2;y,y_*)\ dy_*.
\]
$\cB^1$ is the infinitesimal generator of 
an {\em inhomogeneous Markov jump} process $\by(x_1)$. 
 When $\hat v>0$, the operator $\cB^{2}$ also generates a
Markovian jump process.

\ms\noi
{\bf(iii)}
We write $\cB^{i*}$
 for the adjoint of the operator $\cB^i$
which acts on measures. As before, when
the measure $\nu$ is absolutely continuous with respect to the Lebesgue measure with a $C^1$ Radon-Nykodym derivative, then $\cB^{i*}\nu$ is also absolutely continuous with respect to the Lebesgue measure. By a
 slight abuse of notation, we write 
$\cB^{i*}$ for the corresponding operator that now acts on $C^1$ functions. More precisely, for a probability density $\nu$, we have
\begin{align*}
\big(\cB_x^{i*}\nu\big)(y)=&\left[\int_{-\i}^y 
g^i(x,t,y_*,y ) 
      \nu(y_*)\ dy_*\right]-\hat A(g^i)(x,t,y)\nu(y),    
 \end{align*}
where 
\[
\hat A(g)(y)=\int_y^\i g(y,y_*)\ dy_*.
\]

\ms\noi
{\bf(iv)} Given $\by:[a_-,a_+]\to\bR$, we define
\[
\rho(x,t;\by,s):=M\big(x,t;\by(x),s\big).
\]
\qed

\bs
According to our second main result, if $\rho=u_x$ solves \eqref{eq1.2}
with an initial condition which comes from a Markov process associated with a kernel $g^0$, then at later times $x\mapsto\rho(x,t)$ also comes from
a Markov process associated with a kernel $g$ which satisfies a kinetic equation in the form
\be\la{eq1.25}
g_t-(\hat v g)_x=\hat Q(g)=\hat Q^+(g)-\hat Q^-(g)=\hat Q^+(g)-g\hat J(g),
\ee
where
\begin{align*}
&\hat Q^+(g)(y_-,y_+)=\int_{y_-}^{y_+} \big(\hat v(y_*,y_+)-\hat v(y_-,y_*)\big)
g(y_-,y_*)g(y_*,y_+)\ dy_*,\\
&\hat J(g)(y_-,y_+)=\big(\hat A(\hat vg)(y_+)-\hat A(\hat vg)(y_-)\big)
-\hat v(y_-,y_+)\big(\hat A(g)(y_+)- \hat A(g)(y_-)\big).
\end{align*}
Here we have not displayed the dependence of our functions on 
$(x,t)$ for a compact notation.

\bs
We are now ready to state our hypotheses and the second main result.

\bs\noi
{\bf Hypothesis 1.2(i)}  The Lagrangian function $L$ 
 is  $C^2$  function that is strictly concave in $v$.  Moreover, there are positive constants $c_0,c_1$ and $c_2$ such that
\begin{align*}
&-c_0+c_2v^2\le - L(x,t,v)\le c_0+c_1v^2,\\
&-c_0+c_2|v|\le | L_{v}(x,t,v)|\le c_0+c_1|v|,\\
&\ \ \ \ \ | L_{xv}(x,t,v)|+|L_{xx}(x,t,v)|\le c_1,\\  
&\ \ \ \ \  | H_{x}(x,t,v)|+|H_{x\rho}(x,t,\rho)|\le c_1.
\end{align*}

 \ms\noi
 {\bf(ii)} The rate kernel $g(x,t, y_-,y_+)$  is continuous nonnegative
solution of \eqref{eq1.25} which is $C^1$ in $(x,t)$-variable, 
and is supported on
  \[
    \big\{(x,t,y_-,y_+) : \ x\in\bR, \ t\in[t_0,T],\
Y_- \leq y_- \leq y_+ \leq Y_+\big\},
  \]
  for some constants $Y_\pm$. We write $g^0(x,y_-,y_+)$
for $g(x,t_0,y_-,y_+)$

\ms\noi
 {\bf(iii)}  $\rho\mapsto H_\rho(a_-,t,\rho)>0,$ for every $t\in[t_0,T]$
and $\rho\in [M_-,M_+]$, where
\[
M_+=\sup_{t\in[t_0,T]}M(a_-,t;Y_+,s),\ \ \ \ 
M_-=\inf_{t\in[t_0,T]}M(a_-,t;Y_-,s).
\]

\ms\noi
 {\bf(iv)}
 Given $s$ and $t_0$, with $t_0>s$, 
the initial condition $\rho^0(x) = M(x,t_0;y^0,s)$ for $x<a_-$, and 
$\rho(x,t_0)=\rho(x,t_0;\by_{t_0},s)$ for $x\ge a_-$, where $\by_{t_0}$
is a Markov process which
  starts at $\by_{t_0}(a_-) = y^0>a_-$, and
 has an infinitesimal   generator $\cB_{x,t_0}^1$, associated with a kernel 
$g^0(x,y_-,y_+)=g(x,t_0,y_-,y_+)$.

\ms\noi
 {\bf(v)}
Assume that $\ell:[0,\i)\to\cM_1$ satisfies
  $\ell (t_0,dy_0)=\d_{y^0}(dy_0)$, and
  \be\la{eq1.26}
  \frac {d\ell}{dt}=\cB_{a_-,t}^{2}\ell.
  \ee
\qed

\begin{theorem}
  \label{th1.2}  
 When Hypothesis~1.2 holds, the
  entropy solution $\rho$ to \eqref{eq1.2}
   for each fixed $t \in[ t_0,T]$ has $x = a_-$ marginal given by
$M(a_-,t;y_0,s)$, with $y_{0}$ distributed according to
  $\ell(t,dy_{0})$ and for $a_- < x $ evolves as
$\rho(x,t)=\rho(x,t;\by_t,s)$, with $\by_t$  
  a Markov process with the generator $\cB^1_{x,t}$.
\end{theorem}

\bs\noi
{\bf Remark 1.3(i)} Observe that the finiteness of the Lagrangian $L$ implies that Hamilonian function $H$ cannot be monotone $\rho$.
 As a consequence, the velocity 
$\hat v$ can take both negative and positive values, and the process $t\mapsto \rho(x,t)$ many not be a Markov process for every $x$.  However, when
 $x=a_-$, our Hypothesis 1.2{\bf(iii)} would guarantee that the 
process $t\mapsto \rho(a_-,t)$ is Markovian. Indeed 
Hypothesis 1.2{\bf(iii)} is designed to guarantee that no shock discontinuity can cross $a_-$ from left 
to right. This assumption though can be relaxed for the price of replacing
the boundary line segment $\{(a_-,t):\ t\in[t_0,T]\}$ with a suitable line segment 
which is titled to the right. In other words, part {\bf(i)} of Remark 1.2
is applicable. Moreover, part {\bf(iii)} of Remark 1.2
is also applicable to the kernel $g$ 
satisfying \eqref{eq1.25}.

\ms\noi
{\bf(ii)} As we will see in Proposition 5.2{\bf(iii)} in Section 5, 
there exist positive constants $C_0$ and $C_1$ such that  
 $M(x,t;y,s)\ge -C_1x$ for $x\le -C_0$. Our condition 
$|L_v|\le c_1(1+|v|)$ in Hypothesis 1.2{\bf(i)}, means 
\[
|\rho|\le c_1\big(1+|H_\rho(x,t,\rho)|\big).
\]
 Since $H_\rho(x,t,\rho)$
is an increasing function of $\rho$, we deduce that $H_\rho(x,t,\rho)\to\pm\i$
as $\rho\to\pm \i$. From this we learn that there exists a positive constant $C_2$ such that $H_\rho(a_-,t,\rho)>0$ whenever  $\rho\ge C_2$. As a consequence, 
$H_\rho(a_-,t,\rho)>0$ for $\rho\in[M(a_-,t;y_-,s),M(a_+,t;y_+,s) ]$, provided that $a_-\le -C_2C_1^{-1}$. This means that Hypothesis 1.2{\bf(iii)}
is automatically satisfied when $a_-\le -C_2C_1^{-1}$.

\ms\noi
{\bf(iii)} As a concrete example, when $H(x,t,\rho)=\rho^2/2$, then
$L(x,t,v)=-v^2/2$, and $M(x,t;y,s)=(y-x)/(t-s)$. In this case
Hypothesis 1.2{\bf(iii)} holds if and only if $a_-<Y_-$.
\qed

\bs\noi
{\bf Example 1.1} When $H$ does not depend on $(x,t)$, then
\[
W(x,t;y,s)=(t-s)L\left(\frac{x-y}{t-s}\right),\ \ \ \ 
M(x,t;y,s)=L'\left(\frac{x-y}{t-s}\right).
\]

\bs\noi
{\bf Remark 1.4}  As an example for a stochastic growth model, we may consider
$H(x,t,\rho)=H_0(\rho)-V(x,t)$, with $H_0$ convex, and the potential
$V$ given formally as
\be\la{eq1.27}
V(x,t)=\sum_{i\in I}\d_{s_i}(t)\1(x=a_i),
\ee
where $\o=\big\{(a_i,s_i):\ i\in I\big\}$, is a realization of a 
{\em Poisson Point Process} in $\bR^2$. In practice,
we may approximate $V$ by 
\[
V_\e(x,t)=\sum_{i\in I}\e \z\left(\frac{t-s_i}{\e}\right)
\eta\left(\frac{x-a_i}{\d(\e)}\right),
\]
where $\d(\e)\to 0$, in small $\e$-limit, and
$\eta$ and $\z$ are two smooth functions of compact support such that
$\int\z(t)\ dt=1$, and $\eta(x)=1$ in a neighborhood of the origin.
Replacing $V$ with $V_\e$ yields a Hamiltonian function
$H^\e$ for which the equation \eqref{eq1.1} is well-defined
and its solution $u^\e$ has a limit $u$ as $\e\to 0$.
A variational representation as in \eqref{eq1.21} for $u^\e$ would
yield a variational representation for $u$ as well. Indeed the corresponding $W$ still has the form \eqref{eq1.21}, where
$L(x,t,v)=L_0(v)-V(x,t)$, with $L_0$ a concave function given by
\[
L_0(v)=\inf_p\big(p\cdot v+H_0(p)\big).
\]
It is not hard to show that the minimizing path $\xi$ of the variational
problem \eqref{eq1.21} is a concatenation of line segments
between Poisson points of $\o$. In other words, 
\be\la{eq1.28}
W(x,t;y,s)=W(x,t;y,s;\o)=\sup\left(N(\bz)+\sum_{i=0}^{N(\bz)}(s_{i+1}-s_i)
L_0\left(\frac{a_{i+1}-a_i}{s_{i+1}-s_i}\right)\right),
\ee
where the supremum is over sequences
$\bz=
\big((a_0,s_0),(a_1,s_1),\dots, (a_n,s_n), (a_{n+1},s_{n+1})\big)$,
such that $N(\bz)=n$, and
\[
s_0<s_1<\dots<s_{n+1},\ \ \ (a_0,s_0)=(y,s),\ \ \  (a_{n+1},s_{n+1})=(x,t), \ \ \ (a_1,s_1),\dots, (a_n,s_n)\in\o.
\]
This model was defined and studied in Bakhtin [B] and Bakhtin et al. [BCK] when $H_0(p)=p^2/2$ (which leads to $L_0(v)=-v^2/2$).
If $H_0(p)=|p|$, then 
$L_0(v)=-\i\ \1(|v|>1)$. In this case, 
\[
W(x,t;y,s)=W(x,t;y,s;\o)=\sup N(\bz),
\]
where the supremum is over sequences $\bz$ as in \eqref{eq1.28},
with the additional requirement
\[
|a_{i+1}-a_i|\le s_{i+1}-s_i.
\]
The corresponding $u(x,t)$ is a stochastic growth model that is known
as {\em Polynuclear Growth} (We refer to [PS] for more details).

Our Theorem 1.2 does not directly apply to this model because
  Hypothesis 1.2{\bf(i)} fails. Also for Hypothesis 1.2{\bf(ii)} to hold, we need
to assume that the intensity of $\o$ is $0$ outside
$[a_-,\i)\x \bR$. Nonetheless, our method of proof can be adopted to treat
 this model as well. For this model however, it is more natural to consider
a concatenation of fundamental solutions $M(x,t;y_i,\th_i)$, where  
$\{(y_i,\th_i):\ i\in I\}$ is selected randomly. This extension requires developing
new techniques and goes beyond the scope of the present article.
\qed

   \subsection{Unbounded density}

In Theorem 1.1 (respectively 1.2), we assumed that the density $\rho$ 
(respectively $y$) is bounded. 
This assumption is technically convenient for the derivation of 
the forward equation that is carried out in Section 3, and is at the heart of our proofs of Theorems 1.1 and 1.2. Unfortunately 
it excludes many important models encountered in statistical mechanics, especially when we study stochastic growth models. 
As an example, if we take the case of the Burgers equation with white noise initial data, the density at later times would be an unbounded
Markov jump process (see [Gr], [MS], and [OR]).
In this subsection we explain how one can relax this restriction with the aid of an approximation that is related to Doob's $h$-transform. We carry out this idea 
in the case of Theorem 1.2 only, though our method of proof is also applicable 
to the setting of Theorem 1.1. 

Imagine that we have a kernel $g$ which satisfies
the kinetic equation \eqref{eq1.25}, and the arguments $y_\pm$
are not restricted to a bounded interval as in Hypothesis 1.2{\bf(ii)}. 

\bs\noi
{\bf Hypothesis 1.3} We assume that parts {\bf(i)} and {\bf(iii)}-{\bf(iv)} of 
Hypothesis 1.2 hold, but in part {\bf(ii)}, we allow $Y_+=\i$, and assume that
$g(x,t,y_-,y_+)$ is a continuous kernel such that the 
 Markov process
$\by_{t_0}$ associated with the generator $\cB^1_{x,t_0}$ satisfies
\be\la{eq1.29}
\limsup_{x\to\i}x^{-1} y_{t_0}(x)<1,
\ee 
almost surely.
\qed

\bth\la{th1.3} The conclusion of Theorem 1.2 holds even when $g$ is
a kernel which satisfies Hypothesis 1.3.
\et

Our strategy for proving Theorem 1.3 is by approximating the kernel $g$
with a sequences of kernels $g^n$ for which Theorem 1.2 is applicable.
We cannot simply restrict $g$ to a large bounded interval, because the resulting kernel does not satisfy the kinetic equation. However, if $\by$ is a Markov
process with the jump kernel density $g$ (associated with the generator $\cB^1_{x,t}$ as in the Definition 1.3{\bf{(i)}}),
we may condition $\by$ to remain in a bounded interval. The resulting process
is again a Markov process for which the jump kernel $\hat g$ is related to
$g$ via a Doob's $h$-transform. In other words, there exists a suitable
function $h(x,t,y)$ such that 
\be\la{eq1.30}
\hat g(x,t,y_-,y_+)=\frac{h(x,t,y_+)}{h(x,t,y_-)}g(x,t,y_-,y_+)
=:\eta(x,t,y_-,y_+)g(x,t,y_-,y_+).
\ee
Indeed the resulting kernel is again a solution to the kinetic equation as the following result confirms:

\bp\la{pro1.1} Assume $g$ satisfies \eqref{eq1.25} and $h:[a_-,a_+]\x
[t_0,T]\x\bR\to\bR$ is a $C^1$ function such that
\be\la{eq1.31}
h_{x}+\cB_{x,t}^1h=0,\ \ \ \ h_{t}+\cB_{x,t}^2h=0.
\ee
 Then $\hat g$ given by \eqref{eq1.30} also satisfies
\eqref{eq1.25}.
\ep

\bs
As we will see in Subsection 1.4 below, the two equations that appeared 
in \eqref{eq1.31} are compatible whenever $g$ satisfies the equation
\eqref{eq1.25}. This means that one of these equations is redundant:

\bp\la{pro1.2} Assume that $g$ and $h$ are bounded, and  $C^1$ 
in $(x,t)$ variable, and that $g$ satisfies \eqref{eq1.25}.
Also assume that $h$ is uniformly positive, 
satisfies $h_{x}+\cB^1_{x,t}h=0$, and
\be\la{eq1.32}
h_{t} (a_-,t,y)+(\cB_{a_-,t}^2 h)(a_-,t,y)=0.
\ee
(In other words, the second equation in \eqref{eq1.31} holds at
$x=a_-$.) Then the second equation in 
\eqref{eq1.31} holds in $[a_-,a_+]$.
\ep

\bs
The proof of Proposition 1.2 is similar to the proof of Proposition 4.1
of [OR], and is omitted.

\subsection{Heuristics}

According to Theorem 1.2, the process $x_i\mapsto\rho(x_1,x_2)$
is a Markov process with the generator 
\[
\cA^i\psi(\rho)=\cA^i_{x_1,x_2}\psi(\rho)=b^i\rho_{x_i}+\int_{\rho}^\i f^i(x_1,x_2,\rho,\rho_*)
\big(\psi(\rho_*)-\psi(\rho)\big)\ d\rho_*.
\]
Hence, if  $\ell(x_1,x_2,\rho)$ denotes the probability density of 
$\rho(x_1,x_2)$, then $\rho$ must satisfy the {\em forward equation}
\begin{align*}
&\ell_{x_i}=\cA^{i*}\ell=\ell*f^i-A(f^i)\ell-(b^i\ell)_{\rho}=
\ell*f^i-A(\ell\otimes f^i)-(b^i\ell)_{\rho},\ \ \ \ i=1,2,
\end{align*}
where
\[
(\ell*f^i)(\rho)=\int \ell(\rho_*)f^i(\rho_*,\rho)\ d\rho_*.
\]
From differentiating both sides we learn
\begin{align*}
\ell_{x_1x_2}=&\cA^{1*}(\ell_{x_2})
+\ell *f^1_{x_2}-A(f^1)_{x_2}\ \ell-(b^1_{x_2}\ell)_{\rho}
=\cA^{*1}\cA^{2*}\ell
+\ell *f^1_{x_2}-A(f^1)_{x_2}\ \ell-(b^1_{x_2}\ell)_{\rho},\\
\ell_{x_2x_1}=&\cA^{2*}(\ell_{x_1})+\ell *f^2_{x_1}-A(f^2)_{x_1}\ \ell-(b^2_{x_1}\ell)_{\rho}
=\cA^{2*}\cA^{1*}\ell
+\ell *f^2_{x_1}-A(f^2)_{x_1}\ \ell-(b^2_{x_1}\ell)_{\rho}.
\end{align*}
On the other hand,
\begin{align*}
\cA^{i*}\cA^{j*} \ell=&\left(\ell*f^j-A(f^j)\ell-(b^j\ell)_\rho\right)*f^i-A(f^i)\left(\ell*f^j-A(f^j) \ell-(b^j\ell)_\rho\right)\\
&-\left[b^i\left(\ell*f^j-A(f^j)\ell-(b^j\ell)_\rho\right)\right]_\rho\\
=&\ell*\left[f^j*f^i-A(f^j)\otimes f^i-f^j
\otimes A(f^i)+b^j\otimes f^i_{\rho_-}-(f^j\otimes b^i)_{\rho_+}
\right]
\\
&+\ell \left[ A(f^i)A(f^j)+A(f^i)b^j_\rho+(A(f^j)b^i)_\rho\right]
\\
&+\ell_\rho\left[A(f^i)b^j +A(f^j)b^i\right]+\left(b^ib^j_{\rho}\ell \right)_{\rho}
+\left(b^ib^j\ell_{\rho}\right)_{\rho}.
\end{align*}
(Here, we have performed an integration by parts to replace
$(b^j\ell)_\rho*f^i$ with $\ell*\big(b^j\otimes f^i_{\rho_-}\big)$.)
As a result
\begin{align*}
\cA^{1*}\cA^{2*}\ell-\cA^{2*}\cA^{1*}\ell=&\ell*\left[\cQ(f^2,f^1)-\cQ(f^1,f^2)\right]
+\left[\left( b^1b^2_\rho-b^2b^1_{\rho}\right)\ell\right]_\rho\\
&+\ell\left[A(f^2)_\rho \ b^1-A(f^1)_\rho \ b^2\right],
\end{align*}
where $\cQ(f^j,f^i)$ is given by \eqref{eq1.16}.
Hence,
\begin{align}\nonumber
\ell_{x_1x_2}-\ell_{x_2x_1}=&\ell*[\cQ(f^2,f^1)-\cQ(f^1,f^2)+f^1_{x_2}-f^2_{x_1}]
\\
&-\ell\left[A(f^1)_{x_2}-A(f^2)_{x_1}+A(f^1)_\rho \ b^2
-A(f^2)_\rho \ b^1\right]\la{eq1.4.1}\\
&+\left[\left( b^2_{x_1}-b^1_{x_2}+b^1b^2_\rho-b^2b^1_{\rho}
\right)\ell\right]_{\rho}\ .\nonumber
\end{align} 
It is rather straightforward to show
\[
A\left(\cQ(f^j,f^i)\right)=-A(f^j) A(f^i)+b^j A(f^i)_{\rho}.
\]
As a consequence,
\be\la{eq1.4.2}
A\left(R(f^1,f^2)\right)=A(f^1)_{x_2}-A(f^2)_{x_1}+A(f^1)_\rho \ b^2
-A(f^2)_\rho \ b^1.
\ee
where
\[
R=R(f^1,f^2)=f^1_{x_2}-f^2_{x_1}+\cQ(f^2,f^1)-\cQ(f^1,f^2).
\]
From \eqref{eq1.4.1} and \eqref{eq1.4.2} we deduce
\[
\ell_{x_1x_2}-\ell_{x_2x_1}=\ell*R-A(R)\ell+(S\ell)_\rho,
\]
where
\[
S=S(b^1,b^2)=b^2_{x_1}-b^1_{x_2}+b^1b^2_\rho-b^2b^1_{\rho}.
\]
Clearly $R=S=0$ implies the compatibility of the equations
$\ell_{x_i}=\cA^{i*}\ell$, for $i=1,2$. Observe that
$R=S=0$ are exactly our equations \eqref{eq1.6} and \eqref{eq1.7}.

\bs
Various terms in the kinetic equation can be readily explained in terms of the underlying particle system that represents the dynamics
of the shock discontinuities of a solution to the PDE \eqref{eq1.1}.

\ms\noi
{\bf(1)} According to the generator \eqref{eq1.1}, the process
$x\mapsto \rho(x,t)$  satisfies the ODE
\be\la{eq1.4.3}
\rho_x(x,t)=b\big(x,t,\rho(x,t)\big),
\ee
in between shock discontinuities. The PDE \eqref{eq1.6},
governing the evolution of the velocity $b$, follows from the consistency of \eqref{eq1.22} with \eqref{eq1.1}; differentiating these equations with respect to $t$ and $x$ respectively lead to
\begin{align*}
\rho_{xt}&=b_t+b_\rho \ H(x,t,\rho)_x=b_t+b_\rho \big(H_x+H_{\rho}b\big)=b^1_{x_2}+b^2b^1_\rho,\\
\rho_{tx}&=H(x,t,\rho)_{xx}=\big(H_x+H_{\rho}b\big)_x
=H_{xx}+2H_{x\rho}b+H_{\rho\rho}b^2+H_\rho b_x+H_\rho
b_\rho b=b^2_{x_1}+b^1b^2_\rho.
\end{align*}
Matching these two equations yields \eqref{eq1.6}. This calculation
is simply a repetition of the derivation of the equation $S(b^1,b^2)=0$.

\ms\noi
{\bf(2)}
If a shock discontinuity occurs at a location $x(t)$ with
$\rho_{\pm}(t)=\rho\big(x(t)\pm,t\big)$, then by
the classical Rankine-Hugoniot equation
\be\la{eq1.4.4}
\dot x(t)=-v\big(x(t),t,\rho_-(t),\rho_+(t)),
\ee
where $v$ was defined in \eqref{eq1.8}. This equation 
is responsible for the occurrence of
the term $-(vf)_x$ in  \eqref{eq1.7}.

\ms\noi
{\bf(3)} Since $\rho(x,t)$ solves \eqref{eq1.1} classically
away from the jump discontinuities, we have
\begin{align}\nonumber
\dot \rho_+(t)&=-\rho_x\big(x(t)+,t\big)v\big(x(t),t,\rho_-(t),\rho_+(t))+\rho_t\big(x(t)+,t\big)\\
&=-b\big(x(t),t,\rho_+(t)\big)v\big(x(t),t,\rho_-(t),\rho_+(t)) \nonumber
+(H_x+bH_\rho )(x(t),t,\rho_+(t))\\
&=-K\big(x(t),t,\rho_+(t),\rho_-(t)).\la{eq1.4.5}
\end{align} 
As in {\bf(2)}, this equation is responsible for the occurrence of 
$-C^+f$ in \eqref{eq1.7} (see \eqref{eq1.11} for the definition of $C^+$).

\ms\noi
{\bf(4)} A repetition of our calculation in {\bf(3)} yields
\be\la{eq1.4.6}
\dot \rho_-(t)=-K\big(x(t),t,\rho_-(t),\rho_+(t)).
\ee
Based on this, we are tempted to guess that $C_-f$ is
$\big[K\big(x,t,\rho_-,\rho_+)f\big(x,t,\rho_-,\rho_+)\big]_{\rho_-}.$
This is not what we have in \eqref{eq1.12}. The reason behind this has to do with the fact that we regard 
$\rho(x,t)$ as a Markov process in $x$ as we increase $x$. As a result, the role of $\rho_-$ and $\rho_+$ cannot be interchanged. In order to explain the form of $C^-f$ in \eqref{eq1.12}, we fix $a\in \bR$, 
and assume that $x(t)$ is the first discontinuity which occurs to the right of $a$. Now, if we set
$\rho_0(t)=\rho(a,t)$, and write 
\[
\rho(x)=\phi_a^x(m_0;t),
\]
 for the flow of the ODE \eqref{eq1.4.3} (in other words $\rho(x)$ solves \eqref{eq1.4.3} subject to the initial condition
$\rho(a)=m_0$), then
\[
\rho_-(t)=\phi_a^{x(t)}\big(\rho_0(t);t\big).
\]
Since $\rho_0(t)$ satisfies $\dot\rho_0=\b(a,t,\rho_0)$, its law $\ell (t,\rho_0)$ obeys the equation
\[
\ell_t(t,\rho_0)+\big(\b(a,t,\rho_0)\ell(t,\rho_0)\big)_{\rho_0}=0,
\]
 away from the shock discontinuity. As it turns out, the function
\[
k(x,t,\rho_0,\rho_+):=\ell(t,\rho_0)
f\big(x,t,\phi_a^x(\rho_0;t),\rho_+\big),
\]
satisfies the identity
\[
k_t-(wk)_x-(\b k)_{\rho_0}=\ell\big(f_t-(vf)_x-C^-f\big),
\]
where 
\[
w(x,t,\rho_0,\rho_+)=v\big(x,t,\phi_{a}^{x}(\rho_0;t),\rho_+\big).
\]

\ms\noi
{\bf (5)} Observe that if a solution $\rho$ has two jump discontinuities
at $x=x(t)$ and $y=y(t)$, with $x<y$, and
\[
\rho_-=\rho(x-,t),\ \ \ \ \rho_*=\rho(x+,t),
\ \ \ \ \rho'_*=\rho(y-,t),\ \ \ \ \rho_+=\rho(y+,t),
\]
then the relative velocity of these two discontinuities is exactly
\[
v\big(x,t,\rho_-,\rho_*\big)-v\big(y,t,\rho'_*,\rho_+\big).
\]
As $y(t)$ catches up with $x(t)$, $\rho'_*$ converges
to $\rho_*$ and the relative velocity becomes
\[
v\big(x,t,\rho_-,\rho_*\big)-v\big(x,t,\rho_*,\rho_+\big).
\]
This explains the form of $Q^+$ in \eqref{eq1.9}.
\qed

\subsection{Bibliography and the outline of the paper}

Most of the earlier works on stochastic solutions of Hamilton-Jacobi PDEs have been carried out in the Burgers context. For example,
 Groeneboom [Gr] determined the statistics of
  solutions to Burgers equation ($H(p)=p^2/2$, $d=1$)
with white noise initial data.
Recently Ouaki [O] has extended this result to arbitrary convex Hamiltonian function $H$. The special cases of 
$H(p)=\i \1(p\notin[-1,1])$, and $H(p)=p^+$ were already studied in the references Abramson-Evans [AE], Evans-Ouaki [EO], and Pitman-Tang [PW].

 Carraro and Duchon [CD1-2] considered
  \emph{statistical} solutions, which need not coincide with genuine
  (entropy) solutions, but realized in this context that L\'evy
  process initial data  should
  interact nicely with Burgers equation.  Bertoin [Be]
  showed this intuition was correct on the level of entropy solutions,
  arguing in a Lagrangian style.

Developing an alternative treatment to that given by Bertoin, which
relies less on particulars of Burgers equation and happens to be more
Eulerian, was among the goals of  Menon and Srinivasan [MS].
Most notably, [MS] formulates an interesting conjecture for the evolution
of the infinitesimal generator of the solution $\rho(\cdot,t)$ which is
equivalent with our kinetic equation \eqref{eq1.7} when $H$ is independent of $(x,t)$.
When the initial data
$\rho(x,0)$ is allowed to assume values only in a fixed, finite set of
states, the infinitesimal generators of the processes $x\mapsto \rho(x,t)$ and $t\mapsto \rho(x,t)$ 
can be represented by triangular matrices.  
The integrability of this matrix evolution
has been investigated by Menon [M2] and Li [Li].  For generic matrices---where the genericity
assumptions unfortunately exclude the triangular case---this evolution
is completely integrable in the Liouville sense.  
The full treatment of Menon and Srinivasan's conjecture was achieved
in papers [KR1] and [KR2] (we also refer to [R] for an overview).
The work of [KR1] have been recently extended to higher dimensions in [OR1]. In [OR2], the main result of [KR2] has been used to give a new 
proof of Groeneboom's results [Gr].

\bs
We continue with an outline of the paper:

\ms\noi
{\bf(i)} In Section 2, we show that the evolution of the PDE \eqref{eq1.1} for piecewise smooth solutions is equivalent to a particle system in $\bR\x[P_-,P_+]$. We restrict this particle system to a large finite interval
$[a_-,a_+]$ and introduce a stochastic boundary condition at $a_+$. This restriction allows us to reduce Theorem 1.1 to a finite system; the precise statement can be found in Theorem 2.1 of Section 2. 

\ms\noi
{\bf(ii)} The strategy of the proof of Theorem~2.1 will be described in Section~3. Our strategy is similar to the one that was utilized in our previous work [KR1-2]: Since we have a candidate for the generator of the process $x\mapsto \rho(x,t)$, we have a candidate measure, say $\mu(\cdot,t)$ for the law of $\rho(\cdot,t)$.
  We establish Theorem~2.1 by showing that this candidate measure satisfies the {\em forward equation} associated with Markovian dynamics of the underlying
particle system (see the equation \eqref{eq3.2} in Section 3). The particle system has a
 deterministic evolution inside the interval and a stochastic (Markovian) dynamics at the right end boundary point. 
The rigorous derivation of the forward equation will be carried out in Section~3.
 
\ms\noi
{\bf(iii)} Section~4 is devoted to the proof of Theorem~2.1.

\ms\noi
{\bf(iv)} Section~5 is devoted to the proof of Theorem~1.2.

\ms\noi
{\bf(v)} In Section~6, we establish Theorem 1.3 and Proposition 1.1.
\qed

\section{Particle System}

We assume that the initial condition $\rho^0$, 
in the PDE \eqref{eq1.21} is of the following form
\bi
\item  $\rho^0(x)=m_0$ for $x\le x_0=a_-$.
\item There exists a discrete set $I^0
=\{x_i: i\in\bN\}$, with $a_-<x_1<\dots<x_i<\dots$
such that for every $x>0$ with $x\notin I^0$, we have $\rho_x^0(x)=b^0(x,t_0,\rho^0(x))$. 

\item If $\rho_i^\pm=\rho^0(x_i\pm)$ denote the right and left
values of $\rho^0$ at $x_i$, then $\rho_i^-<\rho^+_i$.
\ei
Now if $\rho$ is an entropic solution of \eqref{eq1.2} with initial
$\rho^0$, then we may apply the {\em method of characteristics}
to show
that for each $t\ge 0$, the function $\rho(\cdot ,t)$
 has a similar form. To explain this, consider the ODE
\be\la{eq2.1}
 \frac d{dx}\rho(x)=b(x,t,\rho(x)),
 \ee
 where $b$ is the  solution to \eqref{eq1.6}, subject to the initial
 condition $b(x,t_0,\rho)=b^0(x,\rho)$. 
Recall that we are write $\phi_a^z(m;t)$ for the flow of the ODE \eqref{eq2.1}.
 In other words, if $\rho(x)=\phi_a^x(m;t)$, then \eqref{eq2.1} holds, and $\rho(a)=m$.
 Then there are pairs 
 $\bq(t)=\big((x_i(t),\rho_i(t)):\ i=0,1,\dots\big)$, with
 \[
 a_-=x_0(t)<x_1(t)<\dots<x_i(t)<\dots,
 \]
 such that for $x\ge a_-$, we can write
 \be\la{eq2.2}
 \rho\big(x,t\big)=\sum_{i=0}^{\i}\phi_{x_i(t)}^x
 \big(\rho_i(t);t\big)
1\!\!1\big(x_i(t)\le x<x_{i+1}(t)\big).
\ee
Note that $\rho(x_i(t)+,t)=\rho_i(t)$, and the data 
 $\bq(t)$ determines $\rho(\cdot,t)$ completely.    
  Because of this, we can fully describe the evolution of $\rho(\cdot,t)$ by describing the evolution of the particle system $\bq(t)$. Indeed from the PDE \eqref{eq1.1} and the
Rankine-Hugoniot Formula, we have
$\dot\rho_0(t)=\b(a,t,\rho_0(t))$, $\rho_0(t_0)=m_0$, and 
\be\la{eq2.3}
\dot x_i(t)=-v(x_i(t),t,\hat\rho_{i-1}(t),\rho_i(t)),\ \ \ \ \ 
\dot\rho_i(t)=-K(x_i(t),t,\rho_i(t),\hat\rho_{i-1}(t)),
\ee
for $i\in\bN$,
where $\hat\rho_{i-1}(t)=\phi_{x_{i-1}(t)}^{x_i(t)}(\rho_{i-1}(t),t)$
(we refer to Subsection 1.4, especially \eqref{eq1.4.4} and \eqref{eq1.4.5}
for explanation).
 Here  \eqref{eq2.3} gives a complete description of
$\bq$ in an inductive fashion; once $(x_{i-1}, \rho_{i-1})$ is determined, then we use \eqref{eq2.3} to write a system of two equations for the pair $(x_i,\rho_i)$.
Moreover \eqref{eq2.3} holds so long as $x_i'$s do not collide. When there is a collision between $x_i$ and $x_{i+1}$, for some $i=0,1,\dots$, we remove $x_{i+1}$ from the system,
replace $\rho_i$ with $\rho_{i+1}$, and relabel $(x_j,\rho_j)$
as $(x_{j-1},\rho_{j-1})$ for $j>i+1$.  As we will see shortly,
the function $\rho(x,t)$, defined by equation \eqref{eq2.2}, with $\bq(t)$ evolving as above, is the unique entropy solution of \eqref{eq1.2}.
 
According to Theorem~1.1 if $\rho(\cdot,t_0)$ is a PDMP 
with drift $b^0$ and jump rate $f^0$, then $\rho(\cdot,t)$ is also
a PDMP with drift $b(x,t,\cdot)$ and $f(x,t,\cdot,\cdot)$. 
We may translate this into a statement about the law of
our particle system $\bq(t)$. However, since the dynamics 
of $\bq$ involves infinite number of particles, we may take advantage of the finiteness of propagation speed in \eqref{eq1.2} and reduce Theorem~1.1 to an analogous claim for a finite interval $[a_-,a_+]$. 

Since $H_\rho>0$ by Hypothesis~1.1{\bf(ii)}, all particles travel to left. Because of this,
 we need to choose appropriate boundary dynamics
at the right boundary $a_+$ only.
The involved analysis will all pertain to the
following result.

\begin{theorem}
  \label{th2.1}
  Assume Hypothesis~1.1.  For
  any fixed $a_+ > a_-$, consider the scalar conservation law
\eqref{eq1.2} in $[a_-,a_+] \times [t_0,T)$ 
   with initial condition $\rho(x,t_0)=\rho^0(x)$ 
(restricted to $[a_-,a_+]$), open
  boundary at $x = a_-$, and random
  boundary $\zeta$ at $x = a_+$.  Suppose the process $\zeta$ has
  initial condition $\zeta(t_0) = \rho^0(a_+)$, and evolves according to the time-dependent rate  kernel $ f^2(a_+,t,\rho_,\rho_+)$
 and drift $b^2(a_+,t,\rho)$, independently of $\rho^0$.
  Then for all $t > t_0$  the law of $\big(\rho(x, t):\ x\in[a_-,a_+]\big)$ is  as follows:
\bi
\item[{\bf(i)}] The $x = a_-$ marginal is $\ell(t,d\rho_0)$, given by 
$\dot \ell=\cA^{2*}_{a_-,t}\ell$.
 
\item[{\bf(ii)}] The rest of the path is a PDMP with generator 
$\cA^1_{x,t}$
 (rate kernel $f(x,t,\rho_-,\rho_+)$  and drift $b(x,t,\rho)$).
\ei
\end{theorem}

To prove our main result Theorem~1.2,
we can send $a_+ \to \infty$, applying
Theorem~2.1 on each $[a_-,a_+]$, and use bounded speed
of propagation.  The argument is straightforward
 and can be found in [KR1].

We prove Theorem~2.1 by showing that the particle system
$\bq(t)$ restricted to the interval $[a_-,a_+]$ has the correct law predicted by this theorem. We now 
give a precise description for the evolution of $\bq$ restricted
to $[a_-,a_+]$.  First we make some definitions.

 \ms\noi
 {\bf Definition 2.1(i)} The configuration space for our particle system $\bq$, is the set $\D=\cup_{n=0}^\i\bar\D_n,$
where $\bar\D_n$ is the topological closure of $\D_n$, with $\D_n$ denoting the set
\[
\big\{\bold{q}=\big((x_i,\rho_i):i=0,1,\dots,n\big):\ x_0=a_-<x_1<\dots<x_n<x_{n+1}=a_+,
\ \ \ \rho_0,\dots,\rho_n\in\bR\big\}.
\]
We write $\bn(\bq)$ for the number of particles i.e.,
$\bn(\bq)=n$ means that $\bq\in\D_n$.
What we have in mind is that $\rho_i(t)=\rho(x_i(t)+,t)$ with $x_1,\dots,x_n$ denoting the locations of all shocks in $(a_-,a_+)$. 

\ms\noi
{\bf(ii)}
Given a realization 
$\bold{q}=\big(x_0,\rho_0,x_1,\rho_1,\dots,x_n,\rho_n\big)\in \bar\D_n,$ we define
\begin{align*}
\rho\big(x,t;\bold{q}\big)&=
R_t(\bq)(x)=\sum_{i=0}^{n}\phi_{x_{i-1}}^{x_i}
\big(\rho_i;t\big)
1\!\!1\big(x_i\le x<x_{i+1}\big).
\end{align*}

\ms\noi
{\bf(iii)} The process $\bq(t)$ evolves according to the following 
rules:
\bi
\item[$\bold{(1)}$] So long as $x_i$ remains in $(x_{i-1},x_{i+1})$, 
for some $i\ge 1$, it satisfies $\dot x_i=-v(x_i,t,\hat\rho_{i-1},\rho_i)$
with $\hat \rho_i(t)=\phi_{x_{i-1}(t)}^{x_i(t)}\big(\rho_{i-1}(t);t\big)$.
\item[$\bold{(2)}$] We have
 $\dot\rho_0=\b(x_0,t,\rho_0),$ and for $i>0$,
we have
$\dot\rho_i=-K(x_i,t,\rho_i,\hat\rho_{i-1}).$
\item[$\bold{(3)}$] With rate $f^2\big(a_+,t,\hat\rho_n,\rho_{n+1}),$
the configuration $\bold{q}$ gains a new particle 
$(x_{n+1},\rho_{n+1})$, with $x_{n+1}=a_+$.
This new configuration is denoted by $\bold q(\rho_{n+1})$.
\item[$\bold{(4)}$] When $x_1$ reaches $a_-$, we relabel 
particles $(x_i,\rho_i), \ i\ge 1$, as $(x_{i-1},\rho_{i-1})$.
\item[$\bold{(5)}$] When $x_{i+1}-x_i$ becomes $0$ for some $i\ge 1$, then 
$\bold{ q}(t)$ becomes
$\bold{ q}^i(t)$, that is obtained from $\bold{ q}(t)$ by omitting
$(\rho_i,x_i)$ and relabeling particles to the right of the $i$-th particle.
\ei
\qed

\bs
As we mentioned before,  the function $\rho(x,t;\bq(t))$ is indeed 
an entropic solution of \eqref{eq1.2}. We also need a
stability inequality for our constructed solutions.

 \bp\la{pro2.1} {\bf(i)} The function $\rho(x,t)=\rho(x,t;\bq(t))$,
 with $\bq(t)$ evolving as above, is an entropy solution of 
$\rho_t=H(x,t,\rho)_x$ in $(a_-,a_+)\x (t_0,T)$.

\ms\noi
{\bf(ii)}
The process $t\mapsto m(t):=\rho(a_+,t)=\rho(a_+,t;\bq(t))$ 
is a Markov process
with generator  $\cA^2_{a_+,t}$.

\ms\noi
{\bf(iii)}  Suppose  $\rho,\rho':[a_-,a_+]\x [t_0,T)
\to [P_-,P_+]$ are two piecewise $C^1$ entropy solutions of 
$\rho_t=H(x,t,\rho)_x$. If $t_0\le s\le t<T$, then
 \begin{align}\nonumber
\int_{a_-}^{a_+}|\rho'(x,t)-\rho(x,t)|\ dx\le &e^{C_0(t-s)}
\int_{a_-}^{a^+}|\rho'(x,s)-\rho(x,s)|\ dx\\ \la{eq2.4}
&+e^{C_0(t-s)}\int_{s}^t
\big|H(a_+,\th,\rho'(a_+,\th))-H(a_+,\th,\rho(a_+,\th))\big|\ d\th,
\end{align}
where 
\[
C_0=\max_{x\in [a_-,a_+]}\max_{t\in[t_0,T]}\max_{\rho\in[P_-,P_+]}
|H_{x\rho}(x,t,\rho)|.
\]
 \ep

\bs\noi
{\bf Remark 2.1} From Proposition 2.1{\bf(iii)} we learn that two solutions
$\rho$ and $\rho'$ equal in $[a_-,a_+]\x[t_0,T]$
if they coincide at $t=0$, and $x=a_+$. This confirms the fact
that under the assumption $H_\rho<0$, the boundary $a_-$ is free.
In particular, $\rho(x,t;\bq(t))$ is the unique solution which satisfies the stochastic boundary condition at $x=a_+$.
\qed

\ms
The proof of Proposition~2.1 will be given at the end of this section.
We continue with a precise description for
the PDMP $\rho(\cdot,t)$ in terms of $\bq(t)$ and some preparatory
steps toward the proof of Proposition 2.1 and Theorem 2.1.

\ms\noi
{\bf Definition 2.2(i)} To ease the notation, we write
$\l(x,t,\rho )$ and $A(x,t,\rho)$, for $(Af^1)(x,t,\rho )$
and $(Af^2)(x,t,\rho)$ respectively. Given $\bq\in\D_n$, we also set
\begin{align*}
\G(x,y,t,\rho)&=\int_{x}^y\l(z,t,\phi_x^z(\rho;t)) \ dz,\\
 \G(\bq,t)&
=\int_{a_-}^{a_+}\l\big(y,t,\rho\big(y,t;\bold{q}\big)\big)\
 dy=\sum_{i=0}^{n}\G(x_{i},x_{i+1},t,\rho_{i}).
 \end{align*}
 
 \ms\noi
 {\bf(ii)} We define a measure $\mu(d\bq,t)$ on the set $\D$ that is our candidate for the law of $\bq(t)$. 
 The restriction of
$\mu$ to $\D_n$ is denoted by $\mu^n(d\bq,t)$. This measure is
explicitly  given by 
\[
\ell(t, d\rho_0)\exp\left\{-\G(\bq,t)\right\}\prod_{i=1}^{n}
\ f\big(x_i,t,\phi_{x_{i-1}}^{x_{i}}(\rho_{i-1};t),\rho_{i})\  dx_{i}d\rho_{i},
\]
where $f$ solves \eqref{eq1.7} and $\ell$ solves \eqref{eq1.18}.
Note that if $\rho(x,t)=R_t(\bq(t))(x)$, with $R$ as in Definition {\bf{2.1(ii)}}, then the process $x\mapsto \rho(x,t), x\ge a_-$ is a Markov process
associated with the generator $\cA_{x,t}^1$, and an initial law $\ell(t,\cdot)$.

\ms\noi
{\bf(iii)}
Let us write $T_x^yg(\rho)=g(\phi_x^y(\rho;t))$ and 
$(\cD_x g)(\rho)=b(x,t,\rho )g'(\rho)$ for its generator
(to simplify the notation, we do not display the dependence of $T_x^y$ and
$\cD_x$ on $t$). It is straightforward to show
\be\la{eq2.5}
T_x^y\circ T_y^z=T_x^z,\ \ \ \ \ \ 
\frac {dT_x^y}{dy}=T_x^y\circ\cD_y ,\ \ \ \ \ \ 
\frac {dT_x^y}{dx}=-\cD_x\circ T_x^y .
\ee
Indeed
\begin{align*}
T_{x}^{y+\d}g&=T_x^y\big(g\circ
\phi_{y}^{y+\d}\big)=T_x^y\big(g+\d\cD_y g+o(\d)\big)
=T_x^yg+\d \big(T_x^y\circ\cD_y\big) g+o(\d ),\\
T_{x-\d }^yg&=T_{x-\d }^x\big(T_x^y g\big)=\big(T_x^y g\big)\circ\phi_{x-\d }^x=
T_x^y g-\d\big( \cD_x\circ T_x^y\big)g+o(\d ).
\end{align*}
\qed 

\bs
In the following Lemma, we derive several identities that we will use
 for the proof of Proposition~2.1 and Theorem~2.1. 

\bl\la{lem2.1} The following identities are true: 
\begin{align}
\la{eq2.6}
&b(x,t,\rho )\G_{\rho}(x,y,t,\rho)=-\G_x(x,y,t,\rho)+\l(x,t,\rho )=
-\int_x^{y}\big[\l\big(z,t,\phi_x^z(\rho;t)\big)\big]_x\ dz,\\
\la{eq2.7}
&b_t=\b_x+b\b_{\rho}-b_{\rho}\b,\\
\la{eq2.8}
&\big[\phi_x^y(\rho;t)\big]_t=\b(y,t,\phi_x^y(\rho;t))-
\b\big(x,t,\rho )
\big[\phi_x^y(\rho;t)\big]_{\rho},\\
\la{eq2.9}
&\l_t(x,t,\rho )+\b(x,t,\rho )\l_\rho(x,t,\rho )
=b(x,t,\rho )A_\rho(x,t,\rho )+A_x(x,t,\rho ),\\
\la{eq2.10}
&\G_t(x,y,t,\rho)+\b(x,t,\rho )
\G_{\rho}(x,y,t,\rho)=A\big(y,t,\phi_x^y(\rho;t)\big)-A(x,t,\rho ),\\
\la{eq2.11}
&\big[\phi_x^y(\rho;t)\big]_x+b(x,t,\rho)\big[\phi_x^y(\rho;t)\big]_{\rho}=0..
\end{align}
\el

\noi\ms
{\bf Proof} 
For the proof of \eqref{eq2.6} use  the definition of $\G$ and 
\eqref{eq2.5} to assert that
the left-hand side of \eqref{eq2.6} equals to
\begin{align*}
\int_x^{y}b(x,t,\rho )
\left[\l\big(z,t,\phi_x^z(\rho;t)\big)\right]_{\rho}\ dz=&-\int_x^{y}
\big[\l\big(z,t,\phi_x^z(\rho;t)\big)\big]_x\ dz\\
&=-\left[\int_x^{y}
\l\big(z,t,\phi_x^z(\rho;t)\big)\ dz\right]_x+\l(x,t,\rho )\\
&=-\G_x(x,y,t,\rho)+\l(x,t,\rho ).
\end{align*}

\ms 
For \eqref{eq2.7} observe that by \eqref{eq1.6},
\begin{align*}
\b_x+b\b_{\rho}-b_{\rho}\b&=\big(b H_\rho+H_x\big)_x+b
\big(b H_\rho+H_x\big)_{\rho}-bb_{\rho} H_\rho-b_{\rho}H_x\\
&=b H_{\rho x}+ b_xH_\rho+H_{xx}
+bb_\rho H_\rho+b^2 H_{\rho\rho}+b
H_{\rho x}-bb_{\rho} H_\rho-b_{\rho}H_x\\
&=2b H_{\rho x}+b_xH_\rho +H_{xx}
+b^2 H_{\rho\rho}-b_{\rho}H_x=b_t.
\end{align*}

\ms
 We now turn to the proof of \eqref{eq2.8}. Set
\[
X(x,y,t,\rho):=\big[\phi_x^y(\rho;t)\big]_t-\b(y,t,\phi_x^y(\rho;t))
+\b\big(x,t,\rho )
\big[\phi_x^y(\rho;t)\big]_{\rho}.
\]
We wish to show  that $X(\rho,x,y,t)=0$ for all $(x,y,t,\rho)$.
This is trivially true when $x=y$.  On the other hand,
\begin{align*}
X_y(x,y,t,\rho)=&\big[b(y,t,\phi_x^y(\rho;t))\big]_t-
(\b_y+b\b_\rho) (y,t,\phi_x^y(\rho;t))
   +\b\big(x,t,\rho )\big[b(y,t,\phi_x^y(\rho;t))\big]_\rho\\
=&b_t(y,t,\phi_x^y(\rho;t))+b_\rho(y,t,\phi_x^y(\rho;t))
\big[\phi_x^y(\rho;t)\big]_t-(\b_y+b\b_\rho) (y,t,\phi_x^y(\rho;t))
\\
& +\b\big(x,t,\rho )b_\rho(y,t,\phi_x^y(\rho;t))
\big[\phi_x^y(\rho;t)\big]_{\rho}\\
=&b_\rho(y,t,\phi_x^y({\rho};t))X(x,y,t,{\rho}),
\end{align*}
where we used \eqref{eq2.7} for the third equality. As a result.
\[
X(x,y,t,{\rho})=X(x,x,t,{\rho})\exp
\left[\int_x^yb_\rho\big(z,t,\phi_x^z
({\rho};t)\big)\ dz\right]=0.
\]
This completes the proof of \eqref{eq2.8}.

\ms
For \eqref{eq2.9}, we first observe
\be\la{eq2.12}
\int_\rho^\i Q(f)(x,t,\rho,\rho_+)\ d\rho_+=0,
\ee
because the left-hand side equals
\begin{align*}
&\int_{\rho}^\i \int 1\!\!1\big(\rho_*\in(\rho, \rho_+)\big)
\big(v(x,t,\rho_*,\rho_+)-v(x,t,\rho,\rho_*)\big)
f(x,t,\rho,\rho_*)f(x,t,\rho_*,\rho_+)
\ d\rho_*d\rho_+\\
& \ \     -\int_{\rho}^\i\left[A(x,t,\rho_+)-A(x,t,\rho )-
v(x,t,\rho,\rho_+)\big(\l(x,t,\rho_+)-\l(x,t,\rho )\big)\right]
f(x,t,\rho,\rho_+)\ d\rho_+\\
&\ =\int_{\rho}^\i 
\big(A(x,t,\rho_*)-v(x,t,\rho,\rho_*)\l(x,t,\rho_*)\big)
f(x,t,\rho,\rho_*)\ d\rho_*\\
&\ \ \    -\int_{\rho}^\i\left[A(x,t,\rho_+)-A(x,t,\rho )-
v(x,t,\rho,\rho_+)\big(\l(x,t,\rho_+)-\l(x,t,\rho )\big)\right]
f(x,t,\rho,\rho_+)\ d\rho_+\\
&=\int_{\rho}^\i\big(A(x,t,\rho )-
v(x,t,\rho,\rho_+)\l(x,t,\rho )\big)
f(x,t,\rho,\rho_+)\ d\rho_+=0.
\end{align*}
We next integrate both sides of \eqref{eq1.7} with respect to
$\rho_+$ and use \eqref{eq2.12} to assert
\begin{align}
\l_t(x,t,\rho )&=\int_\rho^\i \left\{
(Cf)(x,t,\rho,\rho_+)+ \big[(vf)(x,t,\rho,\rho_+)\big]_x\right\}\ d\rho_+\nonumber\\
&=\int_\rho^\i (Cf)(x,t,\rho,\rho_+) \ d\rho_++ A_x(x,t,\rho ).
\la{eq2.13}
\end{align}
On the other hand, 
\begin{align*}
\int \big(C^-f\big)(x,t,\rho,\rho_+) \ d\rho_+=&
b(x,t,\rho )\int_{\rho}^\i 
(vf)_{\rho}(x,t,\rho,\rho_+)\ d\rho_+-\b(x,t,\rho )\int_{\rho}^\i  f_{\rho}(x,t,\rho,\rho_+)\ d\rho_+\\
=&(b A_\rho-\b \l_\rho)(x,t,\rho )
+(bH_\rho-\b)(x,t,\rho)f(x,t,\rho,\rho),\\
\int \big(C^+f\big)(x,t,\rho,\rho_+) \ d\rho_+=&
\int_{\rho}^\i 
\big[K(x,t,\rho_+,\rho)f(x,t,\rho,\rho_+)\big]_{\rho_+}\ d\rho_+\\
=&(\b-bH_\rho)(x,t,\rho)f(x,t,\rho,\rho),
\end{align*}
because $f(x,t,\rho,\i)=0$, and
\[
\lim_{\rho_+\to\rho}v(x,t,\rho,\rho_+)= H_\rho(x,t,\rho).
\]
From this, and \eqref{eq2.13} we deduce \eqref{eq2.9}.

\ms
We now turn to the proof of \eqref{eq2.10}.
We rewrite \eqref{eq2.10} as
\begin{align*}
&\int_x^y\big[\l(z,t,\phi_x^z(\rho;t))\big]_t\ dz+
\b(x,t,\rho )\int_x^y\big[\l(z,t,\phi_x^z(\rho;t))\big]_{\rho}\ dz
=\int_x^y \big[A\big(z,t,\phi_x^z(\rho;t)\big)\big]_z\ dz.
\end{align*}
For this, it suffices to check
\be\la{eq2.14}
\big[\l(z,t,\phi_x^z(\rho;t))\big]_t+
\b(x,t,\rho )\ \big[\l(z,t,\phi_x^z(\rho;t))\big]_{\rho}
= \big[A\big(z,t,\phi_x^z(\rho;t)\big)\big]_z,
\ee
for every $(x,z,t,\rho)$. By \eqref{eq2.8}, the identity \eqref{eq2.14} is equivalent to 
\[
\l_t(z,t,\phi_x^z(\rho;t))+\b(z,t,\phi_x^z(\rho;t))\l_{\rho}(z,t,\phi_x^z(\rho;t))
= \big[A\big(z,t,\phi_x^z(\rho;t)\big)\big]_z.
\]
 This  is an immediate consequence of \eqref{eq2.9} because
\[
\big[A\big(z,t,\phi_x^z(\rho;t)\big)\big]_z=b(z,t,\phi_x^z(\rho;t))A_{\rho}\big(z,t,\phi_x^z(\rho;t)\big)
+A_{z}\big(z,t,\phi_x^z(\rho;t)\big).
\]
The proof of \eqref{eq2.10} is complete.

\ms
Finally, \eqref{eq2.11} is simply the third equation of \eqref{eq2.5} applied to
the function $g(\rho)=\rho$. 
\qed

\bs
We are now ready to establish Proposition~2.1.

\ms\noi
{\bf{Proof of Proposition 2.1(i)}} We first show that $\rho$ solves
\eqref{eq1.2} classically away from the shock curves. For this, 
take a point $(x,t)$ such that $x\in\big(x_i(t),x_{i+1}(t)\big)$,
for some nonnegative integer $i$.
Let us write $\hat \phi_x^y(\rho;t)$ and 
$\tilde \phi_x^y(\rho;t)$ for the partial derivatives
$\big[\phi_x^y(\rho;t)\big]_\rho$ and 
$\big[\phi_x^y(\rho;t)\big]_x$ respectively. From
$\rho(x,t)=\phi_{x_i(t)}^x(\rho_i(t);t)$, we learn
\begin{align}\nonumber
\rho_t(x,t)=&-v\big(x_i(t),t,\hat\rho_{i-1}(t),\rho_i(t)\big)\tilde \phi_{x_i(t)}^x(\rho_i(t);t)\\
\nonumber
&+\b(x,t,\rho(x,t))-\b\big(x_i(t),t,\rho_i(t) )\ 
\hat\phi_{x_i(t)}^x(\rho_i(t);t)\\ \la{eq2.15}
&-K\big(x_i(t),t,\rho_i(t),\hat\rho_{i-1}(t)\big)
\hat\phi_{x_i(t)}^x(\rho_i(t);t)\\ \nonumber
=&-v\big(x_i(t),t,\hat\rho_{i-1}(t),\rho_i(t)\big)
\tilde \phi_{x_i(t)}^x(\rho_i(t);t)
+\b(x,t,\rho(x,t))\\ \nonumber
&-v\big(x_i(t),t,\hat\rho_{i-1}(t),\rho_i(t)\big)b\big(x_i(t),t,\rho_i(t) )
\hat\phi_{x_i(t)}^x(\rho_i(t);t)\\ \nonumber
=&\b(x,t,\rho(x,t))=H(x,t,\rho(x,t))_x,
\end{align}
as desired. Here we used \eqref{eq2.8} for the first equality, and \eqref{eq2.11} for the third equality. 
Since the Rankine-Hugoniot 
Formula is valid at shock curves by our construction, and \eqref{eq1.1} holds classically 
away from the shock curves, we deduce that $\rho$ is a weak 
solution of \eqref{eq1.1}. On the other hand, since 
$\rho(x_i(t)-,t)<\rho(x_i(t)+,t)$
by construction, we deduce that $\rho$ is an entropy solution
in $(a_-,a_+)\x(t_0,T)$. 

\ms\noi
{\bf(ii)} From the way the boundary dynamics is described in {\bf(3)}, the process $m(t)$ 
 depends on the particle system to the left of $a_+$.
Nonetheless we show that if the process $\bar m(t)$ is a Markov process  with generator $\cA_{a_+,t}^2$, and initial state $m_0=
\rho^0(a_+)$, then $m(t)=\bar m(t)$. 
To verify this, let us construct the process 
$t\mapsto \bar m(t)$ with the aid of a sequence
of independent standard exponential random variables $\big(\tau_i:\ i\in\bN\big)$. Let us write 
$\g_s^t(\rho)$ for the flow of the ODE associated with speed 
$\b(a_+,t,\rho)$, and define
\be\la{eq2.16}
\eta(t,m)=\int_m^\i f^2(a_+,t,m,\rho_+)\ d\rho_+.
\ee
Now construct a sequence $\bz=\big((\s_i,m_i):\ i=0,1,\dots\big)$ inductively by the following recipe: $\s_0=0$, and
 given $(\s_i,m_i)$, we set 
\begin{align*}
\s_{i+1}&=\min\left\{s>\s_i:\ \int_{\s_i}^s 
\eta\big(\th,\g_{\s_i}^\th(m_i)\big)\ d\th\ge \tau_{i+1}\right\},\ \ \ \ 
\hat m_{i}=\g_{\s_i}^{\s_{i+1}}(m_i),
\end{align*}
and select $m_{i+1}$ randomly according to the probability measure
\[
\eta\big(\s_{i+1},\hat m_i\big)^{-1}\ 
f^2\big(a_+,\s_{i+1},\hat m_i, m_{i+1}\big)\
dm_{i+1}.
\]
Using our sequence $\bz$, we construct $\bar m(t)$ by
\[
\bar m(t)=\sum_{i=0}^\i \g_{\s_i}^t( m_i)1\!\!1\big(t\in[\s_i,
\s_{i+1})\big).
\]
By the very construction of the processes $m(t)$ and $\bar m(t)$,
the desired equality $m(t)=\bar m(t),\ t\ge t_0$ would follow if we can show  that $m(t)=\bar m(t)$ for $t\in (\s_{i-1},\s_{i})$ for every $i\in\bN$.
This can be checked by induction on $i$.  If
 there are exactly $n$ particle to the left of $a_+$, and 
we already know that $\hat \rho_n(\s_i)= \hat m_{i-1}$, then we can guarantee that
$\rho_{n+1}(\s_i)=m_i$. Moreover, $\hat \rho_{n+1}(t)=\g_{\s_i}^t(m_i)$
for $t\in(\s_i,\s_{i+1})$, because  the function
$\zeta(t)=\phi_{x_n(t)}^{a_+}( \rho_n(t);t)$ satisfies
\be\la{eq2.17}
\dot \zeta (t)=\b(a_+,t,\zeta(t)),
\ee
by
\eqref{eq2.15} in the case of $x=a_+$.
This completes the proof of part {\bf(ii)}.

\ms\noi
{\bf(iii)} The proof of \eqref{eq2.4} is a standard 
application of the celebrated Kruzhkov's inequality
[K], and we only sketch it. It is not hard to show that the piecewise $C^1$ entropy solutions $\rho$ and $\rho'$
can be extended to entropy solutions that are defined on 
a larger domain $(b_-,b_+)\x[t_0,T)$, with
$b_-<a_-<a_+<b_+$. With a slight abuse of notation, we write $\rho$ and $\rho'$ for these extensions.
Given an arbitrary constant $c$, 
the following Kruzkov's entropy 
inequalities hold weakly in $(b_-,b_+)\x[t_0,T)$:
\begin{align*}
|\rho(x,t)-c|_t&\le
|H(x,t,\rho(x,t))-H(x,t,c)|_x+\ sgn\big(\rho(x,t)-c\big)
H_x(x,t,c),\\
|\rho'(x,t)-c|_t&\le
|H(x,t,\rho'(x,t))-H(x,t,c)|_x+\ sgn\big(\rho'(x,t)-c\big)
H_x(x,t,c).
\end{align*}

This allows us to use Kruzkov's standard arguments as in [K] to deduce 
\begin{align*}
|\rho(x,t)-\rho'(x,t)|_t\le &
|H(x,t,\rho(x,t))-H(x,t,\rho'(x,t))|_x\\
&-\ sgn\big(\rho(x,t)-\rho'(x,t)\big)
\big(H_x(x,t,\rho(x,t))-H_x(x,t,\rho'(x,t))\big)\\
\le &|H(x,t,\rho(x,t))-H(x,t,\rho'(x,t))|_x+C_0|\rho(x,t)-\rho'(x,t)|,
\end{align*}
weakly in $(b_-,b_+)\x[t_0,T)$. From this, we can readily deduce
\begin{align}
\left[e^{-C_0t}|\rho(x,t)-\rho'(x,t)|\right]_t\le &
e^{-C_0t}|H(x,t,\rho(x,t))-H(x,t,\rho'(x,t))|_x,\la{eq2.18}
\end{align}
weakly in $(b_-,b_+)\x[0,T)$. We wish to integrate both sides of \eqref{eq2.18} with respect to $x$ from $a_-$ to $a_+$. To perform
such integration rigorously, we take a smooth function
$\g$ of compact support
with $\int\g\ dx=1$,
rescale it as $\g_\e(x)=\e^{-1}\g(x/\e)$, and choose 
$\tau_\e(x)$ so that $\tau_\e\ge 0$,
$\tau'_\e(x)=\g_\e(x-a_-)-\g_\e(x-a_+)$,
and $\tau_\e(b_-)=0$. 
For small $\e$, the function $\tau_\e$ is supported in 
$(b_-,b_+)$. We can
now integrate both sides of \eqref{eq2.18} against $\tau_\e$ to deduce
that weakly
\begin{align*}
&\left[e^{-C_0t}\int |\rho(x,t)-\rho'(x,t)|\tau_\e(x)\ dx\right]_t\le
-\int e^{-C_0t}
|H(x,t,\rho(x,t))-H(x,t,\rho'(x,t))|\ \tau'_\e(x)\  dx.
\end{align*}
We can now send $\e$ to $0$ to arrive at 
\begin{align*}
\left[e^{-C_0t}\int_{a_-}^{a_+} |\rho(x,t)-\rho'(x,t)|\ dx\right]_t\le &
 e^{-C_0t}
|H(a_+,t,\rho(a_+,t))-H(a_+,t,\rho'(a_+,t))| \\
&-e^{-C_0t}
|H(a_-,t,\rho(a_-,t))-H(a_-,t,\rho'(a_-,t))| .
\end{align*}
Integrating both sides over the time interval $[s,t]$ yields
\begin{align*}
e^{-C_0t}\int_{a_-}^{a_+} |\rho(x,t)-\rho'(x,t)|\ dx\le &
e^{-C_0s}\int_{a_-}^{a_+} |\rho(x,s)-\rho'(x,s)|\ dx\\
&+ \int_s^te^{-C_0\th}
|H(a_+,\th,\rho(a_+,\th))-H(a_+,\th,\rho'(a_+,\th))|\ d\th \\
&-\int_s^te^{-C_0\th}
|H(a_-,\th,\rho(a_-,\th))-H(a_-,\th,\rho'(a_-,\th))|\ d\th \\
\le &
e^{-C_0s}\int_{a_-}^{a_+} |\rho(x,s)-\rho'(x,s)|\ dx\\
&+ \int_s^te^{-C_0\th}
|H(x,\th,\rho(a_+,\th))-H(x,\th,\rho'(a_+,\th))|\ d\th .
\end{align*}
This evidently implies \eqref{eq2.4}.
\qed

 \section{Forward Equation}
  
 As a preliminary step for establishing Theorem~2.1,
we derive a Kolmogorov type forward equation for the measure
$\mu(d\bq,t)$.
We first introduce some  notation for the particle dynamics.

 \ms\noi
 {\bf Definition 3.1(i)}
  For $0 \leq s \leq t$ and $\bq \in \D$, we write 
    $\psi_s^t \bq$
    for the deterministic evolution from time $s$ to $t$ of the
    configuration $\bq$ according to the annihilating particle dynamics
    of Definition 2.1{\bf(iii)},
 \emph{without} random entry dynamics at $x =a_+$.
 
 \ms\noi
 {\bf(ii)} Given a configuration $\bq =
   \big((x_0,\rho_0),\dots,(x_n,\rho_n)\big)$ and 
   $\rho_+ \in\bR$,
    write $\epsilon_{\rho_+} \bq$ for the configuration
    $\big((x_0,\rho_0),\dots,(x_n,\rho_n),(a_+,\rho_+)\big)$.
 
  \ms\noi
 {\bf(iii)} Write $\Psi_s^t \bq$ for the \emph{random} evolution of
    the configuration according to deterministic particle dynamics
    interrupted with random entries at $x = a_+$ according to the
    boundary process as in {\bf(3)} in Definition 2.1{\bf(iii)}, where the
    latter has been started at time $s$ with value 
    $\phi_{x_n}^{a_+}(\rho_n;s)$.  In
    particular, if the jumps  between times $s$ and $t$
    occur at times $ \tau_1 < \cdots < \tau_k$ with values
    $m_1, \cdots , m_k$, then
    \begin{equation}\la{eq3.1}
      \Psi_s^t \bq = \psi_{\tau_k}^t \epsilon_{m_k}
      \psi_{\tau_{k-1}}^{\tau_k} \epsilon_{m_{k-1}} \cdots
      \psi_{\tau_1}^{\tau_2} \epsilon_{m_1} 
      \psi_s^{\tau_1} \bq.
    \end{equation}
    
    \ms\noi
   {\bf(iv)} For $n\ge 1$, and $i\in\{0,\dots,n-1\}$,
   we write $\p_i\D_n$ for the portion of the boundary $\D_n$
 such that $x_i=x_{i+1}$. Note that $\bq(t)$ reaches the boundary set   $\p_0\D_n$ at time $\tau$ if at this time $x_1(\tau)=a_-$.
   For time $t$ immediately after $\tau$, the configuration $\bq(t)$
   belongs to $\D_{n-1}$ with $\rho_0$ taking new value. 
Similarly  $\bq(t)$ reaches the boundary set   $\p_i\D_n$ 
for some $i>0$ at time $\tau$ if at this time 
   $x_{i+1}$ collides with $x_{i}$.
   For time $t$ immediately after $\tau$, the configuration $\bq(t)$
   belongs to $\D_{n-1}$. We also set 
\[
\hat\p\D_n=\cup_{i=0}^n\p_i\D_n.
\]

   \ms\noi
   {\bf(v)} We write $\p_{n+1}\D_{n+1}$ for the set of
   points $\bq\in\D_{n+1}$ with $x_{n+1}=a_+$.
   When $\bq\in   \D_n$, and a new particle is created at $a_+$
 at time $\tau$
   by the stochastic boundary dynamics,  
the configuration $\bq(\tau+)$ is regarded as a boundary point in $\p_{n+1}\D_{n+1}$.

 \ms\noi
   {\bf(vi)}
Given a function $G:\D\to\bR$, we  write $G^n$ for the restriction of the function $G$ to the set $\D_n$.
Also, given a measure on $\D$, we write $\nu^n$ for the restriction of a measure $\nu$ to $\D_n$.

 \ms\noi
   {\bf(vii)} We write $\cL=\cL^t$ for the generator of the 
(inhomogeneoys Markov)
process $\bq(t)$. This generator can be expressed 
as $\cL=\cL_0+\cL_b$, where $\cL_0$ is the generator of the deterministic part of dynamics, and $\cL_b$ 
represents the Markovian boundary dynamics. The deterministic and stochastic dynamics restricted to 
 $\D_n$ have generators that are denoted by $\cL_{0n}$
and $\cL_{bn}$ respectively. 
While $\bq(t)$ remains in $\D_n$, its evolution is
governed by an ODE of the form
\[
\frac{d\bq}{dt}(t)=\bb\big(\bq(t),t\big),
\]
with $\bb=\bb_n:\D_n\to\bR^{2n+1}$, that can be easily described with the aid of rules {\bf(1)} and
{\bf(2)} of Definition~2.1{{\bf(iii)}, and \eqref{eq2.3}. 
Given this vector field, the generator $\cL_{0n}$ is given by
\[
\cL_{0n}F=\bb\cdot \nabla F,
\]
where $\nabla F$ is the full gradient of $F$ with respect to variables 
$\big(\rho_0,x_1,\rho_1,\dots,x_n,\rho_n\big)$. We also write 
$\cL_{0n}^{t*}=\cL_{0n}^*$ for the adjoint of $\cL_{0n}$
with respect to the Lebesgue measure:
\[
\cL_{0n}^*\mu=-\nabla\cdot(\mu\bb).
\]
   \qed
 
\bs 
It is not hard to show that  $t\mapsto\Psi_s^t \bq,\ t\ge s$ is indeed
a strong Markov process. This is rather straightforward, and we 
refer to Davis [D] for details.

 We establish Theorem 2.1 by verifying the forward equation
$\dot \mu=\cL^*\mu$, or equivalently
\be\la{eq3.2}
\dot\mu^n= \big(\cL^*\mu\big)^n,
\ee
for all $n\ge 0$,
where $\mu$ was defined in Definition 2.2{\bf(ii)}, and $\cL^*$ is the adjoint of the operator $\cL$. To explain this,
observe that Theorem~2.1 offers a candidate for the law of $\bq(t)$, namely the measure $\mu(d\bq,t)$. Hence
for our Theorem~2.1, it suffices to  show 
  \be\la{eq3.3}
   \int G\big(\bq,t\big)\ \mu\big(d\bq,t\big)
  =\bE\int G\big(\Psi_0^t\bq,t\big)\ \mu(d\bq,0\big),
    \ee
 for every function $G$ of the form
\be\la{eq3.4}
G(\bq,t)=\exp\left(i\int_{a_-}^{a_+}\rho(x,t;\bq)\var(x)\ dx\right),
\ee
 for some smooth function $\var$ (we refer to the beginning 
of Section 3 of [KR1] for more details.) Here and below, we write
$\bP$ and $\bE$ for the probability and the expected value for the randomness associated with the boundary dynamics. 
To ease the notation, we set $\hat G( \bq,s)=\bE \ G(\Psi_s^t \bq,t)$.
We establish
\eqref{eq3.3} by verifying 
\begin{equation}\la{eq3.5}
  \frac{d}{ds} \int \hat G( \bq,s)\ \mu(d\bq,s) = 0,
\end{equation}
for $t_0 < s < t$.
 The differentiation of $\mu(d\bq,s)$ can be carried out directly and
 poses no difficulty. As for the contribution of $G(\Psi_s^t \bq,t)$ to the $s$-derivative, we wish to show
 \be\la{eq3.6}
 \int \hat G_s( \bq,s)\ \mu(d\bq,s) =-\int \big(\cL^s  \hat G\big)( \bq,s)\ \mu(d\bq,s) .
\end{equation}
Since the deterministic part of the evolution is discontinuous in time, 
the justification of \eqref{eq3.6} requires some work. Additionally,
to make sense of the right-hand side, we need 
$\hat G$ to be in the domain of the definition of $\cL_0^s$.
We expect $\hat G$ to be weakly differentiable with respect to $\bq$.
To avoid the differentiability question of $\hat G$, we 
would formally apply an integration by parts to the right-hand side of 
\eqref{eq3.6}, so that the differentiation operator 
would act on the density of $\mu$, which is differentiable.
We also have a boundary contribution that correspond to the collisions
between particles.
We establish the following
 variant of the forward equation \eqref{eq3.6}.

\bth\la{th3.1} We have
\begin{align}\nonumber
\lim_{s'\uparrow s}
\int_{\D_n}
\frac{\hat G^n(\bq,s')-\hat G^n(\bq,s)}{s-s'}\ \mu^n(\bq,s)\ d\bq
=&\int_{\D_n} (\cL_{b}^s \hat G)^n(\bq,s)\mu^n(\bq,s)\ d\bq\\
&+\int_{\D_n}  \hat G^n(\bq,s)\big(\cL_{0n}^{s*}
\mu^n\big)(\bq,s)\ d\bq\la{eq3.7}\\
&+\int_{\hat\p \D_{n+1}}  \hat G^{n+1}(\bq,s)\mu^{n+1}(\bq,s)(\bb_{n+1}\cdot {\bold{N}}_{n+1})\ \s(d\bq),
\nonumber
\end{align}
where ${\bold{N}}_{n+1}$ is the outer unit normal vector of $\hat\p \D_{n+1}$, and
$\s(d\bq)$ is the surface measure of $\p \D_{n+1}$.
\et

\bs
Note that for the differentiation in \eqref{eq3.7} we will need to compare
$\hat G(\bq,s)$ and $\hat G( \bq,s')$ 
for $t_0 < s' < s\le t$. 
As a warm-up we verify the Lipschitzness of the function $s\mapsto \bE\ \hat G(\bq,s)$.

\begin{lemma}
  \label{lem3.1} Fix $t>t_0$. 
  There exists a constant $C_1=C_1(\var,H,f)$
  such that 
  \begin{equation}
    \label{eq3.8}
  \big|\hat G(\bq,s') -\hat G(\bq,s)\big| \le C_1(n+1)|s'-s|,
  \end{equation}
 for all $\bq\in\D_n$ and $s,s'\in[t_0,t]$.
\end{lemma}

The proof follows from the $L^1$-stability \eqref{eq2.4} and a coupling argument for the stochastic boundary dynamics. We skip the proof of Lemma 3.1 because it is very similar to the proof of the analogous
Lemma 3.1 that appeared in [KR2].
Armed with \eqref{eq3.8}, we are now ready for the proof of \eqref{eq3.7}.

\ms\noi
{\bf Proof of Theorem 3.1} {\em(Step 1)}
Let $t_0 < s' < s\le t$.
We first show that we can separate
the deterministic and stochastic portions of the dynamics over the
time interval $[s',s]$, when $s-s'$ is small.
   Write $\tau=\tau(\bq,s')$ for the first time a jump occurs at $x=a_+$ after the time $s'$, and let $E$ denote
 the event that  $\tau\in (s',s)$. We also write
\[
\hat\rho_n=R_s(\bq)(a_+)=\phi_{x_n}^{a_+}(\rho_n;s),\ \ \ \ 
\hat\rho'_n=R_{s'}(\bq)(a_+)=\phi_{x_n}^{a_+}(\rho_n;s').
\]
By the Lipschitz regularity of $b$, 
we can show that $\hat\rho_n-\hat\rho'_n=O(s-s')$
(see also \eqref{eq3.14} below).
Recall that $\g$ denotes the  flow associated with the ODE 
\eqref{eq2.17}, and $\eta$ was defined in \eqref{eq2.16}. Observe that by the Lipschitz regularity of $\eta$ (which is the consequence of the Lipschitz regularity of $v$ and $f$),
  \be\la{eq3.9}
    \bP\big(E\big) =
     \int_{s'}^s\eta\big(\th,\g_{s'}^\th(\th,\hat\rho'_n)\big)
    \ d\th+O((s-s')^2)= 
       (s-s') \eta\big(s,\hat\rho_n\big)  +    O((s-s')^2),
         \ee
  with both errors bounded uniformly over $\bq$. 
We claim that there exists  a constant $c_1$
   so that for $\bq \in \D_n$,
  \begin{align}\nonumber
   \hat G(\bq,s')=&    (s-s')    \int_{\hat\rho_n}^\i     
    \left(\bE \left[\hat G\big( \epsilon_{\rho_+} \psi_{s'}^\tau \bq,s\big) 
    \ \big|\ E \right]- 
    \hat G(\psi^s_{s'}\bq,s)\right)      
 f^2(a_+,s,\hat\rho_n,\rho_+)\ d\rho_+\\
      &\ \ \ \ \ \ \  \ \ \ \ +\hat G\big(\psi^s_{s'}\bq,s\big)+
(s-s')^2R(\bq,s',s),\la{eq3.10}
      \end{align}
     with $\big|R(\bq,s',s)\big|\le c_1 (n+1)$.  To prove \eqref{eq3.10},
first observe that by the  Markov property of the random flow $\Psi$, 
  \[
    \hat G(\bq,s') = \bE\  G\big(\Psi_{s}^t \Psi^s_{s'}\bq,t\big) =
     \bE    \ \hat G\big( \Psi^s_{s'}\bq,s\big).
  \]
 On $E^c$ (the complement of $E$),  we see only 
  the deterministic flow $\psi$ over the time interval $(s',s)$:
  \begin{align}\nonumber
   \bE\  \hat G\big(\Psi_{s'}^s\bq,s\big) 1\!\!1_{E^c} 
    &= \hat G\big(\psi_{s'}^s\bq,s\big)\ \bP\big(E^c\big)
    =\hat G\big(\psi_{s'}^s\bq,s\big)-\hat G\big(\psi_{s'}^s\bq,s\big)\ \bP(E)\\
      &=
  \hat G\big(\psi_{s'}^s\bq,s\big)-(s-s')\hat G\big(\psi_{s'}^s\bq,s\big)
   \eta\big(s,\hat\rho_n\big)  +    O((s-s')^2)\la{eq3.11}\\
   &=  \hat G\big(\psi_{s'}^s\bq,s\big)
   -(s-s')\hat G\big(\psi_{s'}^s\bq,s\big)
    \int_{\rho_+}^\i f^2(a_+,s,\hat\rho_n, \rho_+)\ d\rho_+ +    O((s-s')^2),
  \nonumber
  \end{align}
   where \eqref{eq3.9} is used for the third equality.
Moreover, 
 using the strong Markov property for the
  random boundary at the stopping time $\tau$,
  \begin{equation}\la{eq3.12}
    \bE\ \hat G\big(\Psi_{s'}^s \bq,s) 1\!\!1_{E} =\bE\ \hat G\big(\Psi_\tau^{s} \epsilon_{\rho_+} 
    \psi_{s'}^\tau \bq,s\big) 1\!\!1_{E}
    = \bE\
     \hat G\big(\epsilon_{\rho_+} \psi_{s'}^\tau \bq,\tau\big) 
    1\!\!1_{E}.
  \end{equation}
  By \eqref{eq3.8},
  \begin{equation}\la{eq3.13}
    \left|\bE\ \hat G\big(\epsilon_{\rho_+} \psi_{s'}^\tau \bq,\tau\big) 
    1\!\!1_{E}
    -\bE\ \hat G\big(\epsilon_{\rho_+} \psi_{s'}^\tau \bq,s\big) 
    1\!\!1_{E}\right| \le C_1(n+1)(s-s') \bP(E). 
  \end{equation}
Next we modify the distribution from which $\rho_+$ is selected; at
  present, $\rho_+$ is selected according to a random measure
  with density
  \[
      \hat f^2\big(a_+,\tau,\tilde \rho_n, \rho_+\big):= 
\eta\big(\tau,\tilde \rho_n\big)^{-1}
      f^2\big(a_+,\tau, \tilde \rho_n, \rho_+\big) ,    
  \]
  where $\tilde \rho_n:=\g_{s'}^\tau(\hat\rho'_n)$. From $H\in C^2$, and the Lipshitzness of $b^i(x,s,\rho)$ for $i=1,2$, 
 it is not hard to show that there exists a constant $c_2$ such that
\be\la{eq3.14}
\big|\hat\rho'_n-\hat \rho_n\big|\le c_2|s'-s|,\ \ \ \ \big|\hat\rho'_n-\tilde \rho_n\big|\le c_2|s'-s|.
\ee
Let us write $\hat{\rho}_+$ for an independent random variable
  distributed as  $\hat f(a_+,s,\hat\rho_n, \rho_+)\ d\rho_+$.
  Observe 
  \[
  \eta(\th,m)=\int_m^\i f^2(a_+,\th,m,\rho_+)\
  d\rho_+\ge \left[\min_{\th'\in[t_0,T]}H_\rho(a_+,\th',P_-)\right]\int_m^\i  f(a_+,\th,m,\rho_+)\  d\rho_+.
  \]
  From this, \eqref{eq3.14}, and the Lipschitzness of $f^2=vf$ we can readily show
\be\la{eq3.15}
\big|\hat f^2(a_+,s,\hat\rho_n, \rho_+)-
\hat f^2\big(a_+,\tau,\tilde \rho_n, \rho_+\big)\big|\le c_3|s'-s|,
\ee
for a constant $c_3$. 
 We then use \eqref{eq3.14}  and \eqref{eq3.15} to assert that there exists a constant $c_4$ such that the expression
\[
 \left|\bE\ \left[\hat G\big(\epsilon_{\rho_+} \psi_{s'}^\tau \bq,s\big) 
   -  \hat G\big(\epsilon_{\hat\rho_+} \psi_{s'}^\tau \bq,s\big) \right]
   1\!\!1_{E}     \right| ,
\]
 is bounded above by
  \begin{align*}
     &\left|\bE\ 1\!\!1_{E } \int_{\tilde\rho_n\vee \hat\rho_n}^\i \hat G\big(\epsilon_{\rho_+} \psi_{s'}^\tau \bq,s\big) 
\big(\hat f^2\big(a_+,\tau,\tilde \rho_n, \rho_+\big)
-\hat f^2(a_+,s,\hat\rho_n, \rho_+)\big)\ d\rho_+       \right|  \\
&\ \ \ \   +
\left|\bE\ 1\!\!1_{E } \ \1\big(\hat \rho_n\le \tilde\rho_n\big)\int_{\hat \rho_n}^{\tilde\rho_n} \hat G\big(\epsilon_{\rho_+} \psi_{s'}^\tau \bq,s\big) \ \hat f^2\big(a_+,s,\hat \rho_n, \rho_+\big)\ d\rho_+       \right|\\
&\ \ \ \ +
\left|\bE\ 1\!\!1_{E } \
\1\big(\hat \rho_n> \tilde\rho_n\big)\int^{\hat \rho_n}_{\tilde\rho_n} \hat G\big(\epsilon_{\rho_+} \psi_{s'}^\tau \bq,s\big) \ \hat f^2\big(a_+,\tau,\tilde \rho_n, \rho_+\big)\ d\rho_+       \right|  
    \leq c_4(s'-s) \bP(E ).
  \end{align*}
    From this, 
\eqref{eq3.12}, \eqref{eq3.13}, and \eqref{eq3.9}
  we learn 
  \begin{align*}
  \bE\
   \hat G\big(\Psi_{s'}^s \bq,s) 1\!\!1_{E } =&
  \bE \left[ \hat G\big(\epsilon_{\hat\rho_+} \psi_{s'}^\tau \bq,s\big) 
   \ \big|\ E \right]\bP(E )+
  (s-s')^2 R_1 \\
 =& \bE \left[ \int_{\hat\rho_n}^\i
 \hat G\big(\epsilon_{\rho_+} \psi_{s'}^\tau \bq,s\big)
 \  f^2(a_+,s,\hat\rho_n, \rho_+)\ d\rho_+
   \ \big|\ E \right]\eta(s,\hat\rho_n)^{-1}\ \bP(E ) +
   (s-s')^2 R_1\\
 =& (s-s')\int_{\hat\rho_n}^\i
 \bE \left[ \hat G\big(\epsilon_{\rho_+} \psi_{s'}^\tau \bq,s\big)
 \ f^2(a_+,s,\hat\rho_n, \rho_+)\ \big|\ E \right]\ d\rho_+   + (s-s')^2 R_2.
 \end{align*}
  where $R_1$ and $R_2$ are bounded by a constant multiple of
  $n+1$.   
  This and \eqref{eq3.11} complete the proof of \eqref{eq3.10}.

  \ms\noi
{\em(Step 2)} We wish to  establish \eqref{eq3.7} with the aid of  \eqref{eq3.10}. Observe that  $\mu(d\bq,s)$
is the law of a Markov process with a bounded jump rates. 
For such a Markov process,
we can readily show that if $\bn(\bq)$ denotes the number
of jumps/particles of $\bq$  in the interval $[a_-,a_+]$, then 
\be\la{eq3.16}
\sup_{s\in[t_0,T]}\int \bn(\bq)^k\ \mu(d\bq,s)<\i,
\ee
for every $k\in\bN$. Indeed if we choose $\d_0$ so that 
$\l(x,t_0,\rho)\ge \d_0$ for all $(x,\rho)\in [a_-,a_+]\x [P_-,P_+]$, then there exists a Poisson random variable $N_{\d_0}$ of intensity $\d_0^{-1}(a_+-a_-)$ such that 
$\bn(\bq)\le N_{\d_0}$ almost surely.
From \eqref{eq3.16}, and \eqref{eq3.10}, we can write
\[
(s-s')^{-1}\big(\hat G(\bq,s')-\hat G(\bq,s)\big)=
\sum_{r=1}^5\O_r(s',s),
\]
where
\begin{align*}
\O_1(s',s)=&        \int_{\hat\rho_n}^\i      
    \left(\bE \left[\hat G\big( \epsilon_{\rho_+} \psi_{s'}^\tau \bq,s\big) 
    \ \big|\ E \right]- \hat G\big( \epsilon_{\rho_+} \bq,s\big) 
    \right)f^2(a_+,s,\hat\rho_n,\rho_+)\ d\rho_+\\
=&\int_{\hat\rho_{\bn(\bq)}}^\i  \bE  \left[\hat G\big(\e_{\rho_+}\psi^\tau_{s'}\bq,s\big)
-\hat G\big( \e_{\rho_+}\bq,s\big)\big| E \right]
 f^2(a_+,s,\hat \rho_{\bn(\bq)},\rho_+)\ d\rho_+ ,\\
     \O_2(s',s)=&                \int_{\hat\rho_n}^\i      
    \left(\hat G\big( \epsilon_{\rho_+}  \bq,s\big) 
    - \hat G(   \bq ,s)
    \  \right)  f^2(a_+,s,\hat\rho_n,\rho_+)\ d\rho_+=(\cL^s_{bn}\hat G)
(\bq,s),\\         
  \O_3(s',s)=  &\int_{\hat\rho_n}^\i  \left(
 \hat G(\bq,s)-\hat G\big(\psi^s_{s'}\bq,s\big)\right) f^2(a_+,s,\hat\rho_n,\rho_+)\ d\rho_+
= \left( \hat G(\bq,s)-\hat G\big(\psi^s_{s'}\bq,s\big)\right) \eta(s,\hat\rho_n),
\\
\O_4(s',s)=  &(s-s')^{-1}\left(
 \hat G\big(\psi^s_{s'}\bq,s\big)-\hat G(\bq,s)\right),
      \end{align*}
and $|\O_5(s',s)|=(s-s')|R(\bq, s',s)|\le c_1(n+1)(s-s')$. 
(For the second equality, we used the fact that the event $E $ depends only 
on the stochastic boundary that is independent from the law of
$\rho_+$.)
By \eqref{eq3.16},
\[
\int\big| \O_5(s',s)\big|\ \ \mu(d\bq,s)=O(s-s').
\]
As a result, \eqref{eq3.7} would follow if we can show
\begin{align}\la{eq3.17}
\lim_{s'\uparrow s}\int \O_1(s',s)\ \mu(d\bq,s)=0,\\
\la{eq3.18}
\lim_{s'\uparrow s}\int \O_3(s',s)\ \mu(d\bq,s)=0.
\end{align}
and that the limit
\begin{align}
&\lim_{s'\uparrow s}\ 
\int \O_4(s',s)\ \mu(d\bq,s),\la{eq3.19}
\end{align}
equals to the sum of the last two terms on the right-hand side of \eqref{eq3.7}. 
The proof of this will be carried out in the last step. 
A slight modification of this proof can be carried out to establish \eqref{eq3.18}.

\ms\noi
{\em(Step 3)} We turn out attention to \eqref{eq3.17}.
Recall that $\mu(d\bq,s)$ is the law of a Markov process
$(\rho(x,s):\ x\in[a_-,a_+])$ with generator $\cA^1_{x,s}$. 
For our proof we will need a lower bound on the density $f$. Since
$f(x,s,\rho_-,\rho_+)>0$ only when $(\rho_-,\rho_+)$ is in the interior
of $\L(P_-,P_+)$, we wish to estimate the probability of the 
set $B(\d,s)$ consisting of those $\bq$ such that for some $x\in[a_-,a_+]$, we  have
either $\rho(x,s)=R_s(\bq)(x)\in [P_+-\d,P_+]$, or $\rho(x+,s)-\rho(x-,s)\in (0,\d)$.
If we write $\bP^m_s$ for the law of our Markov process $\rho(x,s)$ associated with
the generator $\cA^1_{x,s}$, and the initial condition $\rho(a_-,s)=m$, and
if $m<P_+-\d$, then it is not hard to show
\begin{align*}
\bP_s^m\big(B(\d,s)\big)=&\bE^m_s\int_{a_-}^{a_+}\left[
\int_{\rho(x,s)}^{\rho(x,s)+\d}+\int_{\rho(x,s)\vee(P_+-\d)}^{P_+}
\right]
f(x,s,\rho(x,s),\rho_+)\ d\rho_+\ dx\le c_5 \d .
\end{align*}
From this, we learn that \eqref{eq3.17} would follow if 
we can show
\be\la{eq3.20}
\lim_{s'\uparrow s}\int \O_1(s',s)\ \hat\mu(d\bq,s)=0,
\ee
where
\[
\hat\mu(d\bq,s)=\1\big(\bq\notin B(\d,s)\big)\ \mu(d\bq, s).
\]

\ms\noi
{\em(Step 4)}
To verify \eqref{eq3.20},
write $\s(\bq,s')$ for the first time $\s>s'$ at which $\psi_{s'}^\s\bq$
experiences a collision between particles of $\bq$. 
We claim
\be\la{eq3.21}
\int 1\!\!1\big(\s(\bq,s')\le s\big)\ \mu(d\bq,s)
\le c_5(s-s')\int \bn(\bq)\ \mu(d\bq,s)\le c_6(s-s'),
\ee
for constants $c_6$ and $c_5$.
This is an immediate consequence of
\eqref{eq3.16} and the following fact:
 If $\bq=(x_0,\rho_0,x_1,\rho_1,\dots,x_n,\rho_n)$,
and $\s(\bq,s')\le s$, then for some $i$, we have 
$|x_i-x_{i+1}|\le 2c_6|s-s'|$, where  $c_6$ is an upper bound on the speed of particles.
 Because of  \eqref{eq3.21},
the claim \eqref{eq3.20} is equivalent to 
\be\la{eq3.22}
\lim_{s'\uparrow s} \left|\sum_{n=0}^\i X_n(s')\right|=0,
\ee
where $X_n(s')$ is the expression
\[
\int  \int_{\hat\rho_n}^\i\bE \ \big[\hat 
G\big(\e_{\rho_+}\psi^\tau_{s'}\bq,s\big)-
\hat G\big( \e_{\rho_+}\bq,s\big) |\ E \big]1\!\!1\big(\s(\bq,s')> s\big)\
 f(a_+,s,\hat \rho_n,\rho_+)  \hat\mu^n(\bq,s)\ d\rho_+d\bq.
\]
On account of \eqref{eq3.9}, the claim \eqref{eq3.22} would follow if we can show 
\be\la{eq3.23}
\lim_{s'\uparrow s} (s-s')^{-1}
\left|\sum_{n=0}^\i Y_n(s')\right|=0,
\ee
where $Y_n(s')=Y_n^+(s')-Y_n^-(s')$, with
\begin{align*}
Y^+_n(s')&=\int\int_{\hat\rho_n}^\i\bE \ 
\hat G\big(\e_{\rho_+}\psi^\tau_{s'}\bq,s\big)
1\!\!1\big(\s(\bq,s')> s>\tau(\bq,s')\big)\
 f^2(a_+,s,\hat \rho_n,\rho_+)
  \hat\mu^n(\bq,s)\ d\rho_+d\bq,\\
Y_n^-(s')&=\int\int_{\hat\rho_n}^\i\bE \ 
\hat G\big( \e_{\rho_+}\bq,s\big)1\!\!1\big(\s(\bq,s')> s>\tau(\bq,s')\big)\
 f^2(a_+,s,\hat \rho_n,\rho_+)   \hat\mu^n(\bq,s)\ d\rho_+d\bq.
\end{align*}

\ms\noi
{\em(Step 5)}
The expected value in the definition of $Y_n^\pm$ is for the random variable $\tau=\tau(\bq,s')$. As was explained in the proof of Proposition 2.1{\bf(ii)},
the variable $\tau$ can be expressed in terms of $\hat \rho_n$ and a standard exponential random variable. 
More precisely,
\[
\tau=\tau(\bq,s')=\ell(r,\hat\rho_n,s'),
\]
with $r>0$ a random variable with distribution $e^{-r}\ dr$, and 
$\ell(r,\hat\rho_n,s')$ denoting the inverse of the map
\[
\tau\mapsto r=\int_{s'}^\tau \eta\big(\th,\g_{s'}^{\th}(\th,\hat\rho_n)\big)\ d\th,\ \ \ \ 
\tau\in(s',\i).
\]
As a result, we may replace the expected values in \eqref{eq3.23} with an integration with respect to $e^{-r}\ dr$.
On the other hand,
\[
1\!\!1(r>0)\ e^{-r}\ dr=1\!\!1(\tau>s')\ e^{-r}
\eta\big(\tau,\g_{s'}^{\tau}(\hat\rho_n)\big)\ d\tau
=1\!\!1(\tau>s')\ \big(\eta\big(s',\hat\rho_n\big)+O(\tau-s')\big)\ d\tau,
\]
by the Lipschitz regularity of $\eta$. 
Because of this, \eqref{eq3.23} would follow if we can show
\be\la{eq3.24}
\lim_{s'\uparrow s} (s-s')^{-1}
\left|\sum_{n=0}^\i Z_n(s')\right|=0,
\ee
where $Z_n(s')=Z_n^+(s')-Z_n^-(s')$, with
\begin{align*}
Z_n^+(s')&=\int \int_{\hat\rho_n}^\i\int_{s'}^s \hat G\big(\e_{\rho_+}\psi^\tau_{s'}\bq,s\big)
1\!\!1\big(\s(\bq,s')> s\big)\eta\big(\hat\rho_n,s'\big)
f^2(a_+,s,\hat \rho_n,\rho_+)\hat\mu^n(\bq,s)
  \ d\tau   d\rho_+ d\bq,\\
Z_n^-(s')&=\int\int_{\hat\rho_n}^\i\int_{s'}^s
\hat G\big( \e_{\rho_+}\bq,s\big)1\!\!1\big(\s(\bq,s')> s\big)
  \eta\big(\hat\rho_n,s'\big) f^2(a_+,s,\hat \rho_n,\rho_+)\hat\mu^n(\bq,s)\ d\tau d\rho_+d\bq .
\end{align*}
To prove \eqref{eq3.24}, we  carry out the $d\bq$ integration first. 
  Fix $\tau>0$ and $\rho_+$,
and make a change of variables $\bq'=\psi_{s'}^\tau\bq$ for $d\bq$ integration
in $Z_n^+(s')$.
For this, we wish to replace $\hat\mu^n(\bq,s)$, with 
$\hat\mu^n\big( \psi_{s'}^\tau\bq,s\big)$.
 Observe that by Hypothesis 1.1{\bf(iii)},
then kernel $f(x,s,\rho_-,\rho_+)>0$ in the interior of
$\L(P_-,P_+)$. As a result, we can find $\d_0>0$ such that 
\be\la{eq3.25} 
\bq=(x_0,\rho_0,\dots,x_n,\rho_n)\notin B(\d,s)\ \ \ 
\implies\ \ \ f(x_i,\hat\rho_i,\rho_{i+1})\ge \d_1.
\ee
Since $\hat\mu^n$ is supported 
on the complement of the event $B(\d,s)$, we use 
\[
\left[\log \hat\mu^n\big( \psi_{s'}^\tau\bq,s\big)\right]_\tau=
\bb\big( \psi_{s'}^\tau\bq,\tau\big)\cdot\nabla\left[\log \hat\mu^n\big( \psi_{s'}^\tau\bq,s\big)\right],
\]
our assumption $f\in C^1$, and \eqref{eq3.25} to assert 
\[
\left[\log \hat\mu^n\big( \psi_{s'}^\tau\bq,s\big)\right]_\tau=O(n),
\]
which in turn implies
\be\la{eq3.26}
\hat\mu^n\big( \psi_{s'}^\tau\bq,s\big)=\mu^n(\bq,s)\big(1+(s'-s)O(n)\big).
\ee
 Since the map $\bq\mapsto\psi_{s'}^\tau\bq$ is the flow of the  ODE associated with vector field $\bb$, its Jacobian 
 has the expansion 
\[
1+(\tau-s')div(\bb)+\bn(\bq)\ o(\tau-s').
\] 
Since $div(\bb)=O(\bn(\bq))$, 
a change of variable $\bq'=\psi_{s'}^\tau\bq$ causes
a Jacobian factor of the form 
\[
1+\bn(\bq)O(\tau-s')=1+\bn(\bq)O(s-s').
\]
From this, \eqref{eq3.26}, and \eqref{eq3.14} we learn
\[
 \eta\big(\hat\rho_n,s'\big) f^2(a_+,s,\hat \rho_n,\rho_+)
\mu^n(\bq,s)\ d\bq
 =\eta\big(\hat\rho'_n,s'\big) f^2(a_+,s,\hat \rho'_n,\rho_+)\mu^n(\bq',s)
 \big(1+nO(s-s')\big)\ d\bq'.
 \]
From all this we deduce that $Z_n^+(s')=\hat Z_n^+(s')+R_n$,
where $\hat Z_n^+(s')$ is given by
\[
\int\int\int_{\hat\rho_n}^\i \int_{s'}^s \hat G\big(\e_{\rho_+}\bq',s\big)
1\!\!1\big(\s(\psi^{s'}_{\tau}\bq',s')> s\big)\eta\big(\hat\rho'_n,s'\big)
f^2(a_+,s,\hat \rho_n,\rho'_+)\hat\mu^n(\bq',s)
  \ d\tau   d\rho_+ d\bq'.
\]
and the $R_n$ is an error term that satisfies
\[
\sum_{n=0}^\i|R_n|\le c_2(s-s')^2\int \bn(\bq)^2\ \mu(d\bq,s)
=c_3(s-s')^2.
\]
By $\psi^{s'}_{\tau}$ we mean the inverse of $\psi_{s'}^{\tau}$.
After renaming $\bq'$ as $\bq$ and comparing $\hat Z_n^+(s')$
with $ Z_n^-(s')$, we learn that $ \hat Z_n^+(s')-Z_n^-(s')$ equals to
\[
\int\int\int_{\hat\rho_n}^\i\int_{s'}^s
\hat G\big( \e_{\rho_+}\bq,s\big)\chi(\bq;s',\tau,s)
  \eta\big(\hat\rho_n,s'\big) f^2(a_+,s,\hat \rho_n,\rho_+)
\mu^n(\bq,s)\ d\tau d\rho_+d\bq,
  \]
 where 
 $\chi(\bq;s',\tau,s)=\big|1\!\!1\big(\s(\psi^{s'}_{\tau}\bq,s')> s\big)-
  1\!\!1\big(\s(\bq,s')> s\big)\big|$. 
 After replacing $\hat G$ with an upper bound, and carrying out the 
$d\rho_+$ integration, we obtain
 \[
 \sum_{n=0}^\i\big|  \hat Z_n^+(s')-Z_n^-(s')\big|\le c_4\int_{s'}^s\int \chi(\bq;s',\tau,s)\ \mu(d\bq,s)d\tau.
 \]  
  Finally, since $\chi(\bq;s,s)=0$,
  we can readily show 
  \[
  \lim_{s'\uparrow s} 
(s-s')^{-1}\sum_{n=0}^\i\big(\hat Z^+(s')-Z^-(s')\big)=0,
\]
completing the proof of \eqref{eq3.24}, that in turn completes the
proof of \eqref{eq3.17}.

\ms\noi
{\it(Final Step)} It remains to find the limit in \eqref{eq3.19}. 
The proof we present is very general, and works whenever $\hat G^n$
is continuous, $\mu^n$ is $C^1$, the vector field $\bb=\bb_n$ is $C^1$, and the boundary of $\D_n$ is piecewise $C^1$.
Fix $s$ and for $s'<s$ we write
\[
\D_n(s',s)=\big\{\bq\in \D_n:\ \psi_{s'}^s\bq\in\D_n\big\},\ \ \ \ 
\hat \D_n(s',s)=\psi_{s'}^s\big(\D_n(s',s)\big).
\]
We make a change of variables to write
\begin{align}\nonumber
\int_{\D_n}\hat G^n\big(\psi_{s'}^s\bq,s\big)\mu^n(\bq,s)\ d\bq
=&\int_{\D_n\setminus \D_n(s',s)}\hat G^n\big(\psi_{s'}^s\bq,s\big)\mu^n(\bq,s)\ d\bq\\
&+\int_{\hat \D_n(s',s)}\hat G^n\big(\bq,s\big)\mu^n(\psi^{s'}_s\bq,s)\ \det \big(D\psi^{s'}_s\big)(\bq)\ d\bq,\la{eq3.27}
\end{align}
where $\psi^{s'}_s$ denotes the inverse of the function $\psi_{s'}^s$.
For $s-s'$ small, the volume $|\D_n\setminus \D_n(s',s)|$
is of order $O(n(s-s'))$. From this, and the continuity of $\hat G$ 
we learn 
\begin{align}\nonumber
\int_{\D_n\setminus \D_n(s',s)}\hat G^n\big(\psi_{s'}^s\bq,s\big)\mu^n(\bq,s)\ d\bq=&\int_{\D_n\setminus \D_n(s',s)}\big(\hat G^n\mu^n\big)\big(\bq,s\big)\ d\bq+o(n(s-s'))\\
=&(s-s')\int_{\p\D_n}\big(\hat G^n\mu^n\big)\big(\bq',s\big)
\big({\bold{N}}_n(\bq')\cdot \bb_n(\bq',s)\big)\ \s(d\bq')\la{eq3.28}\\
&+o(n(s-s')).\nonumber
\end{align}
Here we have used the fact that we may parametrize the set 
$\D_n\setminus \D_n(s',s)$ by the map
\[
\z:\p \D_n\x [s',s]\to \D_n\setminus \D_n(s',s),\ \ \ \ 
\z(\bq',\th)=\psi_s^\th\bq',
\]
with $\1\big(\bq\in\D_n\setminus \D_n(s',s)\big)\ d\bq$ equals to 
\begin{align*}
&
\1\big((\bq',\th)\in\p \D_n\x [s',s]\big)
\left(1+(s-s')\big({\bold{N}}_n(\bq')\cdot \bb_n(\bq',s)\big)
+o(n(s-s'))\right)
\ \s(d\bq')d\th.
\end{align*}
The map $\z$ is one-to-one if $\p\D_n$ is $C^1$, and $s-s'$ is sufficiently small. This is no longer the case when $C^1$ is only piecewise $C^1$. Though the set of
$\bq$ for which $\z^{-1}(\bq)$ is multivalued,
is  of volume $O(n(s-s')^2)$ (this is the set of 
$\bq$ such that for some $i$, we have $x_{i+1}-x_i,x_i-x_{i-1}=
O(s-s')$).

As for the second term on the right-hand side of \eqref{eq3.27},
we use
\begin{align*}
\mu^n(\psi^{s'}_s\bq,s)=&\mu^n(\bq,s)
+(s'-s)\big(\bb_n\cdot \nabla \mu^n\big)(\bq,s)+o(n(s-s'))\\
\det \big(D\psi^{s'}_s\big)(\bq)=&\mu^n(\bq,s)+(s'-s)\big(\mu^n div\ \bb_n\big)(\bq,s)+o(n(s-s')),
\end{align*}
to assert
\begin{align*}
\mu^n(\psi^{s'}_s\bq,s)\ \det \big(D\psi^{s'}_s\big)(\bq)
=&\mu^n(\bq,s)+(s'-s)\big(\mu^n div\ \bb_n\big)(\bq,s)\\
&+(s'-s)\big(\bb_n\cdot \nabla \mu^n\big)(\bq,s)+o(n(s-s'))\\
=&\mu^n(\bq,s)+(s-s')\big(\cL_{0n}^{s*}\mu^n \big)(\bq,s)+o(n(s-s')).
\end{align*}
From this and \eqref{eq3.28} we deduce that the second term on 
 the right-hand side of \eqref{eq3.27} 
equals to
\begin{align*}
\int_{\hat \D_n(s')}&\hat G^n\big(\bq,s\big)\hat G^n\big(\bq,s\big)\big[\mu^n(\bq,s)+(s-s')\big(\cL_{0n}^{s*}\mu^n \big)(\bq,s)\big]\ d\bq+o(n(s-s'))\\
=&\int_{ \D_n}\big[\mu^n(\bq,s)+(s-s')\big(\cL_{0n}^{s*}\mu^n \big)(\bq,s)\big]\ d\bq+o(n(s-s')).
\end{align*}
This, \eqref{eq3.27}, and \eqref{eq3.28} complete the proof.
\qed

\section{Proof of Theorem 2.1}

The proof of Theorem~2.1 is carried out in five steps.

\ms\noi
{\em(Step 1)} As we explained in Section 3, we only need to prove \eqref{eq3.5}.
For this, it suffices to show
\be\la{eq4.1}
\lim_{s'\uparrow  s} (s-s')^{-1}
\big({\mathbb G}_n(s)-{\mathbb G}_n(s')\big)=0,
\ee
where  ${\mathbb G}_n(s)=\int\hat G^n(\bq,s)\ \mu^n(d\bq,s)$.
Evidently
\[
(s-s')^{-1}\big({\mathbb G}_n(s)-{\mathbb G}_n(s')\big)
=\O_1(s')+\O_2(s')-\O_3(s'),
\]
where
 \begin{align*}
  \O_1(s')&=(s-s')^{-1}
  \int \big(\hat G^n(\bq,s)-\hat G^n(\bq,s')\big) \mu^n(d\bq,s)\\
  \O_2(s')&=(s-s')^{-1}\int \hat G^n(\bq,s)\ 
\big(\mu^n(d\bq,s)-\mu^n(d\bq,s')\big)\\
  \O_3(s')&=(s-s')^{-1}\int \big(\hat G^n(\bq,s)-\hat G^n(\bq,s')\big)
  \ \big(\mu^n(d\bq,s)-\mu^n(d\bq,s')\big).
    \end{align*}
We claim 
 \be\la{eq4.2}
  \lim_{s'\uparrow  s}\big|\O_3(s')\big|=0.
  \ee  
By Lemma~3.1,
  \be\la{eq4.3}
 \limsup_{s\to 0} (s-s')^{-1}\big|\O_3(s')\big|\le C_0(n+1)
 \limsup_{s\to 0}\int 
   \left|\big(\mu^n(d\bq,s)-\mu^n(d\bq,s')\big)\right|.
\ee
As we will see in Step 2 below, $\mu_s^n=X^n \mu^n$ for a term $X^n$ that is
explicit. From this and \eqref{eq4.3}, it is not hard to deduce \eqref{eq4.2}.
From \eqref{eq4.2}, and Theorem 3.1 we deduce that \eqref{eq4.1} would follow if we can show
\begin{align}\nonumber
\int_{\D_n}\hat G^n(\bq,s)\ \mu^n_s(d\bq,s)=&\int_{\D_n} (\cL_{b}^s \hat G)^n(\bq,s)\mu^n(\bq,s)\ d\bq+ 
\int_{\D_n}  \hat G^n(\bq,s)\big(\cL_{0n}^{s*}
\mu^n\big)(\bq,s)\ d\bq\\
&+\int_{\hat\p \D_{n+1}}  \hat G^{n+1}(\bq,s)\mu^{n+1}(\bq,s)(\bb_{n+1}\cdot {\bold{N}}_{n+1})\ \s(d\bq).\la{eq4.4}
\end{align}

\ms\noi
{\em(Step 2)}  
To simplify our presentation, we assume that
$\ell$ has a density with respect to the Lebesgue measure. With a slight abuse of notation, we write $\ell(\rho,s)$ for this density: 
$\ell(d\rho,s)=\ell(\rho,s)\ d\rho$. 
To verify \eqref{eq4.4}, we start with finding a tractable expression
for the left-hand side. We claim
\be\la{eq4.5}
\mu_s^n=X^n\mu^n=\left(X_1+X_2+\sum_{i=3}^{11}X^n_i\right)\mu^n,
\ee
 where
\begin{align*}
X_1&= \frac{(\ell* f^2)(x_0,s,\rho_0)}{\ell(s,\rho_0)},
\ \ \ \ \ \ \ \ \ \ \ \ \ \ \ \
X_2=-\frac{\big(\b(x_0,s,\rho_0)\ell(s,\rho_0)\big)_{\rho_0}}{\ell(s,\rho_0)},\\
X^n_3&=\sum_{i=1}^{n}\frac
{Q^+(f)(x_i,s,\hat\rho_{i-1},\rho_{i})}
{f(x_i,s,\hat\rho_{i-1},\rho_{i})},\ \ \ \ \ \
X^n_4=\sum_{i=1}^n 
\frac{f^2_x\big(x_i,t,\hat\rho_{i-1},\rho_i\big)}
{f\big(x_i,t,\hat\rho_{i-1},\rho_i\big)}\\
X^n_5&=\b(x_0,t,\rho_0)
\G_{\rho}(x_0,x_{1},s,\rho_0)
-A(a_+,s,\hat\rho_n)\\
X^n_6&=\sum_{i=1}^{n}\b(x_i,t,\rho_i)
\G_{\rho}(x_i,x_{i+1},s,\rho_i)\\
X^n_7&=\sum_{i=1}^{n}v(x_i,s,\hat \rho_{i-1},\rho_i)
\big(\l(x_i,s,\rho_{i})-\l(x_i,s,\hat\rho_{i-1})\big)\\
X^n_8&=\sum_{i=1}^{n}
\frac{\big[K(x_i,s,\rho_i,\hat\rho_{i-1})f(x_i,s,\hat\rho_{i-1},\rho_i)\big]_{\rho_i}}
{f(x_i,s,\hat\rho_{i-1},\rho_i\big)}\\
X^n_9&=\sum_{i=1}^{n}b(x_i,s,\hat\rho_{i-1})
v_{\rho_-}(x_i,s,\hat \rho_{i-1},\rho_i)\\
X^n_{10}&=\sum_{i=1}^{n}v(x_i,s,\hat \rho_{i-1},\rho_i)
b(x_i,s,\hat\rho_{i-1})\frac{f_{\rho_-}(x_i,s,\hat\rho_{i-1},\rho_{i})}{f(x_i,s,\hat\rho_{i-1},\rho_{i})}\\
X^n_{11}&=-\sum_{i=1}^{n}
\b(x_{i-1},s,\rho_{i-1})\frac{\big[f(x_i,s,\hat\rho_{i-1},\rho_{i})\big]_{\rho_{i-1}}}{f(x_i,s,\hat\rho_{i-1},\rho_{i})},
\end{align*}
where $\hat \rho_{i-1}=\phi_{x_{i-1}}^{x_i}\big(\rho_{i-1};s\big)$.
To verify \eqref{eq4.5},  observe that by
direct differentiation  (see Definition 2.2{\bf(ii)} for the definition of $\mu^n$)
\be\la{eq4.6}
X^n=-\G_s(\bq,s)+\frac 
{ \ell_s(s,\rho_0)}{\ell(s,\rho_0)}
+\sum_{i=1}^{n}\frac{\big[f(x_i,s,\hat\rho_{i-1},\rho_{i})]_s}
{f(x_i,s,\hat\rho_{i-1},\rho_{i})}.
\ee
By \eqref{eq2.10}, \eqref{eq1.18}, and \eqref{eq2.8} (in this order)
\begin{align*}
&-\G_s(\bq,s)=-\sum_{i=0}^{n}
\left\{\big(A(x_{i+1},s,\hat\rho_i)-A(x_i,s,\rho_i)\big)-\b(x_i,t,\rho_i)
\G_{\rho}(x_i,x_{i+1},s,\rho_i)\right\}\\
&\ \ \ \ \ \ \ \ \ \ \ \ \ =-\sum_{i=0}^{n}
\big(A(x_{i+1},s,\hat\rho_i)-A(x_i,s,\rho_i)\big)+\b(x_0,t,\rho_0)
\G_{\rho}(x_0,x_{1},s,\rho_0)+X^n_6\\
&\ \ \ \ \ \ \ \ \ \ \ \ \ =\sum_{i=1}^{n}\big(A(x_i,s,\rho_i)-
A(x_{i},s,\hat\rho_{i-1})\big)+
A(x_0,s,\rho_0)-A(x_{n+1},s,\hat\rho_n)\\
&\ \ \ \ \ \ \ \ \ \ \ \ \ \ \ \ \ +\b(x_0,t,\rho_0)
\G_{\rho}(x_0,x_{1},s,\rho_0)+X^n_6,\\
&\frac{\ell_s(s,\rho_0)}{\ell(s,\rho_0)}=X_1+X_2-A(x_0,s,\rho_0),\\
&\frac{\big[f(x_i,s,\hat\rho_{i-1},\rho_{i})]_s}
{f(x_i,s,\hat\rho_{i-1},\rho_{i})}=\frac{f_s(x_i,s,\hat\rho_{i-1},\rho_{i})}{f(x_i,s,\hat\rho_{i-1},\rho_{i})}+U^n(i),
\end{align*}
where
\begin{align*}
U^n(i)=&\left[\b(x_i,s,\hat\rho_{i-1})-
\b(x_{i-1},s,\rho_{i-1})
\big[\phi_{x_{i-1}}^{x_i}(\rho_{i-1};s)\big]_{\rho_{i-1}}\right]
\frac{f_{\rho_-}(x_i,s,\hat\rho_{i-1},\rho_{i})}
{f(x_i,s,\hat\rho_{i-1},\rho_{i})}\\
=&\b(x_i,s,\hat\rho_{i-1})
\frac{f_{\rho_-}(x_i,s,\hat\rho_{i-1},\rho_{i})}
{f(x_i,s,\hat\rho_{i-1},\rho_{i})}-
\b(x_{i-1},s,\rho_{i-1})\frac{\big[f(x_i,s,\hat\rho_{i-1},\rho_{i})\big]_{\rho_{i-1}}}{f(x_i,s,\hat\rho_{i-1},\rho_{i})}.
\end{align*}
 From this, \eqref{eq4.6}, and the kinetic equation we deduce
\begin{align*}
X^n=&X_1+X_2+X_3^n+X_4^n+X_5^n+X_6^n+X^n_{7}+U^n+W^n,
\end{align*}
where $U^n=\sum_{i =1}^nU^n(i)$, and 
\begin{align*}
W^n&=\sum_{i=1}^{n}\frac{(Cf)(x_i,s,\hat\rho_{i-1},\rho_{i})}
{f(x_i,s,\hat\rho_{i-1},\rho_{i})}=
X^n_8+X^n_9+\sum_{i=1}^{n}K(x_i,s,\hat \rho_{i-1},\rho_i)
\frac{f_{\rho_-}(x_i,s,\hat\rho_{i-1},\rho_{i})}
{f(x_i,s,\hat\rho_{i-1},\rho_{i})}.
\end{align*}
We are done because 
$U^n+W^n=X_8^n+X_9^n+X_{10}^n+X_{11}^n.$

 \ms \noi
 {\em(Step 3)} We now turn our attention to the right-hand side of 
\eqref{eq4.4}.
We certainly have
\be\la{eq4.7}
\int (\cL^s_b \hat G)^n(\bq,s) \mu^n(\bq,s)\ d\bq=Y^n_{b,+}- Y^n_{b,-},
\ee
where
  \begin{align*}  
 Y^n_{b,+} =&
 \int\int_{\hat\rho_n}^\i    f^2(a_+,s,\hat\rho_{n},\rho_+)
       \hat G^{n+1}\big(\e_{\rho_+}\bq,s\big) \mu^{n}(\bq,s)\   d\rho_+d\bq ,    \\
Y^n_{b,-}  =&  \int A(a_+,s,\hat\rho_n)\hat G^n(\bq,s)
   \mu^n(\bq,s)\   d\bq.
   \end{align*}
As for the second term on the right-hand side of \eqref{eq4.4},
we write
$\cL^*_{0n}\mu^n=Z^n\mu^n$, with
  \be\la{eq4.8}
  Z^n=\sum_{j=1}^3Z_{1j}+\sum_{i=2}^3\sum_{j=1}^3Z^n_{ij},
  \ee
 where
\begin{align*}
Z_{11}&=\b(x_0,s,\rho_0)\G_{\rho}(x_0,x_1,s,\rho_0)\\
Z_{12}&=-\b(x_0,s,\rho_0)\frac{\big[f\big(x_1,s,\hat\rho_0,\rho_{1}\big)\big]_{\rho_0}}
{f\big(x_1,s,\hat\rho_0,\rho_{1}\big)} ,\ \ \ \ \ \ \ \ 
Z_{13}=-\frac{\big(\b(x_0,s,\rho_0){\ell(s,\rho_0)}\big)_{\rho_0}}
{\ell(s,\rho_0)},\\
Z^{n}_{21}&=-\sum_{i=1}^{n} 
K(x_i,s,\rho_i,\hat\rho_{i-1})
\G_{\rho}(x_i,x_{i+1},s,\rho_i),\\
Z^n_{22}&=\sum_{i=1}^{n-1}
\frac{\big[K(x_i,s,\rho_i,\hat\rho_{i-1})
f(x_i,s,\hat\rho_{i-1},\rho_i)
f(x_{i+1},s,\hat\rho_{i},\rho_{i+1})\big]_{\rho_i}}
{f(x_i,s,\hat\rho_{i-1},\rho_i\big)
f(x_{i+1},s,\hat\rho_{i},\rho_{i+1})}\\
Z^n_{23}&=\frac
{\big[K(x_n,s,\rho_n,\hat\rho_{n-1})
f(x_n,s,\hat\rho_{n-1},\rho_n)\big]_{\rho_n}}
{f(x_n,s,\hat\rho_{n-1},\rho_n)},\\
Z^n_{31}&=\sum_{i=1}^{n-1} 
\frac{\big[f^2\big(x_i,s,\hat\rho_{i-1},\rho_i\big)
f\big(x_{i+1},s,\hat\rho_{i},\rho_{i+1}\big)\big]_{x_i}}
{f\big(x_i,s,\hat\rho_{i-1},\rho_i\big)
f\big(x_{i+1},s,\hat\rho_{i},\rho_{i+1}\big)},\ \ \ \ 
Z^n_{32}=
\frac{\big[f^2\big(x_n,s,\hat\rho_{n-1},\rho_n\big)\big]_{x_n}}
{f\big(x_n,s,\hat\rho_{n-1},\rho_n\big)},\\
Z^n_{33}&=-\sum_{i=1}^n v\big(x_i,s,\hat\rho_{i-1},\rho_i\big)
\big[\G(x_{i-1},x_i,s,\rho_{i-1})+\G(x_i,x_{i+1},s,\rho_i)\big]_{x_i}.
\end{align*}
Recall that $\cL^*_{0n}$ is the adjoint of $\cL_{0n}$, and is obtained by an integration by parts. More specifically, 
\bi
\item 
The sum $Z_{11}+Z_{12}+Z_{13}$ comes from an integration by parts with respect to the variable $ \rho_0$, and the $i$-terms in $Z^n_{21}$, $Z^n_{22}$ 
come from an integration by parts with respect to the variable
 $ \rho_i$ for $i\in\{1,\dots,n-1\}$, and $Z^n_{23}$ comes from an integration by parts with respect to the variable
 $ \rho_n$. The dynamics of $\rho_i$ as in rule {\bf(2)} of 
 Definition 2.1{\bf(iii)} is responsible for these contributions.
 \item The $i$-th terms in $Z^n_{31}$, $Z^n_{32}$ 
and $Z^n_{33}$
come from an integration by parts with respect to the variable $x_i$.
 The dynamics of $x_i$ as in rule {\bf(1)} of Definition 2.1{\bf(iii)}
 is responsible for this contribution. 
\ei

\ms\noi
{\em(Step 4)}
 We next focus on the third term on the right-hand side
of \eqref{eq4.4}. This term can be expressed as 
  \begin{align}
  Y_0^n&=  \sum_{i=0}^{n} Y^n_{0i}+\hat Y^n_{0},
  \la{eq4.9}
  \end{align}
  where $Y^n_{0i}$ is the boundary contribution coming from
  the condition $x_i=x_{i+1}$,  and $\hat Y^n_{0}$
  is the boundary contribution coming from
  the condition $x_{n+1}=x_{n+2}=a_+$.  
For $i=0,\dots,n$
\be\la{eq4.10}
Y^n_{0i}=\int_{\D_n}
 \hat G^n(\bq,s) W_{i}(\bq,s)\mu^n(\bq,s)\ d\bq,
 \ee
 where 
\begin{align*}
 W_{0}&=\frac{\int f^2(x_0,s,\rho_*,\rho_0)
\  \ell(s,d\rho_*)}{\ell(s,\rho_0)},\ \ \ \ \ \ 
W_{i}=\frac{Q^+(f)\big(x_i,s,\hat\rho_{i-1},\rho_{i})}
{f\big(x_i,s,\hat\rho_{i-1},\rho_i\big)},
\end{align*}
for $i=1,\dots,n$. Here,
\bi
\item The term $W_{0}$ comes from the boundary term
 $x_1=x_0=a_-$ in the  integration by parts with respect to the variable $x_1$. This boundary condition represents the event that $x_1$ has reached $x_0$ after which $\rho_0$ becomes $\rho_1$,
 and $(x_i,\rho_i)$ is relabeled as $(x_{i-1},\rho_{i-1})$ for $i\ge 2$.
 \item The term $W_{i}$ comes from the boundary term
  $x_i=x_{i+1}$ .
 The relative distance $x_{i+1}-x_i$ travels with speed 
 \[
 -\big[v(x_{i+1},s,\hat\rho_i,\rho_{i+1})-
 v(x_i,s,\hat\rho_{i-1},\rho_i)\big],
 \]
  As $x_{i+1}$ catches up with $x_i$, the particle $x_i$ disappears and its density
 $\rho_i=\hat \rho_{i}$ is renamed $\rho_*$, and is integrated out.
 (The resulting integral is 
 $Q^+(f)\big(x_i,s,\hat\rho_{i-1},\rho_{i})$.) 
 We then relabel $(x_j,\rho_j)$, $j>i$,
as $(x_{j-1},\rho_{j-1})$.
\ei
As for $\hat Y^n_{0}$, we simply have 
\be\la{eq4.11}
\hat Y^n_{0}=-Y^n_{b,+},
 \ee
where $Y^n_{b,+}$ was defined in \eqref{eq4.7}.

\ms\noi
{\em(Step 5)}
Recall that we wish to establish \eqref{eq4.4}.  
The identities \eqref{eq4.5}, and \eqref{eq4.7}-\eqref{eq4.11} allow us to rewrite
 \eqref{eq4.4} as
\[
X_1+X_2+\sum_{i=3}^{11}X_i^n=\sum_{j=1}^3Z_{1j}+\sum_{i=2}^3\sum_{j=1}^3Z^n_{ij} -A(a_+,s,\hat\rho_n)+W_0+W^n,
\]
where $W^n=\sum_{i=1}^nW^n_{i}$.
For this we only need to verify
\be\la{eq4.12}
X_4^n+\sum_{i=6}^{11}X_i^n=Z_{12}+\sum_{i=2}^3\sum_{j=1}^3Z^n_{ij},
\ee
because
\[
X_1=W_0,\ \ \ \ X_2=Z_{13},\ \ \ \ X_3^n=W^n,\ \ \ \ X_5^n=Z_{11}-A(a_+,s,\hat\rho_n).
  \]
 We use the definition of $K$ to write $Z^n_{21}=Z^n_{211}+Z^n_{212}$, where $Z^n_{212}=X^n_6$, and 
\[
Z^n_{211}=-\sum_{i=1}^{n} b(x_i,s,\rho_{i})
v(x_i,s,\hat\rho_{i-1},\rho_i)\big)
\G_{\rho}(x_i,x_{i+1},s,\rho_i).
\]
Hence \eqref{eq4.12} is equivalent to
\be\la{eq4.13}
X_4^n+\sum_{i=7}^{11}X_i^n=Z_{12}+
Z^n_{211}+Z^n_{22}+Z^n_{23}+
\sum_{j=1}^3Z^n_{3j}.
\ee
Also observe that the expression 
\[
\big[\G(x_{i-1},x_i,s,\rho_{i-1})+\G(x_i,x_{i+1},s,\rho_i)\big]_{x_i},
\]
equals
\begin{align*}
&\left[\int_{x_{i-1}}^{x_i}\l\big(z,s,\phi_{x_{i-1}}^z(\rho_{i-1};s)\big)\ dz+
\int^{x_{i+1}}_{x_i}\l\big(z,s,\phi_{x_{i}}^z(\rho_{i};s)\big)\ dz\right]_{x_i}\\
&\ \ \ \ \ \ \ \ \ =\l(x_i,s,\hat\rho_{i-1})-\l(x_i,s,\rho_{i})
+\int^{x_{i+1}}_{x_i}\left[\l\big(z,s,\phi_{x_{i}}^z(\rho_{i};s)\big)\right]_{x_i}\ dz
\end{align*}
From this and \eqref{eq2.6} we learn
\[
Z^n_{211}+Z^n_{33}=-\sum_{i=1}^nv\big(x_i,s,\hat\rho_{i-1},\rho_i\big)
\big(\l(x_i,s,\hat\rho_{i-1})-\l(x_i,s,\rho_{i})\big)=X^n_7.
\]
This reduces \eqref{eq4.13} to
\be\la{eq4.14}
X^n_4+\sum_{i=8}^{11}X_i^n=Z_{12}+Z^n_{22}+Z^n_{23}+Z^n_{31}+Z^n_{32}.
\ee

 Observe that 
$Z^n_{22}+Z^n_{23}=\widehat Z^n_{22}+\widehat Z^n_{23}$, and  $Z^n_{31}+Z^n_{32}=Z^n_{311}+Z^n_{312}+Z^n_{313}$, where
\begin{align*}
\widehat Z^n_{22}&=\sum_{i=1}^{n}\frac{
\big[K(x_i,s,\rho_i,\hat\rho_{i-1})
f(x_i,s,\hat\rho_{i-1},\rho_i)\big]_{\rho_i}}
{f(x_i,s,\hat\rho_{i-1},\rho_i\big)},\\
\widehat Z^n_{23}&=\sum_{i=1}^{n-1}
K(x_i,s,\rho_i,\hat\rho_{i-1})
\frac{\big[f(x_{i+1},s,\hat\rho_{i},\rho_{i+1})\big]_{\rho_i}}
{f(x_{i+1},s,\hat\rho_{i},\rho_{i+1})},\\
Z^n_{311}&=\sum_{i=1}^n 
\frac{f^2_x\big(x_i,s,\hat\rho_{i-1},\rho_i\big)}
{f\big(x_i,s,\hat\rho_{i-1},\rho_i\big)},\\
Z^n_{312}&=\sum_{i=1}^n b\big(x_i,s,\hat\rho_{i-1}\big)
\frac{f^2_{\rho_-}\big(x_i,s,\hat\rho_{i-1},\rho_i\big)}
{f\big(x_i,s,\hat\rho_{i-1},\rho_i\big)},\\
Z^n_{313}&=-\sum_{i=1}^{n-1} v\big(x_i,s,\hat\rho_{i-1},\rho_{i}\big)
b\big(x_i,s,\rho_{i}\big)
\frac{\big[f\big(x_{i+1},s,\hat\rho_{i},\rho_{i+1}\big)\big]_{\rho_i}}
{f\big(x_{i+1},s,\hat\rho_{i},\rho_{i+1}\big)},
\end{align*}
where we used \eqref{eq2.5} for the last equation. 
Observe that by the definition of $K$,
\[
\hat Z^n_{23}+Z^n_{313}=-\sum_{i=1}^{n-1}
\b(x_i,s,\rho_i)
\frac{\big[f(x_{i+1},s,\hat\rho_{i},\rho_{i+1})\big]_{\rho_i}}
{f(x_{i+1},s,\hat\rho_{i},\rho_{i+1})}.
\]
The equation
\eqref{eq4.14} follows because
\[
X_8^n=\hat Z^n_{22},\ \ \ \ X_9^n+X_{10}^n= Z^n_{312}, \ \ \ \ X_4^n=Z^n_{311},
\ \ \ \ X_{11}^n=Z_{12}+\hat Z^n_{23}+Z^n_{313}.
\]
\qed

\section{Proof of Theorem 1.2}

According to Theorem~1.2 if $\rho(\cdot,t_0)=\rho(\cdot,t_0;s,\by^0)$,
with  $\by^0$ a Markov jump process with jump rate
$g^0(x,y_-,y_+)$, 
 then for $t>t_0$, we can express
$\rho(\cdot,t)=\rho(\cdot,t;s,\by_{t})$, where $\by_t$ is also a jump process with a jump rate $g(x,t;y_-,y_+)$, where $g$ is a solution to the kinetic equation \eqref{eq1.25}. There is a one-to-one correspondence
between the realization
\[
\by(x)=\sum_{i=0}^\i y_i\1\big(x\in [x_i,x_{i+1})\big),
\]
and the particle configuration
\[
\bq=\big((x_0,y_0),(x_1,y_1),\dots\big),
\]
with 
\[
x_0=a_-<x_1<\dots<x_n<\dots,\ \ \ \ \ y_0<y_1<\dots<y_n<\dots.
\]
We may translate this into a statement about the law of
our particle system $\bq(t)$. As before, it suffices to establish a variant of Theorem 1.2 for a finite interval $[a_-,a_+]$. The condition
$H_\rho(a_-,t,\rho)>0$ means that  particles can cross $a_-$ only from left to right. Because of this, we can treat $a_-$ as a free boundary point. As it turns out, the point $a_+$ will be free boundary point 
(no particle can cross $a_+$ from right to left)
if $a_+$ is sufficiently large. Indeed as we will see in Proposition 5.2{\bf(iii)},
there are positive constants $C_0$ and $C_1$ such that $M(x,t;y,s)\le -C_1x$
for $x\ge C_0$. On the other hand, by Hypothesis 1.2{\bf(i)}, we know
\[
|\rho|\le c_1\big(1+|H_\rho(x,t,\rho)|\big),
\]
which in turn implies that $H_\rho(x,t,\rho)\to -\i$ as $\rho\to -\i$.
As a result, there exists a positive constant $C_2$ such that
$H_\rho(x,t,\rho)\le 0$, whenever $\rho\le -C_2$. From this we deduce
that $\hat v(a_+,t,y_-,y_+)>0$ if $a_+\ge \max\{C_2C_1^{-1},C_0\}=:C_3$. 
From all this we learn that Theorem 1.2 would follow if we can establish the
following result.
\begin{theorem}
  \label{th5.1}
  Assume Hypothesis~1.2.  For
  any fixed $a_+ $ such that $a_+> \max\{a_-,C_3\}$, consider the scalar conservation law
\eqref{eq1.2} in $[a_-,a_+] \times [t_0,T)$ 
   with initial condition  $\rho(x,t_0) =M(x,t_0; \by^0(x),s)$
(restricted to $[a_-,a_+]$), open
  boundary at $x = a_\pm$.
  Then for all $t > t_0$, we have  
$\rho(x,t)=M(x,t;\by_t(x),s)$, where the law of $\big(\by_t(x):\ x\in[a_-,a_+]\big)$ is  as follows:
\bi
\item[{\bf(i)}] The $x = a_-$ marginal is $\ell(t,dy_0)$, given by 
$\dot \ell=\cB^{2*}_{a_-,t}\ell$.
 
\item[{\bf(ii)}] The rest of the path is a PDMP with generator 
$\cB^1_{x,t}$.
 \ei
\end{theorem}

To prove our main result Theorem~1.2,
we send $a_+ \to \infty$. 

We continue with a preparatory definition.

 \ms\noi
 {\bf Definition 5.1(i)} The configuration space for our particle system $\bq$, is the set $\D=\cup_{n=0}^\i\bar\D_n,$
where $\bar\D_n$ is the topological closure of $\D_n$, with $\D_n$ denoting the set
\[
\big\{\bold{q}=\big((x_i,y_i):i=0,1,\dots,n\big):\ x_0=a_-<x_1<\dots<x_n<x_{n+1}=a_+,
\ \ \ y_0<\dots<y_n\big\}.
\]
We write $\bn(\bq)$ for the number of particles i.e.,
$\bn(\bq)=n$ means that $\bq\in\D_n$.

\ms\noi
{\bf(ii)}
Given a realization 
$\bold{q}=\big(x_0,y_0,x_1,y_1,\dots,x_n,y_n\big)\in \bar\D_n,$ we define
\begin{align*}
\rho\big(x,t;\bold{q}\big)&=
R_t(\bq)(x)=\sum_{i=0}^{n}M(x,t;y_i,s)
1\!\!1\big(x_i\le x<x_{i+1}\big).
\end{align*}

\ms\noi
{\bf(iii)} The process $\bq(t)$ evolves according to the following 
rules:
\bi
\item[$\bold{(1)}$] So long as $x_i$ remains in $(x_{i-1},x_{i+1})$, it satisfies $\dot x_i=-\hat v(x_i,t,y_{i-1},y_i)$.
\item[$\bold{(2)}$] When $x_1$ reaches $a_-$, we relabel 
particles $(x_i,y_i), \ i\ge 1$, as $(x_{i-1},y_{i-1})$.
\item[$\bold{(3)}$] When $x_n$ reaches $a_+$, a particle is lost and $\bq$
enters $\D_{n-1}$.
\item[$\bold{(4)}$] When $x_{i+1}=x_i$, then 
$\bold{ q}(t)$ becomes
$\bold{ q}^i(t)$, that is obtained from $\bold{ q}(t)$ by omitting
$(x_i,y_i)$ and relabeling particles to the right of the $i$-th particle.
\ei
\qed

\bs
Some care is needed for the rule {\bf(1)} because $\hat v$ given by 
\eqref{eq1.24} is not a continuous function of $x$. Recall that $x_i(t)$
represents the location of a shock discontinuity that separates two 
fundamental solutions. However, the fundamental solution 
$\big(M(x,t;y_i,s):\ x\in (x_{i-1}(t),x_i(t))$ may also include
some shock discontinuities. When $x_i(t)$ catches up with 
a shock discontinuity of $M(\cdot,t;y_i,s)$, or $M(\cdot,t;y_{i+1},s)$,
$\dot x_i(\cdot)$ fails
to exists. Nonetheless off of such a discrete set of moments, the ODE
of {\bf(1)} is  well-defined, and this is good enough to determine
the evolution of $x_i$ fully. 

\bs
We write $\widehat\cL=\widehat\cL^t$ for the generator of the 
(inhomogeneous Markov)
process $\bq(t)$. This generator can be expressed 
as $\widehat\cL=\widehat\cL_0+\widehat\cL_b$, where $\widehat\cL_0$ is the generator of the deterministic part of dynamics, and $\widehat\cL_b$ 
represents the Markovian boundary dynamics. The deterministic and stochastic dynamics restricted to 
 $\D_n$ have generators that are denoted by $\widehat\cL_{0n}$
and $\widehat\cL_{bn}$ respectively. 
While $\bq(t)$ remains in $\D_n$, its evolution is
governed by an ODE of the form
\[
\frac{d\bq}{dt}(t)=\hat\bb\big(\bq(t),t\big),
\]
where $\hat\bb$ can be easily described with the aid of rule {\bf(1)}  of Definition~5.1{{\bf(iii)}. 
  We establish Theorem 5.1 by verifying the forward equation
\be\la{eq5.1}
\dot\mu^n= \big(\widehat\cL^*\mu\big)^n,
\ee
for all $n\ge 0$,
where  $\widehat\cL^*$ is the adjoint of the operator $\widehat\cL$. We follow our 
strategy as in Section 3 and use a test function $G(\bq,t)$ with is the analog of what we had in \eqref{eq3.4}. Again, our Theorem 5.1 would follow if we can show the analog of \eqref{eq3.5}. We follow our notation
as in \eqref{eq3.2}, and the analog of Theorem 3.1 is also valid
when $\cL$ is replaced with $\widehat\cL$.

\bs
The following variant of Proposition~2.1 ensures that our particle system produces the unique entropy solution of \eqref{eq1.1} in the interval $[a_-,a_+]$.

 \bp\la{pro5.1} The function $\rho(x,t)=\rho(x,t;\bq(t))$,
 with $\bq(t)$ evolving as in Definition 5.1{\bf(iii)}, is the unique entropy solution of 
$\rho_t=H(x,t,\rho)_x$ in $(a_-,a_+)\x (0,\i)$.
\ep

\ms\noi
{\bf Proof}
As in Section 2, we can readily check that $\rho(x,t;\bq(t))$ is a weak solution of
\eqref{eq1.2} because the Rankin-Hugoniot condition is satisfied. To satisfy the entropy condition, we need to make sure that $\rho(x-,t;\bq(t))<
\rho(x+,t;\bq(t))$ at each discontinuity point. This is an immediate consequence 
of the monotonicity of the fundamental solution that is stated in Proposition 5.2{\bf(ii)} below. The uniqueness of the entropy solution follows from the fact that the end points $a_\pm$ are both free. To see this, assume that $\rho$ and $\rho'$ are two solutions that are both concatenations of fundamental solutions. We use
Kruzkov's inequality [K] (as in the proof of Proposition 2.1{\bf(iii)}) to assert
that weakly,
\begin{align*}
|\rho(x,t)-\rho'(x,t)|_t\le &
\big(Q(x,t,\rho(x,t),\rho'(x,t))\big)_x\\
&-\ sgn\big(\rho(x,t)-\rho'(x,t)\big)
\big(H_x(x,t,\rho(x,t))-H_x(x,t,\rho'(x,t))\big)\\
\le &\big(Q(x,t,\rho(x,t),\rho'(x,t))\big)_x+c_1|\rho(x,t)-\rho'(x,t)|,
\end{align*}
where $Q(x,t,\rho,\rho')= sgn(\rho-\rho')
\big(H(x,t,\rho)-H(x,t,\rho')\big).$ Here we have used Hypothesis 1.2{\bf(i)}
for the second inequality.
As a consequence 
\[
\left[e^{-c_1t}|\rho(x,t)-\rho'(x,t)|\right]_t\le e^{-c_1t}
\big(Q(x,t,\rho(x,t),\rho'(x,t))\big)_x.
\]
As in the proof of Proposition 2.1{\bf(iii)}, we can integrate over $[a_-,a_+]$ to assert
\begin{align*}
\left[e^{-c_1t}\int_{a_-}^{a_+}|\rho(x,t)-\rho'(x,t)|\ dx\right]_t\le& e^{-c_1t}
Q(a_+,t,\rho(a_+,t),\rho'(a_+,t))\\
&- e^{-c_1t}
Q(a_-,t,\rho(a_-,t),\rho'(a_-,t)).
\end{align*}
We claim that our free boundary conditions at $a_\pm$ imply that the right-hand side is nonpositive. Indeed, $\mp H_\rho(a_\pm,t,M(a_\pm,t;y_-,y_+))\ge 0$
implies 
\[
Q(a_\pm,t,\rho(a_\pm,t),\rho'(a_\pm,t))=
\mp\big|H(a_\pm,t,\rho(a_\pm,t))-H(a_\pm,t,\rho'(a_\pm,t))\big|.
\]
This allows us to assert
\[
\left[e^{-c_1t}\int_{a_-}^{a_+}|\rho(x,t)-\rho'(x,t)|\ dx\right]_t\le 0.
\]
As an immediate consequence we learn that if $\rho(x,t_0)=\rho'(x,t_0)$
for all $x\in[a_-,a_+]$, then  $\rho(x,t)=\rho'(x,t)$ for all
$(x,t)\in[a_-,a_+]\x[t_0,T]$.
\qed

\ms

\bp\la{pro2.3}  {\bf(i)} If $x_1<x_2$, and  $\xi(\th;x_i,t;y,s)$ is 
a maximizing path in \eqref{eq1.21} for $x=x_i$, then
$\xi(\th;x_1,t;y,s)<\xi(\th;x_2,t;y,s)$  for $\th\in(s,t]$.

\ms\noi
{\bf(ii)}
The fundamental solution $M(x,t;y,s)$ is increasing in $y$.

\ms\noi
{\bf(iii)} Given $s<T$, and $\d\in (0,1)$, there exist positive 
constants $C_0=C_0(s,\d,T)$, and $C_1=C_1(s,\d,T)$ such that
if $|x|\ge C_0$, and $|y|\le (1-\d)|x|$, then $M(x,t;y,s)$ and $-x$ have the same sign, and 
\be\la{eq5.2}
C_1  |x|\le |M(x,t;y,s)|.
\ee
\ep

\ms\noi
{\bf Proof(i)} It is well known that under Hypothesis 2.1{\bf(i)} the following statements are true (see for example [Go]):
\bi
\item [(1)] In \eqref{eq1.21} we may take the supremum over those
$\xi:[t_0,T]$ such that $\xi(s)=y,\ \xi(t)=x$, and $\xi$ is weakly differentiable
 with $\dot\xi\in L^2\big([s,t];\bR\big)$.
\item [(2)] If the supremum in attained at $\xi$, then necessarily $\xi\in C^2$.
\item [(3)] The maximizing path  $\xi$ satisfies the Euler-Lagrange
equation 
\be\la{eq5.3}
\big(L_v(\xi(\th),\th,\dot\xi(\th))\big)_\th=L_x(\xi(\th),\th,\dot\xi(\th)).
\ee
Equivalently, if $p(\th)=L_v(\xi(\th),\th,\dot\xi(\th))$, then the pair
$(\xi,p)$ satisfies the Hamiltonian ODE
\be\la{eq5.4}
\dot \xi(\th)=-H_\rho(\xi(\th),\th,p(\th)),\ \ \ \ \dot p(\th)=H_x(\xi(\th),\th,p(\th)).
\ee
\ei 
We now take $x_1<x_2$, and write $\xi^1,\xi^2:[s,t]\to\bR$ for
the maximizing paths in \eqref{eq1.21} for $x=x_1$ and $x=x_2$ respectively.
We wish to show that $\xi^1(\th)\neq \xi^2(\th)$ for every $\th\in (s,t)$. We argue by contradiction. Suppose to the contrary $\xi^1(\th_0)= \xi^2(\th_0)$ 
for some $\th_0\in (s,t)$. We define
\be\la{eq5.5}
\eta^2(\th)=\begin{cases}\xi^1(\th),\ \ \ \ &\th\in[s,\th_0],\\
\xi^2(\th),\ \ \ \ &\th\in[\th_0,t],\end{cases},\ \ \ \ \ \ 
\eta^1(\th)=\begin{cases}\xi^2(\th),\ \ \ \ &\th\in[s,\th_0],\\
\xi^1(\th),\ \ \ \ &\th\in[\th_0,t].\end{cases}
\ee
Since $\xi^i$ maximizes the action, and $\eta^i$ is weakly
differentiable with square integrable derivative for $i=1,2$,
we learn 
\begin{align}\nonumber
\int_s^tL\big(\eta^1(\th),\th,\dot\eta^1(\th)\big)\ d\th&
\le \int_s^tL\big(\xi^1(\th),\th,\dot\xi^1(\th)\big)\ d\th,\\ 
\int_s^tL\big(\eta^2(\th),\th,\dot\eta^2(\th)\big)\ d\th&
\le \int_s^tL\big(\xi^2(\th),\th,\dot\xi^2(\th)\big)\ d\th.\la{eq5.6}
\end{align}
Expressing the integrals on the left in terms of $\xi^i$ would lead to
\[
 \int_s^{\th_0}L\big(\xi^1(\th),\th,\dot\xi^1(\th)\big)\ d\th
=\int_s^{\th_0}L\big(\xi^2(\th),\th,\dot\xi^2(\th)\big)\ d\th.
\]
This in turn implies that we have equality in \eqref{eq5.6}.
As a result, $\eta^i$ is also a maximizing path with $\eta^i(t)=x_i$.
Hence by (2) above, $\eta^i$ must be $C^1$. This means that we must
have $\dot\xi^1(\th_0)=\dot\xi^2(\th_0)$. Since we also have 
$\xi^1(\th_0)=\xi^2(\th_0)$, we may use the uniqueness of the solutions
to  Euler-Lagrange equation \eqref{eq5.3}, to deduce that $\xi^1=\xi^2$
on $[s,t]$. This contradicts $x_1=\xi^1(t)<x_2=\xi^2(t)$. As a result, $\xi^1$
and $\xi^2$ cannot intersect in $(s,t]$. Since $x_1<x_2$, we must have
$\xi^1(\th)<\xi^2(\th)$ for $\th>s$.

\ms\noi
{\bf(ii)} We  take $y_1<y_2$, and write $\xi^1,\xi^2:[s,t]\to\bR$ for
the maximizing paths in \eqref{eq1.21} for $z=(y_1,s)$ and $z=(y_2,s)$
 respectively. We wish to show that 
$\xi^1(\th)\neq \xi^2(\th)$ for every $\th\in (s,t)$. We again argue by contradiction. Suppose to the contrary $\xi^1(\th_0)= \xi^2(\th_0)$ 
for some $\th_0\in (s,t)$. We define $\eta^1$ and $\eta^2$ as in \eqref{eq5.5}.
Again, since $\xi^1$
(respectively $\xi^2$) is a maximizer  in \eqref{eq1.21}, and that $\eta^2$ 
(respectively $\eta^1$) is weakly
differentiable with square integrable derivative, 
we learn 
\begin{align}\nonumber
\int_s^tL\big(\eta^2(\th),\th,\dot\eta^2(\th)\big)\ d\th&
\le \int_s^tL\big(\xi^1(\th),\th,\dot\xi^1(\th)\big)\ d\th,\\ 
\int_s^tL\big(\eta^1(\th),\th,\dot\eta^1(\th)\big)\ d\th&
\le \int_s^tL\big(\xi^2(\th),\th,\dot\xi^2(\th)\big)\ d\th.\la{eq5.7}
\end{align}
Expressing the integrals on the left in terms of $\xi^i$ would lead to
\[
 \int_{\th_0}^tL\big(\xi^1(\th),\th,\dot\xi^1(\th)\big)\ d\th
=\int_{\th_0}^tL\big(\xi^2(\th),\th,\dot\xi^2(\th)\big)\ d\th.
\]
This in turn implies that we have equality in \eqref{eq5.7}, and that
$\eta^1$ (respectively $\eta^2$)
is also a maximizing path with $\eta^1(s)=y_2$ (respectively $\eta^2(s)=y_1$).
Hence by (2) above, $\eta^i$ must be $C^1$ for $i=1,2$. This means that we must
have $\dot\xi^1(\th_0)=\dot\xi^2(\th)$. By the uniqueness of the corresponding
 Euler-Lagrange equation, we must have $\xi^1=\xi^2$
on $[s,t]$. This contradicts $y_1=\xi^1(s)<y_2=\xi^2(s)$. As a result, $\xi^1$
and $\xi^2$ cannot intersect in $(s,t]$. Since $y_1<y_2$, we must have
$\xi^1(\th)<\xi^2(\th)$ for $\th\le t$. This in particular implies that
$\dot\xi^1(t)\ge \dot\xi^2(t)$. Moreover $\dot\xi^1(t)= \dot\xi^2(t)$
would imply $\dot\xi^1= \dot\xi^2$ by the uniqueness of the corresponding 
\eqref{eq5.3}. As a result, we must have $\dot\xi^1(t)> \dot\xi^2(t)$.
This, \eqref{eq1.22}, and the strict concavity of $L$ in $v$ imply the desired inequality $M(x,t;y_1,s)<M(x,t;y_2,s)$.

\ms\noi
{\bf(iii)} Recall that the pair $(\xi,p)$ satisfies \eqref{eq5.4}, and
the boundary conditions 
\be\la{eq5.8}
\xi(s)=y,\ \ \ \ \xi(t)=x, \ \ \ \ p(t)=M(x,t;y,s).
\ee
From \eqref{eq5.3} and Hypothesis 1.2{\bf(i)}, we learn 
\[
|\dot p(\th)|=\big|\big(L_v(\xi(\th),\th,\dot\xi(\th))\big)_\th\big|
=\big|L_x(\xi(\th),\th,\dot\xi(\th))\big|\le c_1,
\]
which in turn implies
\be\la{eq5.9}
|p(\th)-p(t)|\le c_1(t-s),
\ee
for $\th\in[s,t]$. This,  and \eqref{eq5.4} imply 
\[
|\dot\xi(\th)|=\big|H_\rho(\xi(\th),\th,p(\th))\big|\le 
c_0c_2^{-1}+c_2^{-1}|p(\th)|\le c(1+|p(t)|),
\]
for a constant $c=c(s,T)$. Here we used  
\be\la{eq5.10}
-c_0+c_2\big|H_\rho(x,\th,\rho)\big|\le |\rho|
\le c_0+c_1\big|H_\rho(x,\th,\rho)\big|,
\ee
which follows from Hypothesis 1.2{\bf(i)}. As a result,
\[
|x-y|=|\xi(t)-\xi(s)|\le c'(1+|p(t)|),
\]
for a positive constant $c'=c'(s,T)$. On the other hand, if
 $|y|\le (1-\d)|x|$, then we deduce 
\be\la{eq5.11}
|x|\le c'\d^{-1}(1+|p(t)|).
\ee
We next claim that there exists a constant $C_0$ such that 
\be\la{eq5.12}
|x|\ge C_0\ \ \ \implies  \ \ \ xp(t)<0. 
\ee
To see this, observe that by \eqref{eq5.10} and the monotonicity 
of $\rho\mapsto H_\rho(x,\th,\rho)$,  we can find a constant
$c''$ such that 
\be\la{eq5.13}
|\rho|\ge c''\ \ \ \implies  \ \ \ H_\rho(x,\th,\rho)\rho\ge 0,
\ee
 for all $(x,\th)$. 
Let us assume that $x\ge C_0$,
for a positive constant $C_0$ (to be determined later). 
Suppose contrary to \eqref{eq5.12}, we have $p(t)\ge 0$. 
From \eqref{eq5.11} we deduce 
\[
p(t)\ge{(c')^{-1}}{\d}x-1.
\]
 This and \eqref{eq5.9} imply 
\be\la{eq5.14}
p(\th)\ge {(c')^{-1}}{\d}x-1-c_1(t-s)\ge {(c')^{-1}}{\d}C_0-1-c_1(T-s) ,
\ee
for all $\th\in (s,t)$. Choose $C_0$ large enough so that the right-hand side of
\eqref{eq5.14} is at least $c''$. This would guarantee  
\[
H_\rho(\xi(\th),\th,p(\th))\ge 0,\ \ \ {\text{ for }}\ \th\in [s,t],
\]
 by \eqref{eq5.13}. From this  and \eqref{eq5.4}
we deduce that $\dot\xi(\th)\le 0 $ for $\th\in[s,t]$. 
As a result, $x-y=\xi(t)-\xi(s)\le 0$. But this is impossible
if $y\le (1-\d)x$. Hence  the condition
$x\ge C_0$ implies that $p(t)> 0$. In the same fashion, we can show that
the condition
$x\le -C_0$ implies that $p(t)< 0$. This completes the proof of \eqref{eq5.12}.
From this, \eqref{eq5.8}, and \eqref{eq5.11},  we can readily deduce
\eqref{eq5.2}.
\qed

\ms\noi
We next give a recipe for the law of the process $\by_t$.

\ms\noi
{\bf Definition 5.2(i)}  We  set
\[
\G(a,b,t,\rho)=\int_{a}^b\hat A(g)(z,t,y) \ dz,\ \ \ \  \G(\bq,t)
=\sum_{i=0}^{n}\G(x_{i},x_{i+1},t,y_i).
 \]
 
 \ms\noi
 {\bf(ii)} We define a measure $\mu(d\bq,t)$ on the set $\D$ that is our candidate for the law of $\bq(t)$. 
 The restriction of
$\mu$ to $\D_n$ is given by 
\[
\mu^n(d\bq,t):=
\ell(t, dy_0)\exp\left\{-\G(\bq,t)\right\}
\prod_{i=1}^{n}
\ g\big(x_i,t,y_{i-1},{y_{i}})\  dx_{i}dy_{i},
\]
where $g$ solves \eqref{eq1.25} and $\ell$ solves \eqref{eq1.26}.
To simplify our presentation, we assume that $\ell(t,dy_0)=\ell(t,y_0)\ dy_0$ is 
absolutely continuous with respect to the Lebesgue measure. Such an assumption would allow us to express $\mu^n(d\bq,t):=\mu^n(\bq,t)\ d\bq$, where
\[
d\bq=dy_0\ \prod_{i=1}^{n}dx_{i}dy_{i},\ \ \ \
\mu^n(\bq,t)=\ell(t, y_0)\exp\left\{-\G(\bq,t)\right\}
\prod_{i=1}^{n}
\ g\big(x_i,t,y_{i-1},{y_{i}}).
\]
\qed 

\ms\noi
\bp\la{pro5.3} Let $g$ be a solution of \eqref{eq1.25}. Then
$\hat A(g^1)_t=\hat A(g^2)_x$.
\ep

\ms\noi
{\bf Proof} From integrating both sides of \eqref{eq1.25} with respect $y_+$ 
we learn
\be\la{eq5.15}
\hat A(g^1)_t-\hat A(g^2)_x=\hat A\big(\hat Q^+(g)\big)-
\hat A\big(g\hat J(g)\big).
\ee
On the other hand,
\begin{align*}
\hat A\big(\hat Q^+(g)\big)(y_-)=&\int g^1(y_-,y_*) \hat A(g^2)(y_*)\ dy_*-\int g^2(y_-,y_*) \hat A(g^1)(y_*)\ dy_*,\\
\hat A\big(g\hat J(g)\big)(y_-)=&\int g^1(y_-,y_*) \hat A(g^2)(y_*)\ dy_*-\hat A(g^2)(y_-)\hat A(g^1)(y_-)\\
&-\int g^2(y_-,y_*) \hat A(g^1)(y_*)\ dy_*+\hat A(g^1)(y_-)
\hat A(g^2)(y_-).
\end{align*}
This implies that the right-hand side of \eqref{eq5.15} is $0$.
\qed

\ms
We are now ready to present the proof of Theorem 5.1, which is similar to the proof of Theorem 2.1. 

\bs\noi
{\bf Proof of Theorem 5.1} We wish to establish the analog of
\eqref{eq4.1} in our setting.  Theorem 3.1, and a repetition of
the first step of the proof of Theorem 2.1 allow us to reduce the proof of 
Theorem 5.1 to the verification of 
an analog of \eqref{eq4.4}, namely
\begin{align}\nonumber
\int_{\D_n}\hat G^n(\bq,s)\ \mu^n_s(d\bq,s)=&
\int_{\D_n}  \hat G^n(\bq,s)\big(\widehat\cL_{0n}^{s*}
\mu^n\big)(\bq,s)\ d\bq\\
&+\int_{\hat\p \D_{n+1}}  \hat G^{n+1}(\bq,s)\mu^{n+1}(\bq,s)(\bb_{n+1}\cdot {\bold{N}}_{n+1})\ \s(d\bq),\la{eq5.16}
\end{align}
with $\hat G$ as in \eqref{eq3.5}. For a more tractable expression
for  the left hand side of \eqref{eq5.16}, we write
\be\la{eq5.17}
\mu_s^n=X^n\mu^n,
\ee
where
\be\la{eq5.18}
X^n=-\G_s(\bq,s)+\frac 
{ \ell_s(s,y_0)}{\ell(s,y_0)}
+\sum_{i=1}^{n}\frac{g_s(x_i,s,y_{i-1},y_{i})}
{g(x_i,s,y_{i-1},y_{i})}.
\ee
On the other hand,  by Proposition 5.3,
\begin{align*}
\G_s(\bq,s)=&\sum_{i=0}^n\int_{x_i}^{x_{i+1}}\hat A(g^1)_s(z,s,y_i)\ dz
=\sum_{i=0}^n\int_{x_i}^{x_{i+1}}\hat A(g^2)_z(z,s,y_i)\ dz\\
=&\sum_{i=0}^n
\big(\hat A(g^2)(x_{i+1},s,y_i)-\hat A(g^2)(x_{i},s,y_i)\big)\\
=&\hat A(g^2)(x_{n+1},s,y_n)-\hat A(g^2)(x_{0},s,y_0)-\sum_{i=1}^n\big(\hat A(g^2)(x_{i},s,y_i)
-\hat A(g^2)(x_{i},s,y_{i-1})\big).
\end{align*}
From this, \eqref{eq5.18}, \eqref{eq1.26}, and \eqref{eq1.25} 
we deduce
\begin{align}\nonumber
X^n=&-\hat A(g^2)(a_+,s,y_n)+\frac 
{ (\ell*g^2)(a_-,s,y_0)}{\ell(s,y_0)}+\sum_{i=1}^{n}\frac{\hat Q^+(g)(x_i,s,y_{i-1},y_{i})}
{g(x_i,s,y_{i-1},y_{i})}\\
&+\sum_{i=1}^{n}\hat v(x_i,s,y_{i-1},y_i)\big(\hat A(g)(x_{i},s,y_i)
-\hat A(g)(x_{i},s,y_{i-1})\big).\la{eq5.19}
\end{align}

We can rewrite the right-hand side of \eqref{eq5.16} as
\[
\int_{\D_n}  \hat G^n(\bq,s) Y^n(\bq)
\mu^n(\bq,s)\ d\bq,
\]
where $Y^n=Y^n_1+Y^n_2$, with
$Y^n_1$ and $Y^n_2$ corresponding to two terms on the right-hand side 
of \eqref{eq5.16}. Indeed, an integration by parts yields
\begin{align}\nonumber
Y^1_n(\bq)=&-\sum_{i=1}^n\hat v_i(x_i,s,y_{i-1},y_i)\G(s,\bq)_{x_i}\\
=&\sum_{i=1}^n\hat v_i(x_i,s,y_{i-1},y_i)\big(\hat A(g)(x_{i},s,y_{i})-
\hat A(g)(x_{i},s,y_{i-1})\big).\la{eq5.20}
\end{align}
 As for $Y^n_2$, we write $Y^n_2=Y^n_{2-}+Y^n_{2*}+Y^n_{2+}$,
where the terms $Y^n_{2-},\ Y^n_{2*}$, and $ Y^n_{2+}$ correspond
to the boundary contributions associated with the conditions
$x_1=a_-$, $x_i=x_{i+1},$ with $i\in\{1,\dots,n\}$, and $x_{n+1}=a_+$,
respectively. More precisely,  
\begin{align*}
Y^n_{2-}(\bq)=&\frac { (\ell*g^2)(a_-,s,y_0)}{\ell(s,y_0)},\ \ \ \ \ 
Y^n_{2+}(\bq)=-\hat A(g^2)(a_+,s,y_n),\\
Y^n_{2*}(\bq)=&\sum_{i=1}^{n}\frac{\hat Q^+(g)(x_i,s,y_{i-1},y_{i})}
{g(x_i,s,y_{i-1},y_{i})}.
\end{align*}
This, \eqref{eq5.17},  \eqref{eq5.19}, and \eqref{eq5.20} complete the proof
of \eqref{eq5.16}.
\qed

\section{Proofs of Proposition 1.1 and Theorem 1.3}

\noi
{\bf Proof of Proposition 1.1} 
Let us write 
\[
\cK(g)= \nabla\cdot (\nu g)-\hat Q^+( g)+\hat Q^-( g).
\]
where $\nu=(-\hat v,1)$, and $\nabla=(\p_x,\p_t)$. 
To ease the notation, we do not display the dependence of
$g,\ \hat g, \eta$, and $h$ on $(x,t)$.
We certainly have 
\begin{align*}
&\frac {\hat  Q^+(\hat g)}{\hat g}-
\frac {\hat Q^+( g)}{ g}=0,\\
&\left(\frac{\nabla\cdot(\nu\hat g)}{\hat g}-
\frac{\nabla\cdot(\nu g)}{ g}\right)(y_-,y_+)=
\left(\nu\cdot\frac{\nabla \eta}{ \eta}\right)(y_-,y_+)
=\nu\cdot\frac{\nabla h}{ h}(y_+)-
\nu\cdot\frac{\nabla h}{ h}(y_-),\\
&\frac {\hat Q^-( g)}{ g}(y_-,y_+)=
 {\hat A(\hat vg)}(y_+)-\hat A(\hat vg)(y_-)
-\hat v(y_-,y_+)\left( {\hat A(g)}(y_+)
-\hat A(g)(y_-)\right)\\
&\frac {\hat Q^-(\hat g)}{\hat g}(y_-,y_+)=
\frac {\hat A(\hat vg\otimes h)}{h}(y_+)
-\frac { \hat A(\hat vg\otimes h)}{h}(y_-)\\
&\ \ \ \ \ \ \ \ \ \ \ \ \ \ \ \ \ \ \ \ \ \ 
-\hat v(y_-,y_+)\left(\frac {\hat A(g\otimes h)}{h}(y_+)
-\frac {\hat A(g\otimes h)}{h}(y_-)\right).
\end{align*}
Hence
\begin{align*}
\left(\frac{ \cK(\hat g)}{\hat g}-\frac{ \cK( g)}{ g}\right)(y_-,y_+)=&
\nu(y_-,y_+)\cdot\left(\frac{\nabla h}{ h}(y_+)-\frac{\nabla h}{ h}(y_-)\right)
+\left(\frac {\hat Q^-(\hat g)}{\hat g}-\frac {\hat Q^-( g)}{ g}\right)(y_-,y_+)\\
=&
\frac{h_{t}+\cB^2h}h(y_+)-\frac{h_{t}+\cB^2h}h(y_-)\\
&-\hat v(y_-,y_+)\left[\frac{h_{x}+\cB^1h}h(y_+)
-\frac{h_{x}+\cB^1h}h(y_-)\right].
\end{align*}
The right-hand side is $0$, when $h$ satisfies \eqref{eq1.31}.
This completes the proof because $g$ (respectively $\hat g$)
solves \eqref{eq1.25}
if and only if $\cK(g)=0$ (respectively $\cK(\hat g)=0$).
\qed

\bs
The proof of Theorem 1.3, uses Doob's $h$-transform that we now recall.

\ms\noi
\bp\la{pro6.1} Let $\bP$ be the law of Markov jump process
$\big(\by(x):\ x\in [a_-,a_+]\big)$, with the jump kernel density $g(x,y_-,y_+)$,  and the generator $\cL_x$. Assume that $g$ is $C^1$ in $x$.
Let $U$ be an interval, and let
$\widehat \bP$ denote the law of $\bP$, conditioned on the event $\by(x)\in U$
for all $x\in[a_-,a_+]$. Then $\widehat \bP$ is the law of a Markov jump process
with a jump kernel density $\hat g$, given by
\be\la{eq6.1}
\hat g(x,y_-,y_+)=\frac{h(x,y_+)}{h(x,y_-)}\  g(x,y_-,y_+),\ \  \ \ \ x\in [a_-,a_+),
\ \ y_\pm\in U,
\ee
where 
\be\la{eq6.2}
h(x,y)=\bP\big(\by(a)\in U\ \  {\text{ for }}\ a\in [x,a_+]\ |\ \by(x)=y\big).
\ee
Moreover, $h$ is $C^1$ in $x$, and satisfies
\be\la{eq6.3}
h_x+\cL_x h=0.
\ee
\ep

\ms\noi
{\bf Proof} {\em(Step 1)} We can write
\be\la{eq6.4}
h(x,y)=\sum_{n=0}^\i h_n(x,y),
\ee
where
\be\la{eq6.5}
h_n(x,y)=\int_{X_n(x,y)}\mu^n(\bq,x,y)\ d\bq_n\\
\ee
where for $n\ge 1$, 
\begin{align*}
&\bq_n=(x_1,y_1,\dots,x_{n},y_{n}),\ \ \ \ \ 
d\bq_n= \prod_{i=1}^{n}dx_{i}dy_{i},\\
&\mu^n(\bq_n,x,y)=\exp\left\{-\G(\bq_n,x,y)\right\}
\prod_{i=1}^{n}\ g\big(x_i,y_{i-1},{y_{i}}),\ \ {\text{ with }}\ \  y_0=y,\\
&\G(\bq_n,x,y)=
\sum_{i=0}^n\G(x_i,x_{i+1},y_i),\ \ \ \ \ {\text{with }}\ \ \ \ x_0=x,\ y_0=y,
\\
&\G(a,b,y)=\int_a^b (\hat A g)(z,y)\ dz,
\end{align*}
and the set $X_n(x,y)$ consists of  $\bq$, satisfying
\[
x< x_1<\dots<x_n<a_+,\ \ \ \ 
 y_1,\dots ,y_{n}\in U.
\]
 When $n=0$, we simply have $h_0(x,y)=\exp\big\{-\G(x,a_+,y)\big\}.$
It is straightforward to verify continuous differentiability of $h$,
and deduce \eqref{eq6.3} from \eqref{eq6.4} and 
\eqref{eq6.5}.

\ms\noi
{\em(Step 2)}
The law  $\hat P$ is simply given by 
\[
\hat P=\sum_{n=1}^\i\hat \mu^n,
\]
where $\hat \mu_n(d\bq_n)=\hat\mu_n(\bq_n)\ d\bq_n,$ with 
\be\la{eq6.6}
\hat \mu_n(\bq_n)=h(a_-,y)^{-1}\  \mu_n(\bq_n)\ 1\!\!1\big(y_1,\dots,y_n\in U\big).
\ee
We wish to show that $\hat P$ is the law of a jump process associated with the jump density $\hat g$. To achieve this, we rewrite $\hat \mu_n$ using the fact that $h$ satisfies \eqref{eq6.3}. Indeed, \eqref{eq6.3} implies
\be\la{eq6.7}
e^{-\int_a^b \hat A( g)(z,y_-)dz}\ \frac {h(b,y_+)}{h(a,y_-)}
=e^{-\int_a^b \hat A(\hat g)(z,y_-) dz}\ \frac {h(b,y_+)}{h(b,y_-)}.
\ee
This is equivalent to asserting
\begin{align}\la{eq6.8}
 \frac{h(b,y_-)}{h(a,y_-)}=&\exp\left(-\int_a^b \big(\hat A(\hat g)-\hat A( g)\big)
(z,y_-)\ dz\right)\\ \nonumber
=&\exp\left(-\int_a^b \left(\frac{\hat A( g\otimes h)-\hat A( g) h}{h}\right)
(z,y_-)\ dz\right)\\ \nonumber
=&\exp\left(\int_a^b \frac{-(\cL_z h)(z,y_-)}{h(z,y_-)}\ dz\right)=
\exp\left(\int_a^b \frac{h_z(z,y_-)}{h(z,y_-)}\ dz\right)\\
=&\exp\left(\int_a^b (\log h)_z(z,y_-)\ dz\right),\nonumber
\end{align}
which is evidently true. We set, $x_0=a_-,\ y_0=y$ as before.
Observe that $\hat \mu_n(\bq_n)$ of \eqref{eq6.6}
can be written as
\begin{align*}
&\frac 1{h(a_+,y_n)}\ e^{-\int_{x_{n}}^{a_+} \hat A(g)(z,y_{n}) dz}\ 
\frac {h(a_+,y_n)}{h(x_n,y_n)}\ 
\prod_{i=1}^{n}
 e^{-\int_{x_{i-1}}^{x_i} \hat A(g)(z,y_{i-1}) dz}\ \frac
{h(x_i,y_i)}{h(x_{i-1},y_{i-1})} g(x_i,y_{i-1},y_i)\\
&\ \ \ \ \ \ \  =\frac 1{h(a_+,y_n)}\ e^{-\int_{x_{n}}^{a_+} \hat A(\hat g)(z,y_{n}) dz}\prod_{i=1}^{n}
 e^{-\int_{x_{i-1}}^{x_i} \hat A(\hat g)(z,y_{i-1}) dz}\ \frac
{h(x_i,y_i)}{h(x_{i},y_{i-1})} g(x_i,y_{i-1},y_i)\\
&\ \ \ \ \ \ \  
=e^{-\int_{x_{n}}^{a_+} \hat A(\hat g)(z,y_{n}) dz}\prod_{i=1}^{n}
 e^{-\int_{x_{i-1}}^{x_i} \hat A(\hat g)(z,y_{i-1}) dz}\ 
 \hat g(x_i,y_{i-1},y_i),
\end{align*}
where we used \eqref{eq6.8} and \eqref{eq6.7} for the first equality,
and for the last equality we used the definition of $\hat g$, and
$h( a_+,y_n)=1$, which follows from the definition of $h$.
The right-hand side is the law of a Markov jump process associated
with the kernel density $\hat g$, as desired.
\qed

\bs\noi
{\bf Proof of Theorem 1.3} {\em(Step 1)} Recall that $\rho(x,t)$ is the solution 
of \eqref{eq1.2} with the initial condition $\rho(x,t_0)=
\rho(x,t_0;\by_{t_0},s)$, where $\by_{t_0}$ is a jump process associated with the kernel $g(x,t_0,y_-,y_+)$. We wish to show that $\rho(x,t)=
\rho(x,t;\by_{t},s)$ for $(x,t)\in[a_-,\i)\x[t_0,T]$. It suffices to verify this for
$(x,t)\in[a_-,a_+]\x[t_0,T]$, where $a_+$ is any large number in $(a_-,\i)$.

Pick any $\d\in (0,1)$. From Proposition 5.2{\bf(iii)}, and \eqref{eq5.13},
we learn that there exist
 constants $C_0=C_0(\d),$ $C_1=C_1(\d)$, and $C_2$ such that 
\begin{align}\la{eq6.9}
|y_{+}|\le (1-\d)a_+,\ \  \  a_+\ge C_0\ \  &\implies\ \  
 M(a_+,t_0;y_+,s)<-C_1a_+,\\
\ (x,\th)\in \bR\x[s,T], \ \ \ \  |\rho|\ge C_2\ \   &\implies \ \ \ 
\rho H_\rho(x,t,\rho)>0,\la{eq6.10}
\end{align}
for every $t\in[t_0,T]$.  Note that \eqref{eq6.9} and Proposition 5.2{\bf(ii)}
imply that if $a_+\ge C_0$, then 
\[
y_-<y_+,\ \ |y_+|< (1-\d)a_+,\ \  \    \implies\ \  \
 M(a_+,t_0;y_-,s)<M(a_+,t_0;y_+,s)<-C_1a_+.
\]
From this, and \eqref{eq6.2} we can readily deduce 
\be\la{eq6.11}
Y_-\le y_-<y_+\le (1-\d)a_+,\ \  \ a_+\ge C_3\ \ \ \implies\ \ \ 
\hat v(x,t,y_-,y_+)>0,
\ee
for every $t\in[t_0,T]$, where 
\[
C_3=\max\big\{C_0,C_1^{-1}C_2,|Y_-|\big\}.
\] 
We pick any $ a_+\ge \max\{a_-,C_3\}$.

\ms\noi
{\em(Step 2)} 
We write $\bW$  for the law of the Markov process $(\bw(t):\ t\in[t_0,T])$, 
associated with the generator $\cB^2_{a_-,t}$, such that $\bw(t_0)=y^0$. 
We also define a family of probability measures
$\big(\bP_{t}: t\in [t_0,T]\big)$ with the following recipe:
For each $t$,  $\bP_{t}$ is the law of the Markov process 
$\by_t:[a_-,a_+]\to[Y_-,\i),$ associated
with the generator $\cB^1_{x,t}$, satisfying the initial condition 
$\by_t(a_-)=\bw(t)$. 
We define $\th_\d$ to be the smallest
$t>t_0$ such that $\bw(t)\notin U(\d)$, where
\[
U(\d):=[Y_-, (1-\d)a_+)=:[Y_-,Y^\d).
\]
We also define $\tau_\d(t)$ to be the smallest
$x>a_-$ such that $\by_t(x)\notin U(\d)$.
We set
\[
\bW^{\d}(A)=\bW(A\ |\  \th_\d>T),
\ \ \ \ \bP_{t}^{\d}(A)=\bP_t(A\ |\  \tau_\d(t)>a_+).
\]
We write $\bw^\d(t),\ t\in [ t_0,T]$
for the  jump process that is distributed according to $\bW^\d$.
For each $t\in[t_0,T]$, we  write $\by_t^\d$ for the jump process that is distributed according to $\bP_{t}^{\d}$.
By Proposition 6.1, the process $t\to \bw^\d(t)$, and the processes
$x\mapsto\by_t^\d(x),\ t\in [t_0,T]$ are again Markov jump processes. 
We set
\begin{align*}
h^\d(x,t,y)&=
\bP^{\d}_t\big( \tau_\d(t)>a_+\ | \ \by_t(x)=y\big),\\
\ell^\d(t,y)&=\bW\big(  \th_\d>T\ |\ \bw(t)=y\big).
\end{align*}
By \eqref{eq6.3}, we have the following equations for $h^\d$ and $\ell^\d$:
\begin{align}\la{eq6.12}
&h^\d_{x}(x,t,y)+
\big(\cB^1_{x,t}h^\d\big)(x,t,y)=0,\ \ \ \ & x\in (a_-,a_+), 
\ y\in U(\d), \ t\in [t_0,T],\\
&\ell^\d_{t}(t,y)+
\big(\cB^2_{a_-,t}\ell^\d\big)(s,y)=0,\ \ \ \ & t\in (t_0,T), \ y\in U(\d).\la{eq6.13}
\end{align}

Since, $h^\d(a_-,t,y)=\ell^\d(t,y)$, the equations \eqref{eq6.12} and \eqref{eq6.13} 
allow us to apply Proposition 1.2 to assert that  $h^\d$ satisfies
\be\la{eq6.14}
h^\d_{x}(x,t,y)+
\big(\cB^1_{x,t}h^\d\big)(x,t,y)=0,\ \ \ \ h^\d_{t}(x,t,y)+
\big(\cB^2_{x,t}h^\d\big)(x,t,y)=0.
\ee

\ms\noi
{\em(Step 3)} Since the  jump process $\by^\d_{t_0}$ takes value in
$[Y_-,(1-\d)a_+)$, our Theorem 1.2 or Theorem 5.1 is applicable. 
More precisely, if the initial data of \eqref{eq1.2} is given by
$\rho(x,t_0;\by^\d_{t_0},s)$, for some $s<t_0$, and 
with $y^\d_{t_0}$ a Markov process distributed according
to $\bP_{t_0}^{a_+,\d}$, then the solution $\rho(x,t)$ at a later  time $t>t_0$
is given by $\bP_{t}^{\d}$.
As a consequence, Theorem 1.3
 is true when the initial process satisfies $\by_{t_0}(a_+)\in U(\d)$.
This condition is true with probability density $h^\d(a_-,t_0,y^0)$.
The condition \eqref{eq1.29} would implies that $h^\d(a_-,t_0,y^0)\to 1$, 
and $\by^\d_t\to\by_t$ in small $\d$ limit, when restricted to the interval
 $[a_-,a_+]$. This completes the proof.
\qed

\bs\noi
%{\bf Acknowledgments} The authors thank Govind Menon for many fruitful discussions and his very helpful comments on the first draft of this article.


\begin{thebibliography}{HR}

\bibitem[AE]{AE} J. Abramson and
S. N. Evans, Lipschitz minorants of Brownian motion and
 L\'evy processes,   Probab. Theory Relat. Fields {\bf 158}, 809--857 (2013). 


\bibitem[B]{B} Y. Bakhtin, Burgers equation with Poisson random forcing,  Ann. of Probab. {\bf 1}, 2961--2989 (2013).

\bibitem[BCK]{BCK} Y. Bakhtin, E. Cator and K. Khanin,
Space-time stationary solutions for the Burgers equation,
J. Amer. Math. Soc. {\bf 27}, 193--238 (2014).

\bibitem[Be]{Be} J. Bertoin, The Inviscid Burgers Equation with Brownian initial velocity,
Commun. Math. Phys. {\bf 193}, 397-406 (1998).


\bibitem[CD1]{CD1} L. Carraro and J. Duchon, Solutions statistiques intrinseques de
l'\'equation de Burgers et processus de L\'evy, Comptes Rendus de l'Acad\'emie
des Sciences. S\'erie I. Math\'matique {\bf 319},  855-858 (1994).



\bibitem[CD2]{CD2} L. Carraro and J. Duchon, Equation de Burgers avec conditions initiales
a accroissements ind\'ependants et homogenes, Annales de l'Institut
Henri Poincare (C) Non Linear Analysis {\bf 15}, 431-458 (1998).


\bibitem[D]{D} M. H. A. Davis,
Piecewise-Deterministic Markov Processes: A General Class of Non-Diffusion Stochastic Models, 
J. Royal Stat. Soc. Series B (Methodological) 
{\bf 46}, 353-388 (1984).  


\bibitem[EO]{EO} S. N. Evans
and M. Ouaki, Excursions away from the Lipschitz minorant of a L\'evy process, Ann. Inst. H. Poincar\'e Probab. Statist. {\bf58},
424--454 (2022).


\bibitem[Go]{Go}  D. Gomes, Viscosity Solutions of
Hamilton-Jacobi Equations, Coloquio Brasileiro de Matematica {\bf 27} (2009).




\bibitem[Gr]{Gr} P. Groeneboom, Brownian motion with a parabolic drift and airy functions,
 Probab. Theory and Relat. Fields {\bf 81}, 79-109 (1989).




 
 

\bibitem[KR1] {KR1} D. Kaspar and F. Rezakhanlou,
Scalar conservation laws with monotone pure-jump Markov initial conditions ,
  Probab. Theory Related Fields {\bf 165}, 867-899 (2016).

\bibitem[KR2] {KR2} D. Kaspar and F. Rezakhanlou,
Kinetic statistics of scalar conservation laws with piecewise-deterministic Markov process data,  Arch. Rational Mech. Anal. {\bf 1}, 259-298 (2020).

\bibitem[K] {K} S.N. Kruzhkov,
First order quasilinear equations in several independent variables,
Mat. Sb. {\bf 81}   217--243 (1970).


\bibitem[L] {L} L-C Li,  An Exact Discretization of a Lax Equation for Shock Clustering and Burgers Turbulence I: Dynamical Aspects and Exact Solvability, Commun. Math. Phys.
{\bf 361}, 415-466 (2018).





\bibitem[M2]{M2} G. Menon, Complete integrability of shock clustering and Burgers turbulence, Arch. Ration. Mech. Anal.
 {\bf 203},  853-882 (2012).

 


\bibitem[MS]{MS} G. Menon and R. Srinivasan, Kinetic theory and Lax equations for
shock clustering and Burgers turbulence, J. Statist. Phys. {\bf140}, 1-29 (2010).

\bibitem[O]{O} M. Ouaki, Scalar conservation laws with white noise initial data . Probab. Theory Relat. Fields {\bf 182}, 955--998 (2022).



\bibitem[OR1]{OR1} M. Ouaki and F. Rezakahnlou,
Random Tessellations and Gibbsian Solutions of Hamilton-Jacobi Equations, Commun. Math. Phys. (2022).

\bibitem[OR2]{OR2} M. Ouaki and F. Rezakahnlou,
A kinetic approach to Burgers equation with white noise initial data,
preprint.

\bibitem[PS]{PS} M. Pr\"ahofer and H. Spohn, Scale invariance of the PNG droplet and the Airy process, Journal of Statistical Physics {\bf 108}, 1071-1106 (2002).



\bibitem[PW]{PW} J. Pitman and W. Tang,
The argmin process of random walks, Brownian motion and L\'evy processes, Electro. J. Probab  {\bf 23}, 1--35 (2018).

\bibitem[R]{R} F. Rezakhanlou, Stochastic Solutions to Hamilton-Jacobi Equations,  Springer Proceedings in Mathematics and
Statistics {\bf 282}, 206-238 (2019).






\end{thebibliography}
\end{document}